\documentclass[a4paper,12pt,twoside]{article}
\usepackage[latin1]{inputenc}                 
\usepackage{amsfonts,amssymb,amsmath,exscale} 
\usepackage[dvips]{graphicx,color,psfrag}     
\usepackage{fancyhdr}                         
\usepackage{caption}                          
\usepackage{rotating}
\usepackage{defcml}                           
\usepackage[numbers]{natbib}                  
\Floatboxname{Box}
\usepackage{verbatim}                         

\usepackage{amsmath,amsthm,amsfonts,amssymb,amscd}
\usepackage{subfigure}
\usepackage{graphicx,bm,color}
\usepackage{algorithm}
\usepackage{algpseudocode}
\usepackage{stmaryrd}
\usepackage{booktabs,ctable,multirow,longtable} 
\usepackage{url}
\usepackage[latin1]{inputenc}                 
\usepackage{amsfonts,amssymb,amsmath,exscale} 
\usepackage[dvips]{graphicx,color,psfrag}     
\usepackage{fancyhdr}                         
\usepackage{caption}                          
\usepackage{rotating}
\usepackage{defcml}                            
\usepackage[numbers]{natbib}                  
\Floatboxname{Box}
\usepackage[framed,numbered,autolinebreaks,useliterate]{mcode}
\definecolor{gray}{gray}{0.6}

\newtheorem{Remark}{Remark}[section]

\newtheorem{form}{Formulation}[section]

\def\R{\mathbb R}

\renewcommand{\tablename}{Table}
\fancypagestyle{plain}{%
	\fancyhead{}%
	\fancyfoot[c]{\sffamily\thepage}%
}
\makeatletter
\def\cleardoublepage{\clearpage\if@twoside \ifodd\c@page\else
	\hbox{}
	\vspace*{\fill}
	\thispagestyle{empty}
	\newpage
	\if@twocolumn\hbox{}\newpage\fi\fi\fi}

\usepackage{scalerel,stackengine}
\stackMath
\newcommand\reallywidecheck[1]{%
	\savestack{\tmpbox}{\stretchto{%
			\scaleto{%
				\scalerel*[\widthof{\ensuremath{#1}}]{\kern-.6pt\bigwedge\kern-.6pt}%
				{\rule[-\textheight/2]{1ex}{\textheight}}
			}{\textheight}%
		}{0.5ex}}%
	\stackon[1pt]{#1}{\scalebox{-1}{\tmpbox}}%
}

\begin{document}
\Titel{
Adaptive Global-Local Approach for Phase-Field Ductile Fracture
      }
\Autor{F. Aldakheel, N. Noii, T. Wick, P. Wriggers}
\Report{02--I--17}
\Journal{

}
%



\thispagestyle{empty}

\ce{\bf\large Multilevel Global-Local techniques for adaptive ductile phase-field fracture}

\vskip .35in

\ce{Fadi Aldakheel\(^{a}\), Nima Noii\(^{a,b,}\)\footnote{Corresponding author.\\[1mm]
		E-mail addresses:
		aldakheel@ikm.uni-hannover.de (F. Aldakheel); noii@ikm.uni-hannover.de (N. Noii); thomas.wick@ifam.uni-hannover.de (T. Wick); allix@lmt.ens-cachan.fr (O. Allix);  wriggers@ikm.uni-hannover.de (P. Wriggers).
	}, Thomas Wick\(^{b,d}\), Olivier Allix\(^{c}\), Peter Wriggers\(^{a,d}\)} \vskip .25in

\ce{\(^a\) Institute of Continuum Mechanics} \ce{Leibniz Universit\"at Hannover, An der Universit\"at 1, 30823 Garbsen, Germany}\vskip .25in

\ce{\(^b\) Institute of Applied Mathematics} \ce{Leibniz Universit\"at Hannover, Welfengarten 1, 30167 Hannover, Germany} \vskip .22in

\ce{\(^c\) LMT, ENS Paris-Saclay/CNRS/Universit\'{e} Paris-Saclay} \ce{61 avenue du Pr\'{e}sident Wilson, F-94235 Cachan Cedex,  France}\vskip .22in

\ce{\(^d\) Cluster of Excellence PhoenixD (Photonics, Optics, and
	Engineering - Innovation} \ce{Across Disciplines), Leibniz Universit\"at Hannover, Germany}\vskip .25in

\begin{Abstract}
    This paper outlines a rigorous variational-based \textit{multilevel Global-Local} formulation for  ductile fracture. Here, a phase-field formulation is used to resolve failure mechanisms by regularizing the sharp crack topology on the local state. The coupling of plasticity to the crack phase-field is realized by a constitutive work density function, which is characterized through a degraded stored elastic energy and the accumulated dissipated energy due to plasticity and damage. Two different Global-Local approaches based on the idea of  multiplicative Schwarz' alternating method are proposed: (i) A global constitutive model  with an elastic-plastic behavior is first proposed, while it is enhanced with a single local domain, which, in turn, describes an elastic-plastic fracturing response. (ii) The main objective of the second model is to introduce an adoption of the Global-Local approach toward the multilevel local setting. To this end, an elastic-plastic global constitutive model is augmented with {two distinct} local domains; in which, the first local domain behaves as an elastic-plastic material and the next local domain is modeled due to the fracture state. To further reduce the computational cost, predictor-corrector adaptivity within Global-Local concept is introduced. An adaptive scheme is devised through the evolution of the effective global plastic flow (for \textit{only elastic-plastic} adaptivity), and through the evolution of the local crack phase-field state (for \textit{only fracture} adaptivity). Thus, two local domains are dynamically updated during the computation, resulting with \textit{two-way} adaptivity procedure. The overall response of the Global-Local approach in terms of accuracy/robustness and efficiency is verified using single-scale problems. The resulting framework is algorithmically described in detail and substantiated with numerical examples.
	
\end{Abstract} 

\tableofcontents \cleardoublepage

\sectpa[Section1]{Introduction}

Variational phase-field modeling is the regularized fractured formulation with a strong capability to simulate complicated failure processes. This includes crack initiation (also in the absence of a crack tip singularity), propagation, coalescence, and branching without additional ad-hoc criteria \cite{BourFraMar08,ambati15,wu2018phase}. Such a feature is particularly attractive for industrial applications, as it minimizes the need for time-consuming and expensive calibration tests \cite{NoiiGL18}. In contrast to these advantages, the finite element treatment of the phase-field formulation is known to be computationally demanding, mainly due to the non-convexity of the energy functional to be minimized with respect to the displacement and the phase-field \cite{Gerasimov16,Wick15Adapt,wick2017modified}. Other challenges for the phase-field fracturing formulation is two-fold. 
\begin{itemize}
	\item First, that is a regularized-based formulation which is strongly linked to the element discretization size $h$ due to the principal parameters: a  small residual scalar $\kappa$ and characteristic length-scale $l$ \cite{Wi20_book}.  Specifically, $\kappa:=\kappa(h)$ and $l:=l(h)$ hold such that $h\ll l$ and $h\ll\kappa$ through discretization error estimates \cite{MaWi19,Wi20_book}. Hence, the equations to be
	minimized for the variational phase-field formulation are strongly related to the element size $h$. Thus, for resolving the crack phase-field, a sufficiently small $h$ is chosen to obtain the experimental resolution \cite{BourFraMar00,BourFraMar08,schroder2020selection,Wick15Adapt}. 
	\item The second challenge is to use the phase-field fracture approach for structures of industrial complexity \cite{NoiiGL18,noii2020adaptive,aldakheel2020global}. This has been the subject of limited investigations, and further studies in this direction will pave the way for the wide adoption of phase-field modeling within legacy codes for industrial applications. 	
\end{itemize}
In fact, when dealing with large structures, the failure behavior is solely analyzed in a (small) local region, whereas  in the surrounding medium, a simplified and linearized system of equations can be solved. Thus,  the  idea  of  a  two-scale  formulation,  in  which the nonlinear displacement/phase-field problem is solved on a lower(local) scale while dealing with a purely linear/homogeneous problem on an upper(global) level, is particularly appealing. These features lead us to use the Global-Local approaches as they make it possible first to compute the homogeneous global model, and then to determine the critical areas to be re-analyzed,  while storing the factorization of the structural stiffness decomposition \cite{NoiiGL18,noii2020adaptive,aldakheel2020global}. The local model is then iteratively substituted within the unchanged/fixed global one, which avoids the reconstruction of the global mesh. Here, we propose an efficient \textit{multilevel Global-Local techniques for adaptive ductile phase-field fracture.}

To formulate the coupling of different levels within Global-Local scheme, a single Lagrange multiplier method leads to redundant interface conditions in the case of many
local domains (i.e., over-constrained condition for more than two domains, thus leads to the linear dependency of the imposed {constraints}), see for instance \cite{park2004partitioned} and references therein. Therefore, inconsistency conditions due
to the over-constrained interface displacement continuity appears and this leads to the non-unique
solutions \cite{park2002simple,park2004partitioned,flemisch2007stable,song2015gap}. But this is not the case for the localized Lagrange multiplier (LLM) approach which provides no redundancy
for the interface conditions and leads to unique and stable solutions \cite{park2004partitioned}.
Additionally, if the non-matching discrete interface is used, depending on which side of interface nodes for the single
Lagrange multiplier method are collocated, one would normally obtain different discrete constraint
equations. This issue has been extensively studied in the context of the mortar methods \cite{puso2003mesh,hansbo2005lagrange,seitz2016isogeometric}. In
contrast, LLM through the introduction of an intermediate surface on which both displacements and forces
are introduced as added variables thus offering a regularization of stiffness mismatch issues \cite{song2015gap}. This
is achieved by enforcing an additional weak from to our system of equations that are designed to satisfy
both the displacement compatibility condition and force equilibrium conditions \cite{song2015gap}.
Hence, a variational-based Global-Local approach is formulated based on the LLM \cite{park2002simple} method,
thus enables a straightforward extension of the proposed method for treating non-matching discrete
interfaces. The choice of the Dirac delta function for the Lagrange multiplier interpolation is also possible
which leads to more equilibrium state compared to the single Lagrange multiplier with the same
Lagrange multiplier interpolation function  \cite{song2015gap}. Thus enforcing a point-wise weak equality between
the global and local displacement fields which gives an advantage for the ease of computing the geometric
operators \cite{hautefeuille2012multi}.

In the past decade, both phase-field and Global-Local  approaches have been extended to deal with a growing number of situations of interest for engineers. The currently available phase-field formulations of brittle fracture encompass static and dynamic models. We mention the papers by Amor et al.\ \cite{amor+marigo+maurini09}, Miehe et al.\ \cite{miehe+welschinger+hofacker10a,MieWelHof10b}, Kuhn and M\"uller \cite{KUHN10}, Pham et al.\ \cite{pham+etal11}, Borden et al.\ \cite{borden+hughes+landis+verhoosel14}, Mesgarnejad et al.\ \cite{MESGARNEJAD2015420}, Kuhn et al.\ \cite{kuhn2015degradation}, Ambati et al.\ \cite{ambati15}, Wu et al. \cite{wu2017phase}, where various formulations are developed and validated. Recently, the framework has been also extended to ductile (elasto-plastic) fracture \cite{Duda2015,Ambati2015,Alessi2015,Borden2016,Miehe2016,Ambati2017,Karlo+fadi21}, fracture in films \cite{LeonBaldelli2014}, inverse problem \cite{khodadadian2019bayesian,noii2020bayesian}, anisotropic settings \cite{gultekin2018numerical,denli2020phase,teichtmeister2017phase}, and shells \cite{Shen2014}. Pressurized and fluid-filled fractures using phase-field modeling was subject in numerous papers in recent years. 
These studies range from 
mathematical modeling 
\cite{BourChuYo12,MiWheWi19,NoiiWick2019,CHUKWUDOZIE2019957,LeeWheWi16,LeeMiWheWi16}, 
mathematical analysis \cite{MiWheWi15b,MiWheWi14,MiWheWi15c}, 
numerical modeling and simulations 
\cite{Miehe2015186,HEIDER2018116,LeeMinWhe2018,Cajuhi2017,LeeWheWiSri17}, and 
up to (adaptive) Global-Local formulations \cite{aldakheel2020global} 
(see here in particular also \cite{GePlTuDo20} and \cite{noii2020adaptive}
for non-pressurized studies)
and high performance parallel computations 
\cite{HeiWi18_pamm,JoLaWi20}. A summary of multiphysics phase-field fracture was complied in \cite{Wi20_book}.

 Non-intrusive Global-Local approaches have also been applied to a quite large number of situations: the computation of the propagation of cracks in a sound model using the extended finite element method (XFEM) \cite{PAS13}, the computation of  assembly of plates introducing  realistic non-linear 3D modeling of connectors \cite{GUG16}, the extension to non-linear domain decomposition methods \cite{Duval2014} and to explicit dynamics \cite{BET14} with an application to the prediction of delamination under impact using ABAQUS \cite{BET17}. 

Recently, an adaptive Global-Local approach enhanced with a predictor-corrector scheme is {designed} in which the local domains are {\it dynamically updated} during the computation \cite{noii2020adaptive}. The predictor-corrector  methodology allows us to track a prior unknown crack paths \cite{Wick15Adapt,NoiiWick2019} . In the Global-Local framework, the fractures are prescribed in the local domain and once the fracture  grows further, the new local domains are predicted. Subsequently, the entire solution is corrected on the new Global-Local boundary value problems \cite{noii2020adaptive,aldakheel2020global} . The key requirement for realizing {this} adaptive Global-Local scheme is a non-matching discretization method on the interface. {To this end,}
a dual mortar method \cite{Wohlmuth,REIS2014168,popp2010dual} was used, thus 
providing sufficient regularity of the underlying {meshes}. Thereafter, 
an adaptive Global-Local formulation for pressurized fractures in the mechanics-step is derived in porous media at finite strain setting \cite{aldakheel2020global}. It has been shown that Global-Local framework for the poroelasticty material was up to 60 times faster than the standard phase-field formulation (single-scale solution), yet an excellent performance of the proposed framework was observed.

In the following, we describe in more detail our main goals.
First, we focus on the development of Global-Local
formulation for ductile fracture. {Specifically}, the continuum phase-field approach to ductile fracture
is employed. The coupling of plasticity to the crack phase-field is realized by a constitutive work density function, which is characterized through a degraded stored elastic energy and the accumulated dissipated energy due to plasticity and damage.  Two different Global-Local formulations based on the idea of multiplicative Schwarz'  method \cite{toselli2006domain} are proposed. In the first model, a global constitutive model behaves as an elastic-plastic response, while it is enhanced with a single local domain, which, in turn, describes an elastic-plastic fracturing model. Thereafter, we focus on the key goal of this contribution, by describing the second Global-Local formulation. The main objective of this extension is to introduce an adoption of the Global-Local approach toward the multilevel local setting.  More precisely, in this setting, a global constitutive model behaves as an elastic-plastic response, which is augmented with \textit{two distinct local domains}. The first local domain behaves as an elastic-plastic material and the next local domain is responsible for the \textit{only} crack phase-field formulation. A successful extension of this setting results in a reduction of the simulation time while preserving the computational accuracy. 

The second objective is to introduce  a predictor-corrector
adaptivity within the Global-Local concept to reduce the computational time. By applying predictor-corrector steps, a better estimation for (i) the elastic-plastic response, and (ii) the fracture state before proceeding to the next time step are achieved (since more elements are locally resolved). Hereby, two proposed Global-Local scenarios are augmented with their own adaptivity scheme. In the first Global-Local framework, an adaptivity procedure is realized through an equivalent global plastic strain. In this regard, global enrichment elements are applied to those regions where deviatoric stress reaches to the specified yield surface. This is referred to the \textit{one way adaptivity scheme}. In the second Global-Local  formulation, we have \textit{two-way} of adaptivity procedure. Specifically, we have one adaptivity step for the plasticity state in which a refinement knowledge coming from a global hardening state. The second  adaptivity step deals with fracture state such that a refinement knowledge is the result of the lowest local level, hence crack phase-field is the source of second adaptivity. The proposed adaptive multilevel Global-Local approach for the ductile phase-field fracture save reasonably computational cost (for solving the crack phase-field) since there is no need to do phase-field adaptivity prior to the onset of fracture.

To examine the proposed predictor-corrector adaptivity techniques, three different indicators are used to verify the overall response of the Global-Local approach in terms of accuracy/robustness and efficiency compared to the single-scale problems. These indicators are (i) load-displacement curve (which is observed globally), (ii) local crack phase-field pattern, and (iii) local hardening evolution (which are observed locally). In summary, this work extends the adaptive Global-Local phase-field fracture approach in \cite{noii2020adaptive,aldakheel2020global} to ductile phase-field fracture applications. The main objective of this extension is to introduce:

\begin{itemize}
	\item An extension of the Global-Local approach for a ductile phase-field fracture; 
	\item A multilevel of the Global-Local formulation through two distinct local domains;
    \item One way adaptivity  scheme based on evolution of global plastic flow;
    \item Two-way adaptivity procedure through the information of global hardening flow along with local crack phase-field state.
\end{itemize}

The paper is structured as follows: For a better insight into Global-Local formulation, in Section 2, one dimensional analysis for a simple horizontal bar is provided. We substantiate our derivation with an {\it open-source code}. Next, in Section 3, we outline the variational phase-field formulation of ductile fracture. Then, in Section 4, the extension to a Global-Local formulation for phase-field formulation of ductile fracture is derived.  Further extension toward multilevel Global-Local formulation is explained. In section 5, a robust and efficient predictor-corrector Global-Local adaptive approach is further
developed. In Section 5, four numerical results are performed in order to demonstrate our algorithmic developments. Finally, the last section concludes the paper with some remarks.

\begin{figure}[!t]
	\centering
	{\includegraphics[clip,trim=1cm 0.8cm 1cm 1cm, width=15cm]{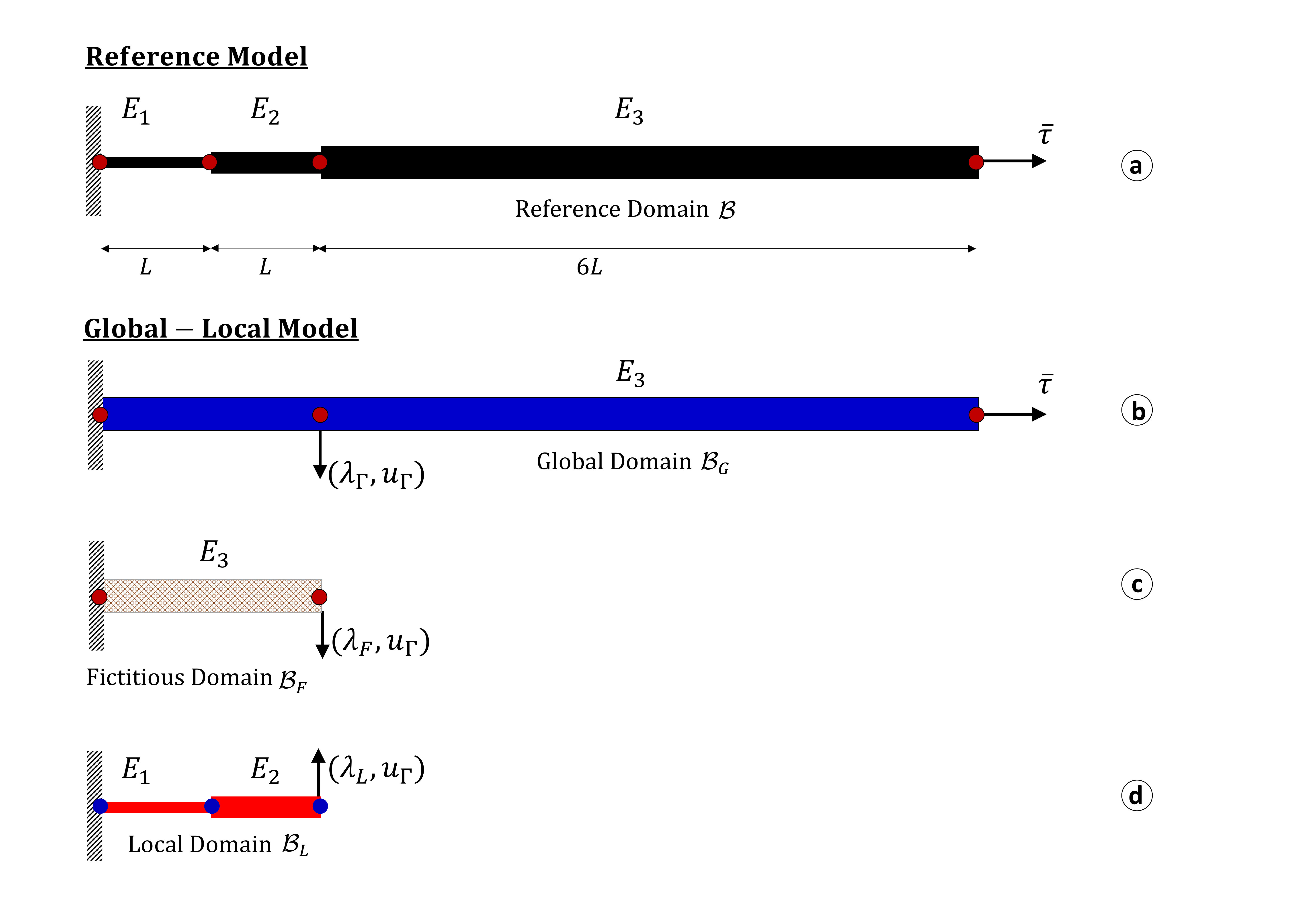}}  
	\caption{ Geometry, loading setup, and discretization for the one-dimensional bar. (a) Reference domain, (b) global domain, (c) fictitious domain, and (d) local domain. Nodal points due to the discretization are depicted for each geometry.}
	\label{Fig_Ch_intro_2}
\end{figure}

\sectpa[Section2]{One-Dimensional Analysis for the Global-Local Formulation}
This section provides a brief illustrative one-dimensional analysis for the Global-Local formulation. Detailed theoretical variational formulation are provided in next sections. To this end, the energy functional $\mathcal{E}(u)$ for linear elasticity is given by
\begin{equation}\label{1D1}
\begin{aligned}
\mathcal{E}(u)=\int_{\calB}\frac{1}{2}E(x) (u^{\prime})^2 \, A\,
\mathrm{d}{{x}}-\int_{{\partial_N\calB }} { {\bar\tau}} \cdot  u\,\mathrm{d}s,
\end{aligned}
\end{equation}
where ${\bar\tau}$ is the applied traction at the Neumann boundary and $E$ is Young's modulus. Let us now consider one-dimensional boundary value problem (BVP) that is shown in Fig. \ref{Fig_Ch_intro_2}a. We {depict} this as \textit{a reference} BVP such that its discretized setting includes three elements and four nodal points with a length of $8L$, see Fig. \ref{Fig_Ch_intro_2}. The cross-sectional area $A$ is used as an identical unit area through the entire bar. Hence, \req{1D1} can be rewritten as
\begin{equation}\label{1D2}
\begin{aligned}
\mathcal{E}(u)=\int^{8L}_{0}\frac{1}{2}E(x)(u^{\prime})^2\,A\,
\mathrm{d}{{x}}-\int_{{\partial_N\calB }} { {\bar\tau}} \cdot  u\,\mathrm{d}s .
\end{aligned}
\end{equation}
The function $E(x)$ is shown in Fig. \ref{Fig_Ch_intro_2}, hence we have
\begin{itemize}
\item $E(x)=E_1$ \quad for \quad $0\le x\le L$\ ,
\item $E(x)=E_2$ \quad for \quad $L< x\le 2L$\ ,
\item $E(x)=E_3$ \quad for \quad $2L< x\le 8L$\ .
\end{itemize}
The minimization of the given one-dimensional linear elasticity \req{1D2} leads to the Euler-Lagrange equation given by
\begin{align*}
{\mathcal E}_{ u}( u;\delta { u}_G)&:=
\int^{8L}_{0}E(x)u^{\prime} \delta u^{\prime}\,  A\,\mathrm{d}{{x}}
-\int_{\partial_N\calB } { {\bar\tau}} \cdot \delta{ u}\,\mathrm{d}s=0,
\label{Reference}
\tag{{R}}
\end{align*}
where ${\mathcal E}_{ u}$ is the directional derivative of the energy functional ${\mathcal E}$ with respect to the displacement $u$. Here, $\delta u \in  {H}_0^1(0,8L)$ is a test function. We now aim to resolve $(\textbf{R})$ using the Global-Local formulation. In this regard, the corresponding Global BVP is given in Fig. \ref{Fig_Ch_intro_2}b. It is depicted as \textit{a global} BVP such that its discretized setting includes two elements and three nodal points with a length of $8L$. Here, a homogenized Young's modulus $E_3$ is considered for the entire global domain, thus $E_G=E_3$ at $0\le x\le 8L$. Accordingly, a local BVP is given in Fig. \ref{Fig_Ch_intro_2}d.  We refer to this as \textit{a local} BVP such that its discretized setting includes two elements and three nodal points with a length of $2L$.
Coarse representation of the local domain within the global level is the so-called
\textit{fictitious} domain; see Fig. \ref{Fig_Ch_intro_2}c. A global variational equation is defined to find $u_G\in{H}_0^1(0,8L)$ through
\begin{align*}
\widetilde{\mathcal E}_{ u_G}( u_G;\delta { u}_G)&:=
\underbrace{\int^{8L}_{0}Eu_G^{\prime} \delta u_G^{\prime}\,\mathrm{d}{{x}}
-\int_{\Gamma_{N,G} } { {\bar\tau}} \cdot \delta{ u}_G\,\mathrm{d}s}_{\text{standard terms}}
\underbrace{-\int_{\Gamma_G} \lambda_\Gamma \cdot \delta { u}_G\,\mathrm{d}s}_{\text{jump term}}=0.
\label{Globalintro}
\tag{G}
\end{align*}
\begin{figure}[t]
	\centering
	{\includegraphics[clip,trim=1cm 4cm 1cm 5cm, width=15cm]{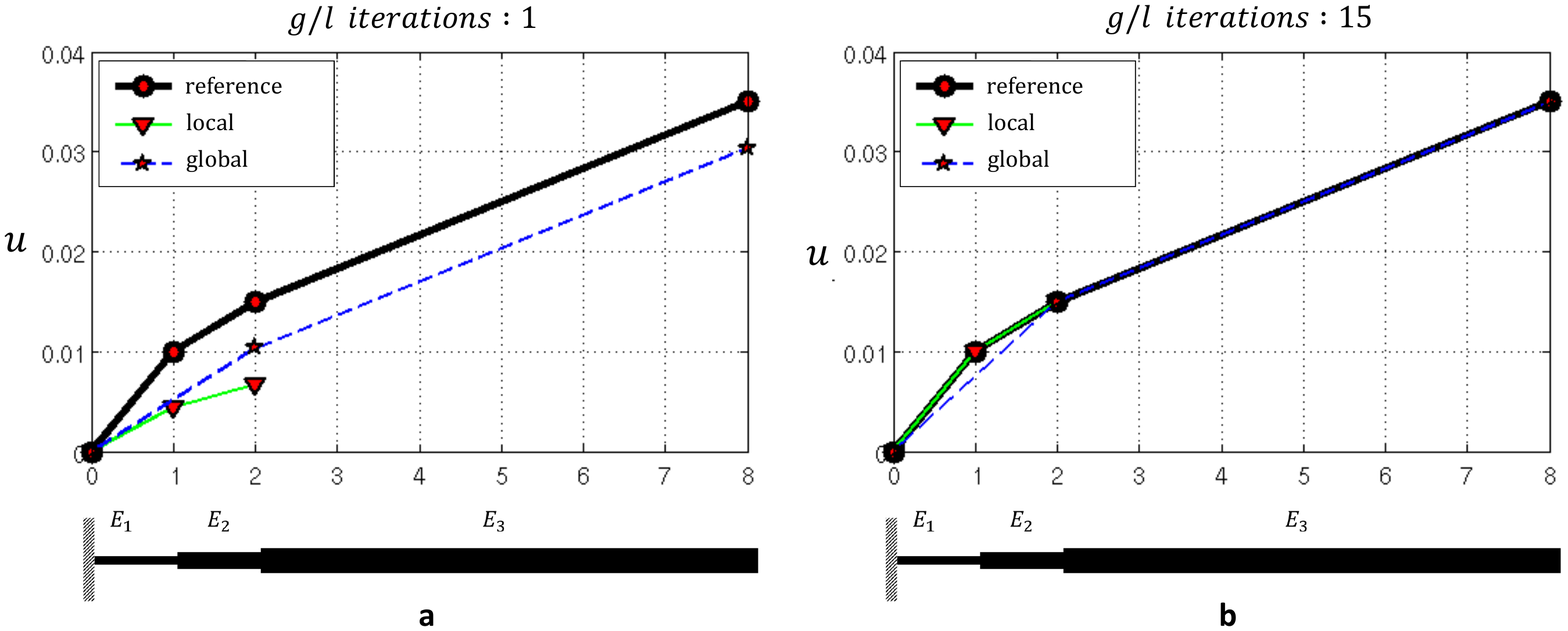}}  
	\caption{ Displacement distribution along the bar: reference, global and local solutions. (a) Global-Local solutions at the first iteration, and (b) Global-Local solutions at the $15^{th}$ iteration.}
	\label{Fig_Ch_intro_3}
\end{figure}
\begin{figure}[b]
	\centering
	{\includegraphics[clip,trim=1cm 9cm 1cm 11cm, width=16cm]{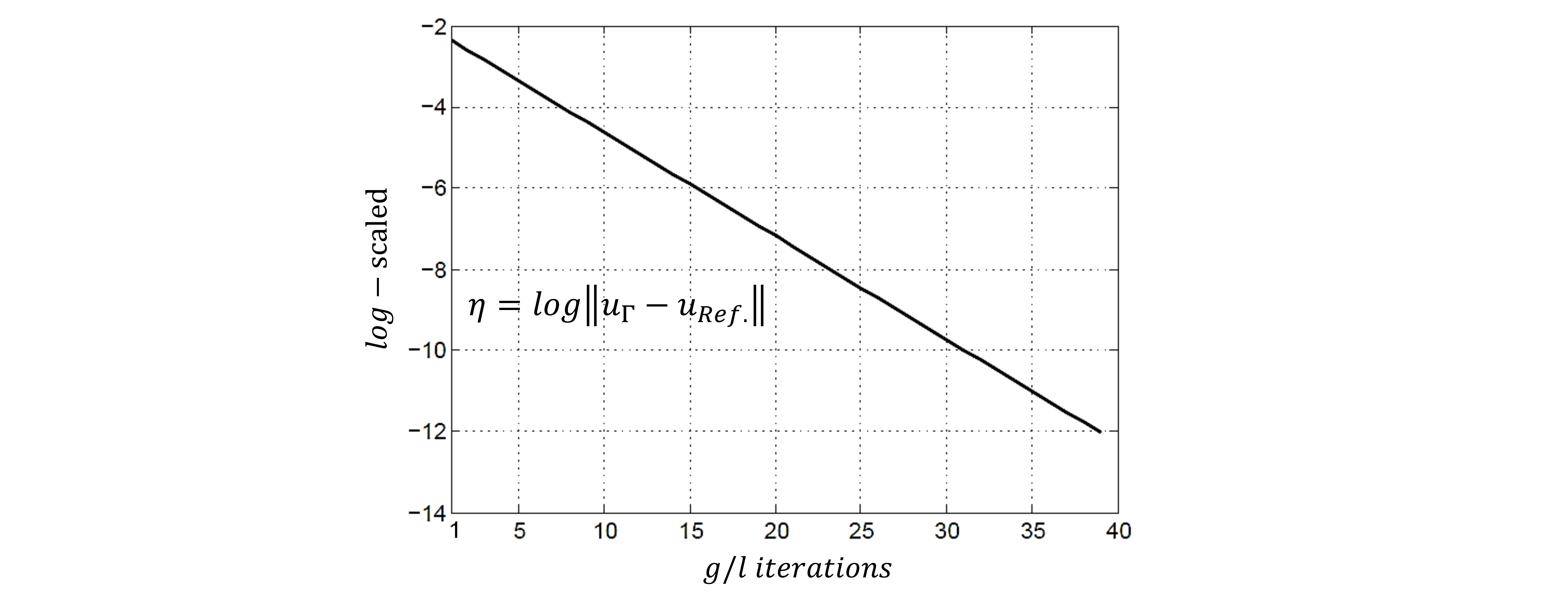}}  
	\caption{ Convergence behavior of the Global-Local formulation for the 1D BVP.}
	\label{Fig_Ch_intro_4}
\end{figure}
Here, $\lambda_\Gamma$ means the interface residual for measuring the discrepancy between global and local solutions at the interface (i.e., global nodal point 2), which in turn enters the global scale problem as a source term, thereby enabling an update of the global solution. An interface residual quantity as a traction jump between the fictitious and local domains takes the following form
\begin{align*}
\lambda_\Gamma(x) =
\lambda_F(x)-\lambda_L(x) \quad \text{at} \quad x_G=2L ,
\label{Global2}
\end{align*}
where $(\lambda_F,\lambda_L)\in L_2$ are given fictitious and local traction quantities at the global level through the previous solution field. To ensure displacement continuity between global and local domains, the resulting global displacement field at the interface, called $u_\Gamma$,  is imposed on the local BVP; hence we have a constrained local BVP to find $(u_L,\lambda_L)\in(H^1_0(0,2L),L_2)$ using
\begin{equation*}\label{alLocLintro}
\left\{
\begin{tabular}{l}
$\widetilde{\mathcal E}_{ u_L}( u_L,\lambda_L;\delta { u}_L):=
\displaystyle\int^{2L}_{0}E(x)u_L^{\prime} \delta u_L^{\prime}\,\mathrm{d}{{x}}
-\int_{\Gamma_G} \lambda_L \cdot \delta { u}_L\,\mathrm{d}s=0,$\\
$\widetilde{\mathcal E}_{\lambda_L}(u_L,\lambda_L;\delta { \lambda}_L):=
u_\Gamma- u^{3}_L=0.$
\end{tabular}
\right.
\tag{L}
\end{equation*}
Here, $u^3_L$ stand for the third node in local BVP in Fig. \ref{Fig_Ch_intro_2}. Two BVPs, namely $(G)$ and $(L)$ have to be solved in an iterative manner such that convergence is ensured. Convergence is achieved when both displacement and traction continuity along the interface are held. To evaluate the Global-Local formulation, the BVP given in Fig. \ref{Fig_Ch_intro_2} is considered. We set $A=1\;m^2$, $L=1\;m$, and $(E_2,E_3)=(2E_1,3E_1)$ with $E_1=10$. The resulting displacement distribution for the reference, global and local BVPs are provided in Fig. \ref{Fig_Ch_intro_3} for different iterations. After 15 iterations, the Global-Local formulation indeed recovers the displacement solutions corresponds to the reference one.
Figure \ref{Fig_Ch_intro_4} illustrates the convergence behavior of the Global-local iterative procedure for the one-dimensional BVP given in Fig \ref{Fig_Ch_intro_2}. We observe that Global-Local formulation is reached to the convergence state, i.e., $\Vert u_{\Gamma}-u_{Ref.}\Vert_2<\texttt{TOL}$, after 39 iterations. Here, we set $\texttt{TOL}=10^{-12}$. The compact open-source code that can be used to reproduce this example will be available online at \textcolor{blue}{\url{https://github.com/IKM-LUH/Noii-Aldakheel}}, and given in Appendix A.

\sectpa[Section3]{Variational Phase-Field Ductile Fracture}
After the introductory motivated 1D elasticity accomplished by introducing GL techniques, we now summarize the material model of phase-field ductile fracture. This will be next analyzed  using the Global-Local approach. The formulation here is based on a minimization of a pseudo-potential energy for the coupled problem undergoing small strains. 

\begin{figure}[!t]
	\centering
	{\includegraphics[clip,trim=8cm 13cm 1cm 5cm, width=18cm]{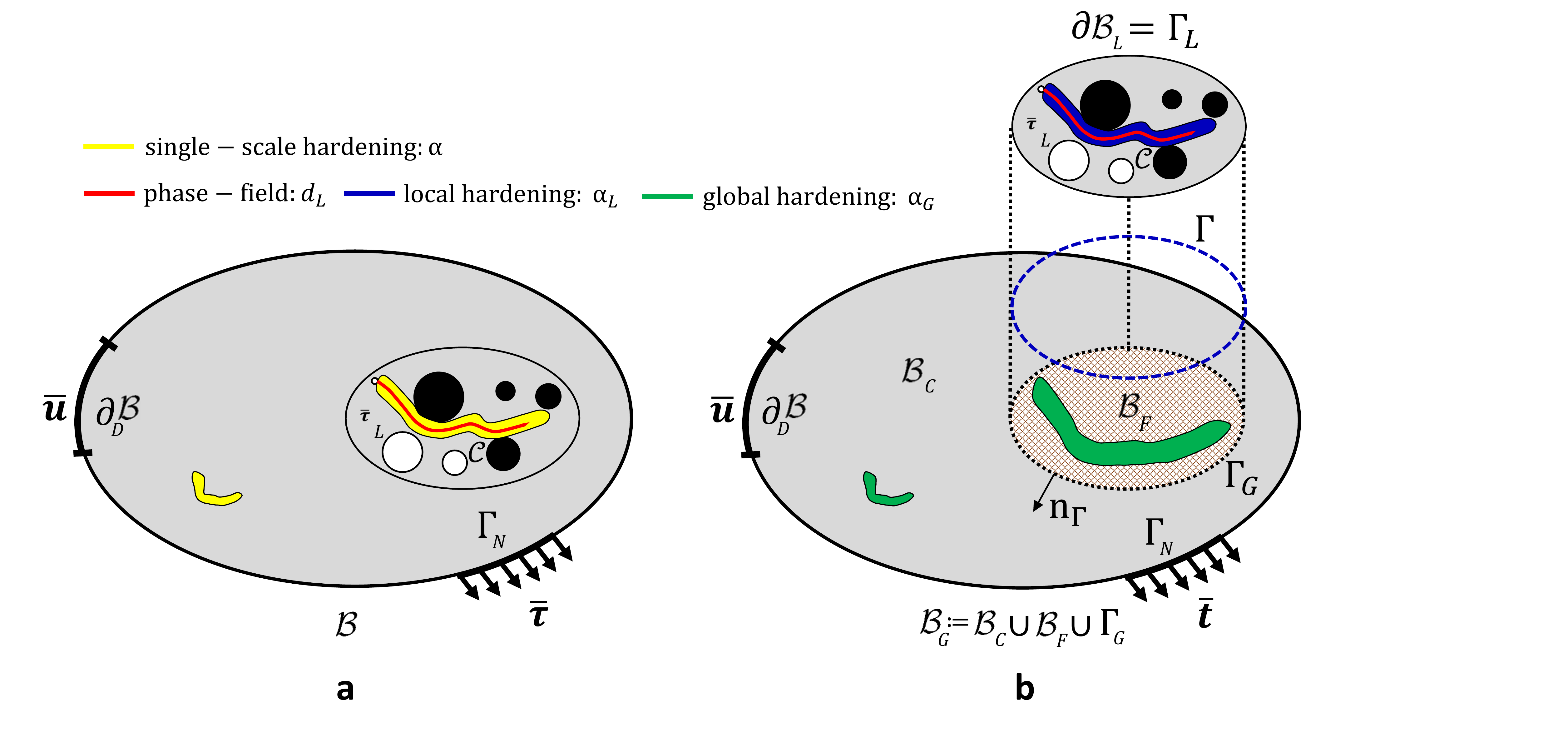}}  
	\caption{Solid with a crack inside of a plastic zone and boundary conditions. (a)
		Single-scale domain, and (b) Global-Local boundary value problem.}
	\label{Figure1}
\end{figure}

\sectpb[Section31]{Basic kinematics}
Let $\calB\in{\calR}^{2}$ be a solid domain with $\partial\calB$ denoted as its boundary. We assume a Dirichlet boundaries conditions
$\partial_D\calB $ and Neumann condition on $\partial_N \calB := \Gamma_N \cup \mathcal{C}$, where $\Gamma_N$  denotes the outer domain boundary and the lower dimensional fracture $\calC\in \calR$ is the crack boundary, as illustrated in Fig. \ref{Figure1}a. The response of fracturing solid at material points $\Bx\in\calB$ and time $t\in \calT = [0,T]$ is described by the displacement field $\Bu(\Bx,t)$ and the crack phase-field $d(\Bx,t)$ as
\begin{equation}
\Bu: 
\left\{
\begin{array}{ll}
\calB \times \calT \rightarrow \calR^2 \\
(\Bx, t)  \mapsto \Bu(\Bx,t)
\end{array}
\right.
\AND
d: 
\left\{
\begin{array}{ll}
\calB \times \calT \rightarrow [0,1] \\
(\Bx, t)  \mapsto d(\Bx,t)
\end{array}
\right.
\WITH
\dot{d} \ge 0\ ,
\label{s2-fields}
\end{equation}
where $d(\Bx,t)=0$ and $d(\Bx,t)=1$ describe the unbroken and fully fractured state of the material, respectively. The fracture surface $\mathcal{C}$ is approximated in $\calB_L\subset\calB$ so-called \textit{local domain}. Thus, $\calB_L$ represent the domain in which the smeared crack phase-field is approximated, and its boundary $\partial \calB_L$ depend on the choice of the phase-field regularization parameter $l>0$. The intact region with no fracture is denoted as \textit{complementary domain} $\calB_C:=\calB \backslash \calB_L\subset\calB$, such that ${\calB}_C\cup{\calB}_L=:\calB$ 
and ${\calB_C}\cap{\calB_L}=\varnothing$. 

The gradient of the displacement field defines the symmetric strain tensor of the geometrically linear theory as
\begin{equation}
\Bve = \nabla_s \Bu = \sym[ \nabla \Bu ] := \frac{1}{2} [\nabla\Bu + \nabla\Bu^T]
\ .
\label{s2-disp-grad}
\end{equation}
Focusing on the isochoric setting of von Mises plasticity theory, the strain tensor is additively decomposed into an elastic $\Bve^e$ and a plastic part $\Bve^p$ as
\begin{equation}
\Bve = \Bve^e + \Bve^p 
\WITH
\tr{[\Bve^p]} = 0
\AND
\tr{[\Bve]} = \tr{[\Bve^e]}
\ ,
\label{s2-strain-e-p}
\end{equation}
where the plastic strain is considered as the first local internal variable. To account for phenomenological hardening/softening response, we define the equivalent plastic strain variable by the evolution equation
\begin{equation}
\dot\alpha = \dot{\gamma}
\WITH
\dot\alpha \ge 0
\ ,
\label{s2-evol-alpha}
\end{equation}
as a second local internal variable, where $\dot\gamma\ge 0$ is the plastic Lagrange multiplier. The hardening variable starts to evolve from the initial condition $\alpha(\Bx,0) = \mathit{0}$.
\begin{Remark}
\label{assume}
In this work the elastic-plastic material behavior is considered in both the {\it complimentary} as well as the {\it local domains} whereas the fracture response lives locally at the lower scale $\calB_L$.  
\end{Remark}
The solid $\calB$ is loaded by prescribed deformations and external traction on the boundary, defined by time-dependent Dirichlet- and Neumann conditions
\begin{equation}
\Bu = \overline{\Bu} \ \textrm{on}\ \partial_D\calB
\AND
\Bsigma \cdot \Bn 
= \overline{\Btau} \ \textrm{on}\ \partial_N\calB\ ,
\label{s2-bcs}
\end{equation}
where $\Bn$ is the outward unit normal vector on the surface $\partial \calB$. The stress tensor $\Bsigma$ is the thermodynamic dual to $\Bve$ and $\bar{\Btau}$ is the prescribed traction vector.

For the phase-field problem, a sharp-crack surface topology $\calC \rightarrow \calC_l$ is regularized by the crack surface functional as outlined in \cite{miehe+hofacker+schaenzel+aldakheel15}
\begin{equation}
\calC_l(d) = \int_{\calB} \gamma_l(d, \nabla d) \, dv
\WITH
\gamma_l(d, \nabla d) =  
\dfrac{1}{2l} d^2 + \dfrac{l}{2} \vert \nabla d \vert^2\ , 
\label{s2-gamma_l}
\end{equation}
based on the crack surface density function $\gamma_l$ per unit volume of the solid and the fracture length scale parameter $l$ that governs the regularization.
Evolution of the regularized crack surface functional \req{s2-gamma_l} can be driven by the constitutive
functions as outlined in \cite{aldakheel+wriggers+miehe18}, postulating a global evolution equation of regularized crack surface as
\begin{equation}
\vphantom{\int_{\calB}}
\frac{d}{dt} \calC_l(d) =
\int_{\calB} {\delta_d}
\gamma_l(d,\nabla d) \, \dot d \, dv := 
\frac{1}{l} \int_{\calB} [\; (1-d) \calH - \eta_f \dot{d}\; ]\;
\dot{d} \, dv \ge 0\ ,
\label{gamma-evol}
\end{equation}
where $\eta_f \ge 0$ is a material parameter that characterizes 
\textcolor{black}{the artificial/numerical viscosity of the crack propagation}. The crack driving force
\begin{equation}
\calH = \max_{s\in [0,t]} D(\Bx,s) \ge 0
\ ,
\label{driving-force}
\end{equation}
is introduced as the third {\it local history variable} that accounts on the irreversibility of the phase-field evolution by filtering out a maximum value of what is known as the crack driving state function $D$.

\sectpb[Section32]{Constitutive work density function}

The ductile failure response of a solid is based on the displacement field $\Bu$ and the crack phase-field $d$ as global primary fields. Hence, the constitutive approach focuses on the set
\begin{equation}
\mbox{Constitutive State Variables}\; 
\BfrakC := \{ \Bve, \Bve^p, \alpha, d, \nabla d \}\ ,
\label{state}
\end{equation}
representing a combination of elasto-plasticity with a first-order gradient damage modeling. It is based on the definition of a pseudo-energy density per unit volume contains
\begin{equation}
{W}(\BfrakC) = 
	{W}_{elas}(\Bve^e, d) + {W}_{plas}(\alpha, d) +
	{W}_{frac}(d, \nabla d)\ ,
\label{pseudo-energy}
\end{equation}
the sum of a degrading elastic ${W}_{elas}$ and plastic energies ${W}_{plas}$ and a contribution due to fracture ${W}_{frac}$, which includes the accumulated dissipative energy. 
\\
\\
\smallskip
The {\it elastic contribution} is assumed to have the simple quadratic form
\begin{equation}
{W}_{elas}(\Bve^e, d)= g(d)\; \psi_e^{+}(\Bve^{e}_+) +  \psi_e^{-}(\Bve^{e}_-)
\WITH
\psi_e^{\pm}= \frac{\kappa}{2} \langle {\tr}[\Bve^e]\rangle^2_\pm + \mu \tr\big[\dev (\Bve^e_\pm)^2\big]\ ,
\label{elas-part}
\end{equation}
in terms of the bulk modulus $\kappa > 0$ and the shear modulus $\mu > 0$; characterizing an isotropic, linear stress response. The function ${g}(d) = (1-d)^2$ models the degradation of the elastic-plastic energy of the solid due to fracture. Hereby, a crack evolution only in tension is enforced by decomposing the stored elastic energy of the solid into a positive part ${\psi}{}^{+}$ due to tension and a negative part ${\psi}{}^{-}$ due to compression. This is given in terms of the two ramp functions $\langle x \rangle_\pm := (x \pm \vert
x\vert)/2$ of $\calR_\pm$, and the positive and negative elastic strain
tensors $\Bve^e_+ := \hbox{$\sum_{a=1}^3$} \langle\varepsilon^e_a
\rangle_+\;\Bn_a \otimes \Bn_a$ and $\Bve^e_- := \Bve^e -
\Bve^e_+$. $\{\varepsilon^e_a\}_{a=1,2,3}$ are the principal strains and $\{\Bn_a\}_{a=1,2,3}$ are the principal strain directions.
\\
\\
\smallskip
According, \cite{miehe+hofacker+schaenzel+aldakheel15}, the {\it plastic contribution} in \req{pseudo-energy} is assumed to have the form
\begin{equation}
{W}_{plas}(\alpha, d) = {g}(d)\; {\psi}_{p}(\alpha)
\WITH
{\psi}_{p} = Y_0\; \alpha + \frac{H}{2} \alpha^2 + (Y_\infty - Y_0)\big(\alpha + \exp[-\delta\alpha]/\delta\big)\ ,
\label{plas-part}
\end{equation}
with the initial yield stress $Y_0$, infinite yield stress $Y_\infty \ge Y_0$, the isotropic hardening modulus $H\ge 0$ and the saturation parameter $\delta$. Furthermore, this elasto-plastic model requires additionally the formulation of a yield function, a hardening law and an evolution equation for the plastic variables. The yield function restricts the elastic region. By assuming $J_2$-plasticity with nonlinear isotropic hardening the yield function has the form
\begin{equation}
\chi^p:=\chi^p(\Bu,\alpha,d)= \hbox{$\sqrt{3/2}$}\;\vert \BF^p \vert - R^p\ ,
\label{yield-fcn}
\end{equation}
with
\begin{equation}
	\BF^p:=\BF^p(\Bu,\alpha,d) = \dev[\Bsigma] = \Bsigma - \frac{1}{3} \mbox{tr} [\Bsigma] \Bone
	\AND
	R^p:=R^p(\Bu,\alpha,d) = \partial_\alpha {W}_{plas}\ ,
	\label{yield-fcn1}
\end{equation}
in terms of the deviatoric plastic driving force $\BF^p$ and the resistance force $R^p$. With the yield function at hand, we define the dual dissipation function for visco-plasticity according to Perzyna-type model as
\begin{equation}
\Phi^\ast(\BF^p, R^p) = \frac{1}{2\eta^p} 
\Big< \hbox{$\sqrt{3/2}$}\; \vert \BF^p \vert - R^p \Big>_+^2\ ,
\label{dissp-fcn}
\end{equation}
with $\eta^p$ being the viscosity parameter of the rate dependent plastic deformation. The evolution equations for the plastic variables are, see e.g. \cite{wriggers06}
\begin{equation}
\dot\Bve^p = \lambda^p  \,\BfrakN
\WITH
\BfrakN:=\frac{\partial \chi^p}{\partial \BF^p}
\AND
\dot\alpha = \lambda^p := \frac{1}{\eta^p} \big<\, \chi^p \,\big>_+\ ,
\label{plast-evol-eqs}
\end{equation}
The Kuhn-Tucker conditions for the elasto-plastic model are
\begin{equation}
\chi^p \le 0 , \quad\quad \lambda^p \ge 0, 
\quad\quad \mbox{and} \quad\quad
\chi^p\; \lambda^p = 0 \ .
\end{equation}
The {\it fracture part} of pseudo-energy density \req{pseudo-energy} takes the form
\begin{equation}
{W}_{frac}(d, \nabla d) =  [1-g(d)] \psi_c + 2 \frac{\psi_c}{\zeta}\; l \;{\gamma}_l(d, \nabla d)\ ,
\label{frac-part}
\end{equation}
where 
${\psi}_c > 0$ is a critical fracture energy and $\zeta$ controls the post-critical range after crack initialization.
\smallskip
\begin{form}[Energy functional for ductile phase-field fracture]
\label{form_1}
The development of a Global-Local approach for ductile phase-field fracture can start from a pseudo potential density functional as.
Let the initial conditions $\bm u_0=\bm u(\bm{x},0)$ and  $d_0=d(\bm{x},0)$ be given. For the loading increments $n=0,1,2,\ldots, N$, find $\bm u:=\bm u^{n+1}$ and $d:=d^{n+1}$ such that the functional for the coupled problem:
\begin{equation*}
{\mathcal{E}(\Bu, d) = \int_{\calB} {W}(\BfrakC) \, dv 
\; - \; 
\vphantom{\frac{d}{dt}}
\mathcal{E}_{ext} (\Bu)
\WITH
\mathcal{E}_{ext} (\Bu) := 
\int_{\calB} \overline\Bf \cdot \Bu\, dv  +
\int_{\partial_N\calB} \bar\Bt \cdot \Bu\, da}\ ,
\label{potential-functional}
\end{equation*}
{with ${W}(\BfrakC)$  given in \req{pseudo-energy}, is minimized.}
\end{form}

\sectpb[Section33]{Governing equations}

The minimization problem for the given energy functional of the inelastic crack topology in Formulation \ref{form_1} takes the following compact form:
\begin{equation}\label{compat_argmin}
\{ \bm{u}, d \} = 
\mbox{arg} \Big\{\substackrel{\bm{u}}{\mbox{min}}
\substackrel{d  }{\mbox{min}} \,
[\; \mathcal{E} (\bm u,d) \; ] \Big\}.
\end{equation}
The stationary points of the energy functional in Formulation \ref{form_1} are characterized
by the first-order necessary conditions, namely the so-called Euler-Lagrange equations, which are obtained by differentiation with respect to ${\bm u}$ and $d$ as follows:

\sectpc[Section331]{Balance of linear momentum}
The first equation is the stress equilibrium or the quasi-static form of the
balance of linear momentum defined as
\begin{equation}
\div\,\Bsigma + \overline\Bf = \Bzero\ ,
\label{s2-equil:defo}
\end{equation}
where dynamic effects are neglected and $\overline\Bf$ is the given body force. Following the Coleman-Noll procedure, the stress tensor is obtained from the potential ${W}_{elas}$ in \req{elas-part} by
 \begin{equation}
 \Bsigma := \partial_{\Bve^e} {W}_{elas}
 = g(d) \widetilde\Bsigma_{+} + \widetilde\Bsigma_{-}
 \WITH
 \widetilde\Bsigma_{\pm} = \kappa \langle\mbox{tr}[\Bve^e]\rangle_{\pm} \Bone 
 + 2 \mu \dev(\Bve^e_\pm)\ ,
 \label{vari-bf-23}
 \end{equation}
 where $\widetilde\Bsigma$ is the effective stress tensor.

\sectpc[Section332]{The fracture phase-field equation}

The evolution statement \req{gamma-evol} provides the second governing equation, representing the
evolution of the crack phase-field in the domain $\calB$ along with its homogeneous Neumann boundary condition as
\begin{equation}
\eta_f \dot{d}
= 
(1-d) {\calH} 
-
[ \, d - l^2 \Delta d \, ]
\qquad \mbox{with} \qquad
\nabla d \cdot \Bn = 0
\quad \mbox{on} \quad \partial\calB\ ,
\label{euler-eq-d}
\end{equation}
where the history field $\calH$ is defined by
\begin{equation}
 \calH := \max_{s\in [0,t]} D(\BfrakC ; s) \ge 0
 \WITH
 D := \zeta \bigg< \frac{{\psi}_e^{+} + {\psi}_{p}}{\psi_c} - 1 \bigg>_+\ ,
 \label{history-field}
\end{equation}
as outlined in \cite{aldakheel16}, with the Macaulay bracket $\langle x \rangle_+ := (x + \vert x\vert)/2$, that ensures the irreversibility of the crack evolution. 

\sectpa[Section4]{Extension Towards Global-Local Formulations}

For efficient and robust numerical solution procedures, a multi-scale approach is developed within this section, where the characteristic length of the lower scale is of the same order as its global part; see \cite{noii2020adaptive}. This is accomplished by introducing the so-called Global-Local approach for solving the above introduced system of equations obtained from Formulation \ref{form_1} for the coupled problem. The Global-Local (GL) method is rooted in the domain decomposition approach \cite{DDRey06}. In this regard, we proposed two different GL methods and compare them with the standard single scale formulations (Section \ref{Section3}) to illustrate their efficiency and capability for solving fracture mechanics problems numerically. 

Let the material body $\calB$ is decomposed into a global domain $\calB_G$ illustrating an elastic-plastic material behavior and a local domain $\calB_L$ reflecting the ductile fracture region. The global domain $\calB_G:= \calB_C\cup\calB_f\cup\Gamma$ is further split into a complementary domain $\calB_C$ corresponds to the intact area, a fictitious domain $\calB_f$ depicts a coarse projection of the local domain into the global one and an interface $\Gamma$
between the unfractured and the fractured domains. The {fictitious domain} $\calB_{f}$ is a prolongation of $\calB_{C}$ towards ${\calB}$, i.e. {\it recovering the space} of ${\calB}$ that is obtained by removing $\calB_{L}$ from its continuum domain, see Fig. \ref{Figure1}b. This gives the same constitutive modeling used in $\calB_{C}$ for $\calB_{f}$. We also use the same discretization space for both $\calB_{f}$ and $\calB_{C}$, which results in identical element size i.e. $h_f:=h_C$. The external loads are applied on $\calB_{C}$ and hence $\calB_{L}$ is assumed to be divergence-free. Such assumption is standard for the multi-scale setting, see \cite{Fish2014}.

At the interface $\Gamma$, Global and Local interfaces denoted as ${\Gamma}_{G}\subset\calB_{G}$ and $\Gamma_{L}\subset\calB_{L}$ are defined, such that in the continuum setting we have ${\Gamma}=\Gamma_{G}=\Gamma_{L}$. Hence, the displacement field $\Bu$ for both Global and Local domains do exactly coincide in the strong sense at the interface, yielding
\begin{equation}
{\Bu}_L(\Bx,t)\overset{!}{=}{\Bu}_G(\Bx,t) 
\quad  \mbox{at} \quad {\Bx}\in\Gamma\ .
\label{uCont}
\end{equation}
However, in a discrete setting we might have ${\Gamma}\neq\Gamma_{G}\neq\Gamma_{L}$ due to the presence of different meshing schemes (i.e. different element size/type used in $\calB_{G}$ and $\calB_{L}$ such that $h \neq h_L \neq h_G$ on ${\Gamma}$). As outlined in \cite{Farhat1991,NoiiGL18}, such strong continuity requirements defined above are too restrictive from the computational standpoint. To overcome these difficulties, we introduce the displacement interface $\Bu_\Gamma(\Bx,t)$ and its corresponding traction forces $\{\Blambda_L, \Blambda_C\}$ that are introduced as Lagrange multipliers. This results in a set of equations at the interface:
\begin{align*}
	\left\{
	\begin{tabular}{ll}
	${\Bu}_L(\Bx,t)  = {\Bu}_{\Gamma}(\Bx,t) $ & \qquad \mbox{at} ${\Bx}\in{\Gamma}_{L}$, \\[0.15cm]
	${\Bu}_G(\Bx,t)  = {\Bu}_{\Gamma}(\Bx,t) $ & \qquad \mbox{at} ${\Bx}\in{\Gamma}_{G}$, \\[0.15cm]
	$\Blambda_L(\Bx,t) + \Blambda_C(\Bx,t) = \Bzero$ & \qquad \mbox{at} ${\Bx}\in{\Gamma}$.
	\end{tabular}
	\right.
\end{align*}
Accordingly, the {\it single-scale} displacement field ${\Bu}(\Bx,t)$ in Section \ref{Section3} is decomposed as
\begin{equation}
{\Bu}(\Bx,t) =\left\{
\begin{tabular}{ll}
${\Bu}_L(\Bx,t)$ & for ${\Bx}\in\calB_{L}$, \\[0.15cm]
${\Bu}_G(\Bx,t)$ & for ${\Bx}\in\calB_{G}$, \\[0.15cm]
${\Bu}_{\Gamma}(\Bx,t)$ & for ${\Bx}\in{\Gamma}$.
\end{tabular}
\right.
\label{u-DD}
\end{equation}
The fracture surface lives only in $\calB_{L}$. Hence we can introduce scalar-valued function $d_L(\Bx,t):\calB_{L}\rightarrow[0,1]$. The {\it single-scale phase-field} $d$ is then decomposed in the following representation
\begin{equation}
d(\Bx,t)  :=\left\{
\begin{tabular}{ll}
$d_L$ & for ${\Bx}\in\calB_{L}$, \\[0.1cm]
$0$ & for ${\Bx}\in\calB_{G}$.
\end{tabular}
\right.
\label{dL-DD}
\end{equation}

\sectpb[Section41]{Constitutive formulations for the Global-Local coupling system}
Now the multi-physics problem for the Global-Local approach is based on six primary fields to characterize the ductile fracture in elastic-plastic solids as
\begin{equation}
\mbox{Extended Primary Fields}: \BfrakP := \{ \Bu_G, \Bu_L, d_L, \Blambda_C, \Blambda_L, \Bu_\Gamma \}
\label{gl-fields}
\ .
\end{equation}
Based on the above introduced decompositions and the governing equations in Section \ref{Section3}, we describe here the variational formulation for the Global-Local coupling system. To this end, the {\em Global-Local} approximation of the single-scale energy functional $\mathcal{E}$ indicated in Formulation \ref{form_1} by
\begin{equation}
\begin{aligned}
\widetilde{\mathcal E}(\BfrakP):
&= \underbrace{\int_{\calB_G} W(\Bve_G, \Bve^p_G, \alpha_G, 0, \Bzero) \,\mathrm{d}{v}
	- \int_{\calB_F} W(\Bve_G, \Bve^p_G, \alpha_G, 0, \Bzero) \,\mathrm{d}{v}
	- 	\vphantom{\frac{d}{dt}}\mathcal{E}_{ext} (\Bu_G)}_{\text{Global terms}} \,\\
&+\underbrace{\int_{\calB_L} W(\Bve_L, \Bve^p_L, \alpha_L, d_L, \nabla d_L) \,\mathrm{d}{v}}_{\text{Local term}}\\
&+\underbrace{\int_\Gamma \Big\{ \bm\lambda_C\cdot(\bm u_\Gamma-\bm u_G) + \bm\lambda_L\cdot(\bm u_\Gamma-\bm u_L) \Big\} \mathrm{d}a}_{\text{Coupling terms}},
\label{EGL}
\end{aligned}
\end{equation}
where the approximation ${\mathcal E} \equiv \widetilde{\mathcal E}$ holds. By neglecting volume forces $\overline\Bf$, the external load functional is defined as $\mathcal{E}_{ext} (\Bu_G) =\int_{\Gamma_{N,G} } {\bm {\bar\tau}} \cdot \delta{\bm u}_G\,\mathrm{d}a$. Next the same procedure applied in Section \ref{Section33} is followed here, such that the Global-Local energy functional \req{EGL} is minimized. The outcome minimization problem for the Global-Local energy functional that is applied to the phase-field modeling of fracture in ductile solids takes the following compact form,
\begin{equation}\label{GL}
\BfrakP= 
\begin{aligned}
\mbox{arg} \Big\{ 
\substackrel{\bm u_G, \bm u_L, \bm u_\Gamma, {d_L} \; \; }{\mbox{min}}
\substackrel{\bm\lambda_C,\bm\lambda_L}{\mbox{max}} \,
\big[\; \widetilde{\mathcal{E}}(\Bu_G, \Bu_L, d_L, \Blambda_C, \Blambda_L, \Bu_\Gamma) \big] \Big\}\, ,
\end{aligned}
\end{equation}
where $\BfrakP $ represents the primary fields vector. Now Euler-Lagrange equations can be obtained by differentiation with respect to the unknowns. Hence, the directional derivatives of the functional $\widetilde{\mathcal E}$ with respect to $\bm u_G$ yield the global weak form as
\begin{align*}
\widetilde{\mathcal E}_{\bm u_G}(\BfrakP;\delta {\bm u}_G)&:=
\int_{\calB_G}\bm\sigma(\bm u_G):\bm\varepsilon(\delta {\bm u}_G)\,\mathrm{d}{v}
-\int_{\calB_F}\bm\sigma(\bm u_G):\bm\varepsilon(\delta {\bm u}_G)\,\mathrm{d}{v}\\
&-\int_{\Gamma_G} \bm\lambda_C \cdot \delta {\bm u}_G\,\mathrm{d}a
-\int_{\Gamma_{N,G} } {\bm {\bar\tau}} \cdot \delta{\bm u}_G\,\mathrm{d}a=0,
\label{Global}
\tag{G}
\end{align*}
where $\bm\sigma(\bm u_G):=\partial_{\Bve} {W}(\bm\varepsilon(\bm u_G), \Bve^p_G(\bm u_G), \alpha_G(\bm u_G), 0, \Bzero)=\widetilde\Bsigma$ and $\delta {\bm u}_G\in\{ {\bf H}^1(\calB_G): \delta {\bm u}_G=\bm 0 \; \mathrm{on} \; \partial_D\calB \}$ is the test function. The local weak formulations assume the form
\begin{equation*}\label{alLocL}
\left\{
\begin{tabular}{l}
$\widetilde{\mathcal E}_{\bm u_L}(\BfrakP;\delta {\bm u}_L):=\displaystyle \int_{\calB_L}\bm\sigma(\bm u_L,d_L,\alpha_L):\bm\varepsilon(\delta {\bm u}_L)\,\mathrm{d}{v}
-\int_{\Gamma_L} \bm\lambda_L \cdot \delta {\bm u}_L\,\mathrm{d}a=0$, \\[0.1cm]
$\displaystyle \widetilde{\mathcal E}_{d_L}(\BfrakP;\delta d_L):=
\int_{\calB_L} \Big( d_L + \frac{\eta_f}{\Delta t} (d_L - d^{n}_L) + (d_L-1) \;\calH(\bm u_L,\alpha_L) \Big) \delta d_L \;\mathrm{d}{v}$\\[0.1cm]
\qquad \qquad \; \: $ \displaystyle + \int_{\calB_L} l^2\; \nabla d_L.\nabla(\delta d_L)\,\mathrm{d}{v}= 0$,\\[0.1cm]
\end{tabular}
\right.
\tag{L}
\label{alLoc}
\end{equation*}
where $\bm\sigma(\bm u_L,d_L)=\partial_{\bm\varepsilon}W(\bm\varepsilon(\bm u_L), \Bve^p_L(\bm u_L), \alpha_L(\bm u_L), d_L, \nabla d_L)$, $\delta {\bm u}_L\in {\bf H}^1(\calB_L)$ is the local test function, $\delta d_L\in \text{H}^1(\calB_L)$ is the local phase-field test function and $\Delta t := t - t_n > 0$ denotes the time step. The variational derivatives of $\widetilde{\mathcal E}$ with respect to $(\bm u_\Gamma,\bm\lambda_C,\bm\lambda_L)$ provide kinematic equations due to weak coupling between Global and Local form
\begin{equation*}
\widetilde{\mathcal E}_{\bm u_\Gamma}(\BfrakP;\delta {\bm u}_\Gamma):=\int_\Gamma (\bm\lambda_C+\bm\lambda_L) \cdot \delta {\bm u}_\Gamma \,\mathrm{d}a=0,
\label{Coupl1}
\tag{C$_1$}
\end{equation*}
\begin{equation*}
\widetilde{\mathcal E}_{\bm\lambda_C}(\BfrakP;\delta {\bm \lambda}_C):=\int_\Gamma (\bm u_\Gamma-\bm u_G) \cdot \delta {\bm \lambda}_C \,\mathrm{d}a=0,
\label{Coupl2}
\tag{C$_2$}
\end{equation*}
\begin{equation*}
\widetilde{\mathcal E}_{\bm\lambda_L}(\BfrakP;\delta {\bm \lambda}_L):=\int_\Gamma (\bm u_\Gamma-\bm u_L) \cdot \delta {\bm \lambda}_L \,\mathrm{d}a=0,
\label{Coupl3}
\tag{C$_3$}
\end{equation*}
herein, $\delta {\bm u}_\Gamma\in{\bf H}^1(\Gamma)$ and $\delta {\bm \lambda}_C,\delta {\bm \lambda}_L\in{\bf L}^2(\Gamma)$ are the corresponding test functions.

\sectpb[Section36]{Different Global-Local scenarios}

Equations (G), (L) and (C1)--(C3) specify the entire system of the Global-Local approach. To solve this system of equations an alternate minimization scheme is employed, where first the Local boundary value problem is computed in terms of given data from the upper scale then the Global boundary value problem will be solved in line with our previous works \cite{aldakheel2020global,NoiiGL18,noii2020adaptive}. In this regard, the global and local level are computed in a multiplicative manner according to the idea of Schwarz' alternating method \cite{Gander2007, Magoules2007}. Hereby, Dirichlet-Neumann type boundary conditions are considered to solve for the above introduced equations, as discussed in authors previous work, see \cite{NoiiGL18}. To relax the stiff local response that is observed at the global level (due to the local non-linearity), the Global-Local formulation is enhanced further using Robin-type boundary conditions, as outlined in \cite{aldakheel2020global,noii2020adaptive}.  Next we focus on the two Global-Local methods and compare their efficiency with respect to the single-scale solution.

\begin{Remark}
	\label{plasticity-choice}
Note, when the accumulated plastic strain is reached to its critical value (see Section \ref{Section4}), thus a refine approximation of the elastic-plastic response is needed. For instance, the right green plastic flow which is indicated in Fig. \ref{Figure1}, is accumulated enough to be considered for the fine approximation, which, in turn, the left plastic flow (which is shown with the small green region) is approximated with its coarse representation (in the sense of the Global-Local framework). 
\end{Remark}

In the following, let $n$ and $k$ indicate loading time step and Global-Local iterations, respectively. Note, we assumed the Global-Local formulation presented at the fixed loading time step, thus for the sake of simplicity, we omit index $n$ (but, we added where it is required).

\sectpc[Section361]{Global domain ($\texttt{G-EP}$) with a single local domain ($\texttt{L-EPD}$) }

In this part, the first Global-Local formulation denoted as $g/l-1$ will be discussed in details. Herein, a global constitutive model behaves as an elastic-plastic response, abbreviated as $E(elastic)-P(plastic)$, which is augmented with a \textit{single local domain} which behaves as an elastic-plastic material at fracture, abbreviated as $E(elastic)-P(plastic)-D(damage)$.

\sectpd[Section3611]{Robin-type boundary conditions at the local level}

Following our recent works (see \cite{noii2020adaptive,aldakheel2020global}), for the mechanical deformation field at the local level, a new coupling terms reads
\begin{equation*}
	\int_\Gamma \Blambda^{k}_L\cdot \delta {\Bu}_\Gamma \,\mathrm{d}a+\KIA_L\int_\Gamma \Bu_\Gamma^{k,\frac{1}{2}} \cdot \delta {\Blambda}_C \,\mathrm{d}a={\BLambda}^{k-1}_L,
	\label{RBC3}
	\tag{$\widetilde {\text{C}}_1$}
\end{equation*}
\begin{equation*}
	\int_\Gamma (\Bu_\Gamma^{k,\frac{1}{2}}-\Bu^{k}_L) \cdot \delta {\Blambda}_L \,\mathrm{d}a=0,
	\label{RBC4}
	\tag{$\widetilde {\text{C}}_2$}
\end{equation*}  
with
\begin{equation}
	{\BLambda}^{k-1}_L:=\BLambda_L(\Blambda^{k-1}_C,\Bu_G^{k-1};\KIA_L)=
	\KIA_L\int_\Gamma \Bu_G^{k-1} \cdot \delta {\Blambda}_C \,\mathrm{d}a-\int_\Gamma \Blambda^{k-1}_C\cdot \delta {\Bu}_\Gamma \,\mathrm{d}a \ .
	\label{RBC5}
\end{equation}
Here, the set $({\bf \Lambda}_L, \KIA_L)$ represent the local Robin-type parameters. 
To complete the local BVP, the principle of maximum dissipation (see e.g. \cite{wriggers06}) for a local elastic-plastic response  leads to the so-called local evolution equations for the plastic variables through
\begin{equation}
	\dot\Bve_L^p = \lambda_L^p  \,\BfrakN_L
	\WITH
	\BfrakN_L:=\frac{\partial \chi^p_L}{\partial \BF^p_L}
	\AND
	\dot\alpha_L = \lambda_L^p := \frac{1}{\eta^p_{L}} \big<\, \chi^p_L \,\big>_+\ .
	\label{plast-evol-eqsL}
\end{equation}
Along with the Kuhn-Tucker conditions for the elasto-plastic model as a local inequality constraints to our model by
\begin{equation}
	\chi^p_L \le 0 \quad\quad , \quad\quad \lambda_L^p=\dot\alpha_L\ge 0 
	\quad\quad \mbox{and} \quad\quad
	\chi^p_L\; \lambda_L^p = 0 \ ,
	\label{kktP}
\end{equation}
such that
\begin{equation}
	\chi^p_L=\chi^p(\Bu_L,\alpha_L,d_L) \AND	
	\BF^p_L=\BF^p(\Bu_L,\alpha_L,d_L),
\end{equation}
see \req{yield-fcn} and \req{yield-fcn1}. Accordingly,  the Kuhn-Tucker conditions for the gradient damage model reads
\begin{equation}
	\chi^d_L \le 0 \quad\quad , \quad\quad \lambda_L^d=\dot d_L\ge 0 
	\quad\quad \mbox{and} \quad\quad
	\chi^d_L\; \lambda_L^d = 0 \ .
	\label{kktD}
\end{equation}
such that the so-called damage yield function defined as
\begin{equation}
	\chi^d_L:=- \delta_{d_L} W_L\le 0 \WITH 
	W_L=W(\Bve_L, \Bve^p_L, \alpha_L, d_L, \nabla d_L),
\end{equation}
as outlined in \cite{aldakheel16,miehe2017phase}. The new coupling conditions ($\widetilde {\text{C}}_1$) and ($\widetilde {\text{C}}_2$), along with (\ref{alLocL}) and evolution equations \req{plast-evol-eqsL} which are imposed by two sets of inequality constraints in \req{kktP} and \req{kktD} introduce an {\it enhanced} local BVP. The local system of equations has to be solved for $(\bm u^{k}_L,d^{k}_L,\alpha^{k}_L,\bm\lambda^{k}_L,\bm u_\Gamma^{k,\frac{1}{2}})$ for given local Robin-type parameters $({\bf \Lambda}^{k-1}_L, \KIA_L)$. In summary, local BVP has the following abstract form,
\begin{equation}
	\texttt{L-EPD:}\quad(\bm u^{k+1}_L,d^{k+1}_L,\alpha^{k+1}_L,\bm\lambda^{k+1}_L,\bm u_\Gamma^{k+1,\frac{1}{2}})
	=\texttt{L}(\bm u^{k}_L,d^{k}_L,\alpha^{k}_L,\bm\lambda^{k}_L,\bm u_\Gamma^{k,\frac{1}{2}},{\bf \Lambda}^{k-1}_L, \KIA_L)\ ,
\end{equation}
such that, at $k=0$ the Global-Local formulation is initialized with previous converged solution at $n-1$.
\sectpd[Section3612]{Robin-type boundary conditions at the global level}

For the mechanical deformation field at the global level, a new set of coupling terms reads
\begin{equation*}
	\int_\Gamma \Blambda^{k}_C\cdot \delta {\Bu}_\Gamma \,\mathrm{d}a+\KIA_G\int_\Gamma \Bu_\Gamma^{k} \cdot \delta {\Blambda}_L  \,\mathrm{d}a={\BLambda}^{k}_G\ ,
	\label{RBC11}
	\tag{$\widetilde {\text{C}}_5$}
\end{equation*}
\begin{equation*}
	\int_\Gamma (\Bu_\Gamma^{k,\frac{1}{2}}-\Bu^{k}_G) \cdot \delta {\Blambda}_C  \,\mathrm{d}a=0\ ,
	\label{RBC12}
	\tag{$\widetilde {\text{C}}_6$}
\end{equation*}  
with
\begin{equation}
	{\BLambda}^{k}_G:=\BLambda_G(\Blambda^{k}_L,\Bu_L^{k};\KIA_G)=
	\KIA_G\int_\Gamma \Bu_L^{k} \cdot \delta {\Blambda}_L  \,\mathrm{d}a-\int_\Gamma \Blambda^{k}_L\cdot \delta {\Bu}_\Gamma \,\mathrm{d}a \ ,
	\label{RBC13}
\end{equation}
where the set  $({\bf \Lambda}_G, \KIA_G)$ are the global Robin-type parameters. Since the global constitutive model is formulated for $EP$ response, thus the global evolution equations for the plastic variables becomes
\begin{equation}
\dot\Bve_G^p = \lambda_G^p  \,\BfrakN_G
\WITH
\BfrakN_G:=\frac{\partial \chi^p_G}{\partial \BF^p_G}
\AND
\dot\alpha_G = \lambda_G^p := \frac{1}{\eta^p_{G}} \big<\, \chi^p_G \,\big>_+\ .
\label{plast-evol-eqsG}
\end{equation}
The Kuhn-Tucker conditions for the elasto-plastic model are
\begin{equation}
	\chi^p_G\le 0 \quad\quad , \quad\quad \lambda_G^p=\dot\alpha_G \ge 0 
	\quad\quad \mbox{and} \quad\quad
	\chi^p_G\; \lambda_G^p = 0 \ ,
	\label{kktPG}
\end{equation}
such that
\begin{equation}
	\chi^p_G=\chi(\Bu_G,\alpha_G,0) \AND	
	\BF^p_G=\BF^p(\Bu_G,\alpha_G,0)\ .
\end{equation}
The new global coupling conditions ($\widetilde {\text{C}}_3$) and ($\widetilde {\text{C}}_4$), together with (\ref{Global}) and evolution equations \req{plast-evol-eqsG} which are imposed by global loading-unloading inequality constraint in \req{kktPG} are introducing an enhanced global BVP. The global system of equations has to be solved for $(\bm u^{k}_G,\alpha^{k}_G,\bm\lambda^{k}_C,\bm u_\Gamma^{k,\frac{1}{2}})$ for given global Robin-type parameters $({\bf \Lambda}^{k}_G, \KIA_G)$.
In a sake of completeness, we write BVP for $\calB_{G}$ which deals with elastic-plastic response through the following abstract form,
\begin{equation}\label{BVPG}
	\texttt{G-EP:}\quad (\bm u^{k}_G,\alpha^{k}_G,\bm\lambda^{k}_C,\bm u_\Gamma^{k})
	=\texttt{G}(\bm u^{k-1}_G,\alpha^{k-1}_G,\bm u_\Gamma^{k,\frac{1}{2}},{\bf \Lambda}^{k}_G, \KIA_G)\ .
\end{equation}
The new coupling formulations for ductile fracture are related to the set of $(\KIA_G,\KIA_L)$, which represent the global and local augmented stiffness matrices and are given by
\begin{equation}\label{eq_aug_stiff}
		\KIA_G:=\KIA_G(\bm u_L,\alpha_L,d_L)={\bm L^T_L}{\bm T^{-T}_L}{\bm {\mathcal{S}}}_L \quad \mbox{and} \quad 
		\KIA_L:=\KIA_L(\bm u_G,\alpha_G)={\bm {\mathcal{S}}}_C\ ,
\end{equation}
which serve as augmented stiffness matrices to regularize the Jacobian matrix. Herein, ${\bf{L}}_L$ and ${\bf{T}}_L$ stand for the coupling terms which arise form the discretization of \req{RBC4}. Additionally, ${\bm{\mathcal{S}}}$ refers to {the} \textit{Steklov-Poincar\'e mapping}  \cite{Steklov,DtN1,DtN2}, which, in turn, returns the outward normal stress derivative with respect to the trace of the displacement. For details on the derivation of those matrices, we refer the interested reader to \cite{noii2020adaptive}.

\sectpc[Section361]{Global domain ($\texttt{G-EP}$) with two local domains ($\texttt{L1-EP}+\texttt{L2-D}$) }

\begin{figure}[!t]
	\centering
	{\includegraphics[clip,trim=0cm 14cm 1cm 6cm, width=16.5cm]{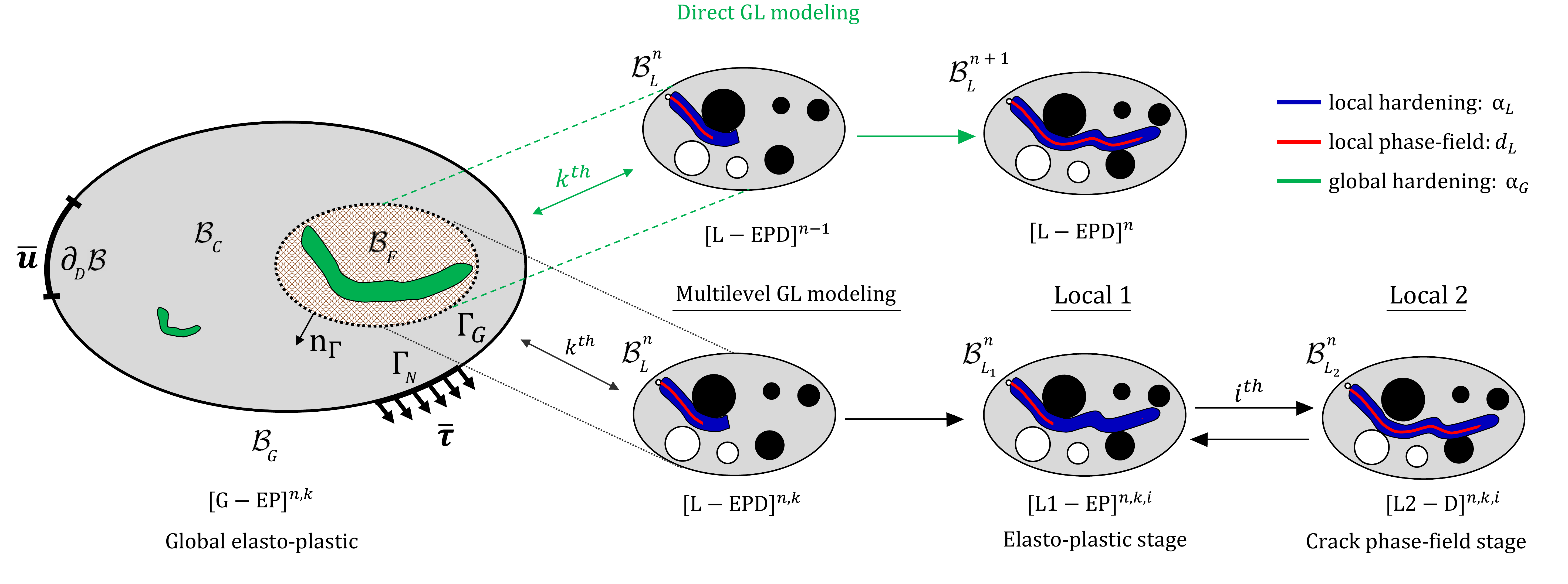}}  
	\caption{Different Global-Local scenarios. First row is the direct Global-Local approach with a single local domain, and second row is the multilevel Global-Local techniques with one domain for plasticity response along with another local domain for crack phase-field resolution.
	}
	\label{fig_multilevel_gl}
\end{figure}

Next, we focus on the key goal of this contribution, by describing the second Global-Local  formulation denoted as $g/l-2$. The main objective of this extension is to introduce an adoption of the Global-Local approach toward the multilevel local setting.  More precisely, in this setting, a global constitutive model behave as a elastic-plastic response (abbreviated as $EP$), which is augmented with \textit{two distinct local domains}; see Fig. \ref{fig_multilevel_gl}. The first local domain behaves as an elastic-plastic material ($\texttt{L1-EP}$) and the next local domain is responsible for the \textit{only} crack phase-field formulation ($\texttt{L2-D}$). 

In fact, additive splitting of the single local response ($\texttt{L: EPD}$) to the two distinct local models ($\texttt{L1-EP}+\texttt{L2-D}$), highlights the role of alternate minimization approach (i.e., iterative staggered modeling) for solving ductile phase-field fracture. Thus, if a ductile phase-field fracture is solved through the monolithic  scheme (e.g.,  a primal-dual active set method \cite{NoiiWick2019,Wick15Adapt}, to name a few among others), one should use $g/l-1$ rather than $g/l-2$.

\begin{figure}[!b]
	\centering
	{\includegraphics[clip,trim=0cm 14cm 5cm 2cm, width=14cm]{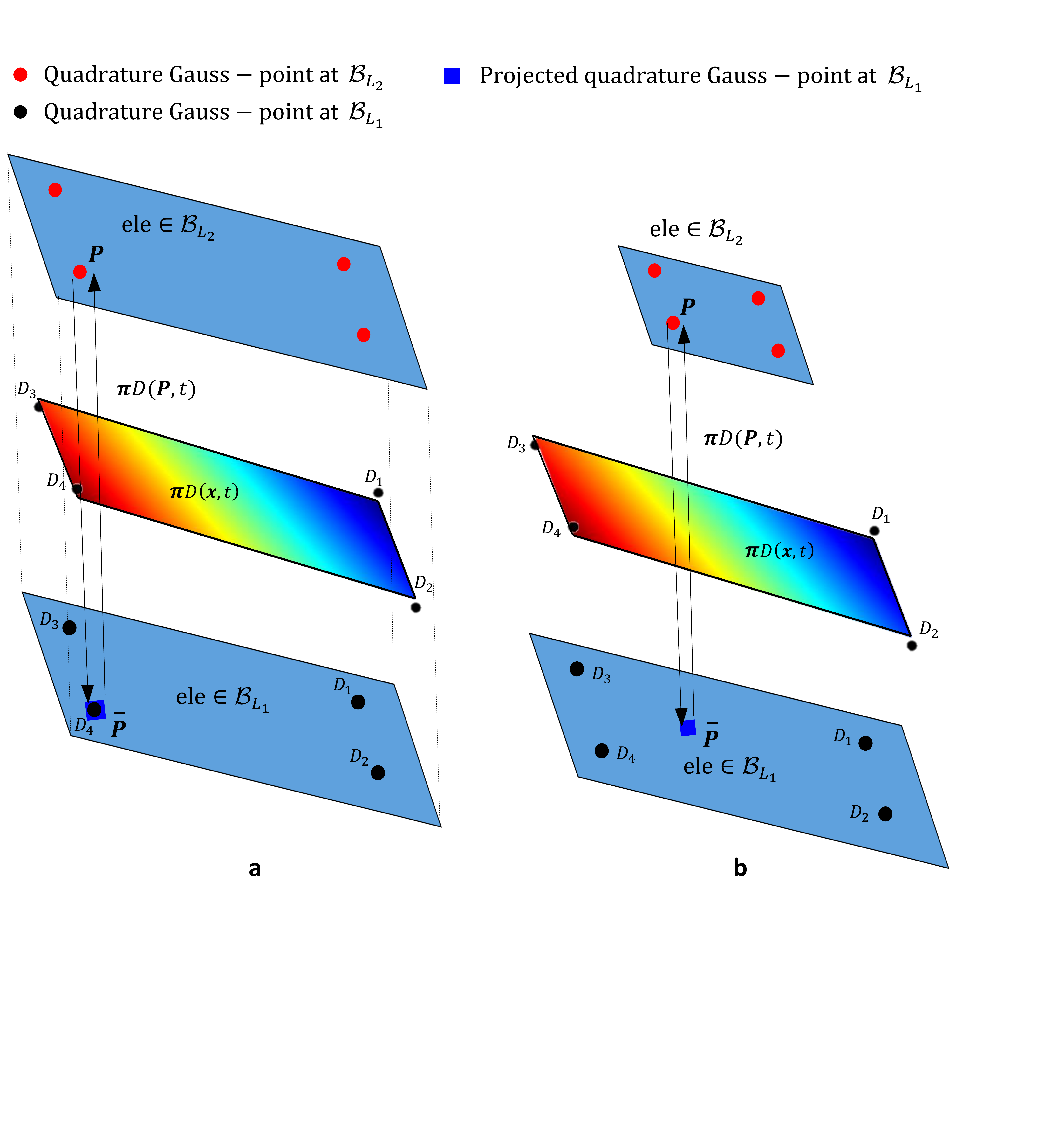}}  
	\caption{Different $g/l-2$ scenarios. (a) One-to-one discretization between $B_{L_1}$ and $B_{L_2}$, and (b) distinct discretization space between $B_{L_1}$ and $B_{L_2}$ such that $B_{L_2}
		\subset B_{L_1}$. The procedure is used to down-scale location $\Bx_{L_2}$ at point $\BP\in\calB_{L_2}$ and next to up-scale crack driving state function $D$ to solve crack phase-field equation. 
	}
	\label{fig_adaptivity_gl2}
\end{figure}

Let us define open and bounded local domains $\calB_{L_1}$ and $\calB_{L_2}$, such that
\begin{equation}\label{cond_2local}
	\calB_{L_2}\subset\calB_{L_1}\subset\calB_{G}   \quad \Leftrightarrow   \quad  l^d_L\le l^p_L\ ,
\end{equation}
where $l^p$ depicts the length-scale of ductility zone which reflects the width of plastic shear bands for hardening/softening response. Also, $l^d_L\equiv l$ in $\calB_{L_2}$ represents the local fracture length-scale. This assumption is almost valid for wide range of materials. So far, we formulated the Global-Local approach through the coupling between a global domain with a single local domain in the following abstract form
\begin{equation}\label{GL1BVP}
	{\bm \BfrakP_n}=\texttt{GL}({\bm \BfrakP_{n-1}})\ .
\end{equation}
As follows, we decompose the single local response to two local responses, resulting in the modified abstract form as
\begin{equation}\label{GL2BVP}
	{\bm \BfrakP_n}=\texttt{GL}({\bm \BfrakP_{n-1}})=
	\texttt{G}({\bm \BfrakP_{n-1}})\widetilde{\texttt{L}}({\bm \BfrakP_{n-1}})
	\WITH
	\widetilde{\texttt{L}}({\bm \BfrakP_{n-1}}):=\texttt{L}_1({\bm \BfrakP_{n-1}})+
	\texttt{L}_2({\bm \BfrakP_{n-1}})\ .
\end{equation}
In addition to that, different solver and discretization space could be used for each local domain, individually. Next, we describe in details the BVPs for each local domain, thus formulating multilevel Global-Local techniques.
%
\sectpd[Section361]{L1-EP: Boundary value problem for $\calB_{L_1}$}
The first local domain behaves as an elastic-plastic material, thus we refer it as 	\texttt{L1-EP}. Hereby, the {\it first} variational equation of \req{alLocL} is aimed to be solve, whereas projected crack phase-field denoted as $\widehat{d}_L$ from {\it second} local domain \texttt{L2-D} enters this equation as follows
\begin{equation*}\label{alLocL1_eff}
	\widetilde{\mathcal E}_{\bm u_L}(\BfrakP;\delta {\bm u}_L):=\displaystyle \int_{\calB_{L_1}}\bm\sigma(\bm u_L,\widehat{d}_L,\alpha_L):\bm\varepsilon(\delta {\bm u}_L)\,\mathrm{d}{v}
	-\int_{\Gamma_L} \bm\lambda_L \cdot \delta {\bm u}_L\,\mathrm{d}a=0\ ,
	\tag{L$_1$}
\end{equation*}
with 
\begin{equation}
	\widehat{d}(\Bx_{L_1},t) \equiv\left\{
	\begin{tabular}{ll}
		$d_{L_2}(\widehat{\Bx}_{L_1},t) $  & if ${\widehat{\Bx}_{L_1}}\in\calB_{L_2}$\ , \\[0.1cm]
		$0$ & if ${\widehat{\Bx}_{L_1}}\notin\calB_{L_2}$\ ,
	\end{tabular}
	\right.
	\label{dL-DD2}
\end{equation}
here $\widehat{\Bx}_{L_1}$ is the nearest point in $\calB_{L_2}$ which is obtained by the projection of $\bar{\BP}: \calB_{L_1}\rightarrow\calB_{L_2}$; see Fig. \ref{fig_adaptivity_gl2}. In case of \textit{no} nearest nodal point that corresponds to $\bar{\BP}$ in $\calB_{L_2}$ (i.e., $\nexists \:\BP\in\calB_{L_2}: \Vert \BP-\bar{\BP}\Vert_2<\texttt{TOL}_{Proj}= 10^{-6}$), we set $	\widehat{d}(\Bx_{L_1},t)=0$; see Fig. \ref{fig_adaptivity_gl2}. 

Equation($\text{L}_1$) together with ($\widetilde {\text{C}}_1$), ($\widetilde {\text{C}}_2$)  and evolution equations for the plastic variables \req{plast-evol-eqsL} and plasticty Kuhn-Tucker conditions \req{kktP} define the BVP for \texttt{L1-EP}. Note that, the plastic yield function and the deviatoric plastic driving force are modified through
\begin{equation}
	\chi^p_L=\chi^p(\Bu_L,\alpha_L,\widehat{d}_L) \AND	
	\BF^p_L=\BF^p(\Bu_L,\alpha_L,\widehat{d}_L)\ .
\end{equation}
Next the BVP for $\calB_{L_1}$ that deals with elastic-plastic response can be rewritten in the following abstract form:
\begin{equation}\label{BVPL1}
	\texttt{L1-EP:}\quad (\bm u^{k,i}_L,\alpha^{k,i}_L,\bm\lambda^{k,i}_L,\bm u_\Gamma^{k,i,\frac{1}{2}})=\texttt{L}_1(\widehat{d}^{k,i-1}_L,{\bf \Lambda}^{k}_L, \KIA_L,\bm u_\Gamma^{k,\frac{1}{2}})\ ,
\end{equation}
with $d^{k,0}_L=d^{k-1}_L$ and the index $i$ corresponds to the iteration process between the two local domains, see Fig. \ref{fig_multilevel_gl}.
%
\sectpd[Section361]{L2-D: Boundary value problem for $\calB_{L_2}$}
The second local domain corresponds to the fracturing response, thus we refer it as \texttt{L2-D}. Here, the \textit{second} variational equation in \req{alLocL} is aimed to be solve, whereas projected crack driving force denoted as $\widehat{\calH}(\Bx,s) $ from \texttt{L1-EP} enters this equation as follows
\begin{equation*}\label{alLocL2_eff}
	\begin{aligned}
		\displaystyle \widetilde{\mathcal E}_{d_L}(\BfrakP;\delta d_L)&=
		\int_{\calB_{L_2}} \Big( d_L + \frac{\eta_f}{\Delta t} (d_L - d^{n}_L) + (d_L-1) \;\widehat{\calH}(\bm u_L,\alpha_L) \Big) \delta d_L \;\mathrm{d}{v}\\[0.1cm]
		 \displaystyle &+ \int_{\calB_{L_2}} l^2\; \nabla d_L.\nabla(\delta d_L)\,\mathrm{d}{v}= 0\ ,
    \end{aligned}
	\tag{L$_2$}
\end{equation*}
such that 
\begin{equation}
	\widehat{\calH} = \max_{s\in [0,t]} \widehat{D}(\Bx,s) \ge 0 \WITH
	\widehat{D}(\Bx_{L_2},t)\equiv\pi{D}(\Bx_{L_2},t)\ .
	\label{driving-force2}
\end{equation}
Here , we defined the linear interpolation surface operator $\pi{D}: \text{ele}\in\calB_{L_1}\rightarrow\calB_{L_2} $ for the given crack driving state function. To clarify this, let us consider Fig. \ref{fig_adaptivity_gl2}. We aim to determine $\widehat{D}(\BP)\in\calB_{L_2}$, such that $\BP=(x=9,y=24.9)\in\calB_{L_2}$. To do so, we assume the Cartesian
coordinates of the quadrature Guass-points for the set of crack driving state function $(D_1=22738,D_2=24033,D_3=42292,D_4=50288) \in ele\in\calB_{L_1} $ are
\begin{equation*}
	\begin{aligned}
		(x_{D_1}, x_{D_2}, x_{D_3}, x_{D_4}) &= (9.064,9.064,8.936,8.936),\\
		(y_{D_1}, y_{D_2}, y_{D_3}, y_{D_4}) &= (24.952, 24.82, 24.95, 24.82)\ .
    \end{aligned}
\end{equation*}
A linear surface denoted by $\pi D(\Bx,t)$ which cross over these points is obtained as,
\begin{equation}\label{linear_surface}
	\pi D(\Bx,t)=c_0+c_1x_1+c_2x_2\  ,
\end{equation}
with
\begin{equation}
(c_0,c_2,c_2)=(2.523\times10^
	6,-1.785\times10^5, -3.541 \times10^4)\ .
\end{equation}
Thus, the crack driving state function at point $\BP$ is approximately obtained as $	\widehat{D}\equiv\pi{D}(\BP,t)=34355$; as sketched in Fig. \ref{fig_adaptivity_gl2}.
\begin{Remark}
\label{rem_adpat_higherorder}
A further improvement of the proposed framework can be done to reduce the error of  interpolation/extrapolation during the determination of the crack driving state function, i.e., $\widehat{D}(\Bx_{L_2},t)$. This can be achieved  by the introduction of higher-order Global-Local approach through the following producers:
	\begin{itemize}
		\item[(1)]  One could use the information of neighboring quadrature Gauss-points (within nearest elements). Thus we have more information to create a higher-order surface for $\pi D(\Bx,t)$ instead of linear one (which is used in \req{linear_surface}). 
		\item[(2)] One could use higher number of quadrature Gauss-points per elements. Hence, we have more information for the $D(\Bx,t)$, resulting with higher-order surface for $\pi D(\Bx,t)$. 		
	\end{itemize}
	{Hence, these type of enhancement for the Global-Local procedure is open for further research.}
\end{Remark}
Equation ($\text{L}_2$) along with the Kuhn-Tucker conditions for the phase-field equation, i.e., \req{kktD} define the BVP for \texttt{L2-D}. Next, the BVP for $\calB_{L2}$ that deals with fracturing response can be rewritten in the following abstract form:
\begin{equation}\label{BVPL2}
	\texttt{L2-D:}\quad d^{k,i}_L=\texttt{L}_2(\widehat{D}^{k,i},d^{n}_L)\ .
\end{equation}
\newpage
The aforementioned Global-Local procedure is summarized in Algorithm \ref{alg_3}.
\begin{algorithm}[H]\small
	\caption{\em Multilevel Global-Local procedure with two distinct local domains.}
	\label{alg_3}
	Let assume two distinct local domains $\calB_{L_1}$ and $\calB_{L_2}$ . We solve $\texttt{L1-EP}$ and $\texttt{L2-D2}$ through the projected $\widehat{d}$ and $\widehat{D}$, respectively, along with $\texttt{G-EP}$, by the following steps:\\[0.1cm]
	{\bf{Initialization:}} We set initial guess through:\\[0.3cm]
	\quad$(\bm u^{n,0}_G,\alpha^{n,0}_G,\bm\lambda^{n,0}_C,\bm u^{n,0}_L,\alpha^{n,0}_L,\bm\lambda^{n,0}_L,\bm u_\Gamma^{n,0}):=(\bm u^{n-1}_G,\alpha^{n-1}_G,\bm\lambda^{n-1}_C,\bm u^{n-1}_L,\alpha^{n-1}_L,\bm\lambda^{n-1}_L,\bm u_\Gamma^{n-1})$.\\[0.3cm]
   \textbf{Loop.} Increment for $k=1$ and $i=1$ until convergence:
	\begin{enumerate}
		\item[\textbf{1.1}.]Solve local 1 through $\texttt{L1-EP:}\quad (\bm u^{k,i}_L,\alpha^{k,i}_L,\bm\lambda^{k,i}_L,\bm u_\Gamma^{k,i,\frac{1}{2}})=\texttt{L}_1(\widehat{d}^{k,i-1}_L,{\bf \Lambda}^{k}_L, \KIA_L,\bm u_\Gamma^{k,i-1,\frac{1}{2}})^{\star}$,
		\item[\textbf{1.2}.] Post-processing step for $\texttt{L1-EP}$: For every quadrature Guass-points located at $\Bx_{L_2}$ in $\calB_{L_2}$ find interpolated/extrapolated crack driving state function denoted as $\widehat{D}^{k,i}$,
		\item[\textbf{1.3}.] Solve local 2 through $\texttt{L2-D:}\quad d^{k,i}_L=\texttt{L}_2(\widehat{D}^{k,i})$,
		\item[\textbf{1.4}.] Post-processing step for $\texttt{L2-D}$: For every nodal points located at $\Bx_{L_1}\in\calB_{L_1}$ find nearest node in $\calB_{L_2}$ (if exist) to determine $\widehat{d}$ ( and if it does not exist, we set $\widehat{d}=0$),
		\item[\textbf{1.5}.]  If \textit{checking criterion}$^{\star\star}$ between local 1 and local 2 is satisfied, thus set\\[0.3cm]
		\hspace{2.7cm}  $(\bm u^{k,i}_L,d^{k,i}_L,\alpha^{k,i}_L,\bm\lambda^{k,i}_L,\bm u_\Gamma^{k,i,\frac{1}{2}}):=
		(\bm u^{k}_L,d^{k}_L,\alpha^{k}_L$ $,\bm\lambda^{k}_L,\bm u_\Gamma^{k,\frac{1}{2}})$,\\[0.3cm]
		\textbf{STOP} and \textbf{GO} to \textbf{2.1}; else increment $i+1\rightarrow i$, and \textbf{GO} to \textbf{1.1},
		\item[\textbf{2.1}.] Solve global BVP through $\texttt{G-EP:}$ $(\bm u^{k}_G,\alpha^{k}_G,\bm\lambda^{k}_C,\bm u_\Gamma^{k})
		=\texttt{G}(\bm u^{k-1}_G,\alpha^{k-1}_G,\bm u_\Gamma^{k,\frac{1}{2}},{\bf \Lambda}^{k}_G, \KIA_G){^\star}$,
		\item[\textbf{2.2}.] If Global-Local procedure converged$^{\star\star\star}$, set\\[0.3cm]
		$(\bm u^{k}_G,\alpha^{k}_G,\bm\lambda^{k}_C,\bm u^{k}_L,\alpha^{k}_L,\bm\lambda^{k}_L,\bm u_\Gamma^{k})=:(\bm u^{n+1}_G,\alpha^{n+1}_G,\bm\lambda^{n+1}_C,\bm u^{n+1}_L,\alpha^{n+1}_L,\bm\lambda^{n+1}_L,\bm u_\Gamma^{n+1})$,\\[0.3cm]
		\textbf{STOP} and \textbf{Go} to \textbf{3};
		else increment $k+1\rightarrow k$, and \textbf{GO} to \textbf{1.1}.
		\item[\textbf{3}.] {\bf Output:} Solution $(\bm u^{n+1}_G,\alpha^{n+1}_G,\bm\lambda^{n+1}_C,\bm u^{n+1}_L,\alpha^{n+1}_L,\bm\lambda^{n+1}_L,\bm u_\Gamma^{n+1},\calH^{n+1}_{L})$.
	\end{enumerate}
	{\parbox{6in}{ {$\star$} index $i$ corresponds to the iteration process between $(\calB_{L_1}, \calB_{L_2})$ and index $k$ corresponds to the Global-Local iterations; see Fig. \ref{fig_multilevel_gl}.}}\\[0.15cm]
	{\parbox{6in}{ {$\star\star$} $\sqrt{\Vert \Bu^{k,i}_{L}-\Bu^{k,i-1}_{L}\Vert_2+\Vert d^{k,i}_{L}-d^{k,i-1}_{L}\Vert_2}<\texttt{TOL}_{stag}:=10^{-4},$}}\\[0.15cm]
	{\parbox{6in}{ {{$\star\star\star$} $\sqrt{ {\left\| \bm u_\Gamma^k-\bm u_L^k \right\|^2_{{\bf L}^2(\Gamma)}} + {\left\| \bm\lambda_F^n-\bm\lambda_F^{k+1} \right\|^2_{{\bf L}^2(\Gamma)}} }<\texttt{TOL}_{\text{GL}}:=10^{-6}$; see \cite{NoiiGL18} Section 3.4.2.}}}\\[0.15cm]
\end{algorithm}

\sectpa[Section4]{Predictor-Corrector Adaptivity Applied to the Global-Local
	Formulation of Ductile Phase-Field Fracture}

To further reduce the computational time, an adaptive Global-Local approach is used.
To this end, a predictor-Corrector concept  is preformed at $t_{n}$  \cite{noii2020adaptive}. Now, let the following prerequisites without loss of generality are taken into account:

\begin{equation}
	(1): 0\ge\dot \chi^p_G\ge\dot\chi^p_L\ge\dot\chi^d_L \AND	(2):  l^p_L\ge l^d_L\ge0\ .
\end{equation}
\begin{itemize}
	\item[(1)] The prior condition results in the global plastic evolution yield earlier than local plastic surface, thus leads us to the prior knowledge of those global elements which have to be refined. In other words, this prior global knowledge provides us the refinement strategy for our Global-Local approach.
	\item[(2)] Later assumption results that plasticity initiated before fracture state, i.e. $\varepsilon_d:=\psi_c/Y0+Y_0/{2E}>\varepsilon_p:=Y_0/E$, in which $\varepsilon_d$ and $\varepsilon_p$ refer to the yield fracture strain and yield plastic strain, respectively. In other words, this hypothesis results from the fact that fracture behaves in a more localized region compare to the plasticity zone thus it is bounded in the plasticity state. In fact, the earlier condition $\calB_{L_2}\subset\calB_{L_1}$ that is presented in \req{cond_2local} is the result of this assumption. Note, this assumption is valid for wide range of materials. For a detail discussion on the different elastic-plastic-fracture scenarios see \cite{alessi2018coupling}.
\end{itemize}

\begin{figure}[!ht]
	\centering
	{\includegraphics[clip,trim=0cm 14cm 0cm 1.8cm, width=16cm]{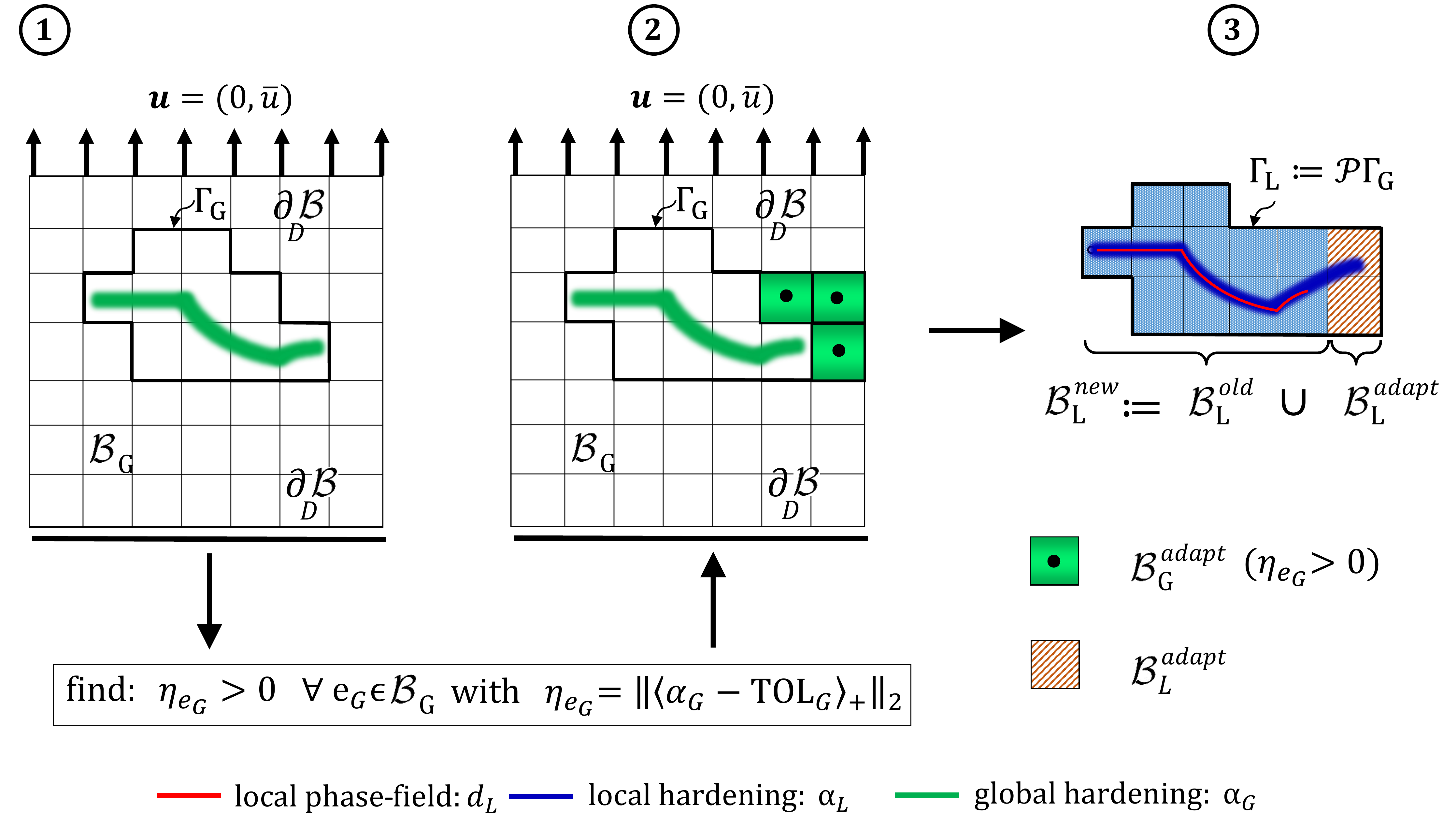}}  
	\caption{Explanation of flagging cells treatment for the adaptive refinement of the ductile phase-field fracture according to $\eta_{e_G}$. 
	}
	\label{adaptivity2}
\end{figure}

Let us now assume the Global-Local formulation is at the converged state, thus results in ${\BfrakP^0_{n}}:=\texttt{GL}({\BfrakP_{n-1}})$. The
Global-Local approach is augmented by a {\it dynamic allocation} of a local state using an adaptive scheme which has to be performed at
time step $t_{n}$. By the adaptivity procedure, we  mean: ({a}) to
determine which global elements need to be refined and identified by
$\calB^{\;adapt}_G\subset\calB_G$, see Fig. \ref{adaptivity2} ;({b}) to create the new fictitious
domain $\calB^{\;new}_F:=\calB^{\;old}_F\cup\calB^{\;adapt}_F$ with $\calB^{\;adapt}_F:= \calB^{\;adapt}_G$ and result in
a new local domain that is $\calB^{\;new}_L:=\calB^{\;old}_L\cup\calB^{\;adapt}_L$; ({c}) to determine a new local interface denoted as $\Gamma_L$;
({d}) to interpolate the old global solution in $\calB^{\;adapt}_L$.  

Hereby, different options are devised for the adaptivity procedure of ductile phase-field fracture.

\begin{itemize}
	\item[\textbf{Option 1}.] First, we assume a global domain $\texttt{G-EP}$ is coupled with a single local domain $\texttt{L-EPD}$. Here, the adaptivity procedure is described based on crack phase-field, which, in turn, $\texttt{TOL}_d$ (denoted as a crack phase-field threshold value) is used to determine the adaptivity procedure;  \cite{noii2020adaptive}. In this case, we observe the great continuity between displacement fields (i.e., $\Bu_G$ and $\Bu_L$ due to (C$_2$) and (C$_3$)) while discontinuity for hardening fields (i.e., $\alpha_G$ and $\alpha_L$ ) drastically increased. This is mainly because hardening fields are defined as an internal field and not as a primary fields (i.e., hardening value lives at the quadrature Guass-points). This significant discontinuity between $\alpha_G$ and $\alpha_L$ cause two problems: (1) serious convergence issues in our numerical treatment., and (2) the plasticity/fracture path could be lost (due to the lack of continuity between two scales). This kind of issue also reported for the computational homogenzation applied to the localaized problems, whereas local scale as a RVE imposed to the Global  quadrature Guass-points thus behaves as an internal state; see \cite{Fish2014,zohdi2008introduction}. 
	\item[\textbf{Option 2}.] Next, we consider a global domain $\texttt{G-EP}$ which is coupled with a single local domain $\texttt{L-EPD}$. Here, the adaptivity procedure is explained based on global hardening value $\alpha_G$ to determine the adaptivity procedure. To do so, we define the following adaptivity indicator for every global element:
	\begin{equation}\label{eta_G}
		\text{find:}\;\; \eta_{e_G}>0\;\; \text{with}\;\;\eta_{e_G}:=\Vert\langle\dot\alpha_G-\texttt{TOL}_{\alpha}\rangle_+\Vert_2
		= \Vert\langle\lambda_G^p-\texttt{TOL}_{\alpha}\rangle_+\Vert_2
		\quad\forall e_G\in\calB_G\ ,
	\end{equation}
	with $\texttt{TOL}_{\alpha}$  to be prescribed. It can be grasped that indicator $\eta_{e_G}$ manifests those elements whose are living in the elasticity or plasticity states. As a result, local domains are defined for those elements if $\eta_{e_G}>0$. 
	Typically, we set $\texttt{TOL}_{\alpha}:=10^{-4}$. The detailed adaptivity algorithm is sketched in Table \ref{alg_4}.
	\item[\textbf{Option 3}.] Finally, we assume a global domain $\texttt{G-EP}$ which is coupled with two distinct local domains as $\texttt{L1-EP}$ and $\texttt{L2-D}$. In this case, we have \textit{two-way} of adaptivity procedure. Specifically, we have one adaptivity step for plasticity state in which a refinement knowledge coming from global solution state. The second  adaptivity step deals with fracture state such that refinement knowledge is the results of the lowest local level, i.e., local 2. For a better insight into this scenario,  two-way multi-level adaptivity scheme, is summarized through the following steps:
    \begin{itemize}
    	\item[(1).] Once $\dot\chi^p_G=0$ holds, a global domain is reached to the plastic yield surface. Thus, we perform Algorithm \ref{alg_4}, to determine \texttt{L1-EP}. We set resulting new global elements which has to be refined as $ele^{n}_{Gp}$ at time step $n$. Next we proceed steps 1-2  from Algorithm \ref{alg_3}.
        \item[(2).] As soon as $\dot\chi^d_L=0$ holds, local 1 reached to the fracture yield state. If $\calB_{L_2}\subset\calB_{L_1}$, we perform Algorithm \ref{alg_5},  to determine \texttt{L2-D}. We set resulting new global elements which has to be refined as $ele^{n}_{Gd}$ at time step $n$.
    	\item[(3).] If $\calB_{L_2}\not\subset\calB_{L_1}$ holds, thus $\exists\; ele^{n}_{Gd}\not\subset\calB_{L_1}$. To do so, for those element we set positive value for $\eta_G$ and go back to step 1. As a result, we will have $\calB_{L_2}\subset\calB_{L_1}$. In other words, for those global elements which are not already refined for $\calB_{L_1}$ but needs to be considered for the $\calB_{L_2}$, we first refined it for $\calB_{L_1}$ and then continue for solving BVP at $\calB_{L_2}$.
    	\item[(4).] Finally, we introduce a corrector step {in which the computation is re-run on the newly determined local mesh.} To this end, we compute Global-Local
    	framework in Algorithm \ref{alg_3} with new interpolated solutions (output from Algorithms \ref{alg_4}-\ref{alg_5}). 
    	\begin{equation*}
    		\text{Corrector step: Compute the Global-Local solution by } {\BfrakP_{n}}=\texttt{GL}({\BfrakP^*_{n}})\;\;\text{in} \;\calB^{new}_{L}.
    	\end{equation*}
    \end{itemize}
\end{itemize} 
In the following, if we are dealing with $g/l-1$ an option 2 is used, while option 3 is applied for $g/l-2$.  Note that, option 1 is not considered in this work.
 \begin{algorithm}[H]\small
 	\caption{\em Predictor steps for the adaptive procedure: plasticity.}	
 	\label{alg_4}
 	Let $0<\texttt{TOL}_{\alpha}<1$ be given. At the fixed time step $n$, predictor-corrector steps for the adaptive procedure is performed using the following steps:
 	\begin{enumerate}
 		\item Compute the Global-Local solution from Algorithm \ref{alg_3} and set ${\BfrakP^0_{n}}=\texttt{GL}({\BfrakP_{n-1}})$,
 		\item Compute $\eta_{e_G}:=\Vert\langle\dot\alpha_G-\texttt{TOL}_{\alpha}\rangle_+\Vert_2$ at every $e_G\in\calB_G$,
 		\item For every $e_G\in\calB_G$  if $\eta_{e_G}>0$ set adaptivity flag 1 otherwise 0,
 		\item Set new local adaptive domain denoted by $\calB^{adapt}_L$ by adding refined global $e_G$ solution results from previous step. We denote interpolated solution in the $\calB^{adapt}_L$ with $\BfrakP^{adapt}_{n}$,
 		\item We define ${\BfrakP^*_{n}}:={\BfrakP^0_{n-1}}$ if $\Bx\in\calB^{old}_L$ otherwise  ${\BfrakP^*_{n}}:={\BfrakP^{adapt}_{n}}$ if $\Bx\in\calB^{adapt}_L$.
 	\end{enumerate}
 Output: $\calB^{new}_L=\calB^{\;old}_L\cup\calB^{\;adapt}_L$ (or $\calB^{new}_{L_1}=\calB^{\;old}_{L_1}\cup\calB^{\;adapt}_{L_1}$) and $\BfrakP^*_{n}$,
 \end{algorithm}

 \begin{algorithm}[H]\small
 	\caption{\em Predictor steps for the adaptive procedure: crack phase-field.}	
 	\label{alg_5}
 	Let $0<\texttt{TOL}_{d}<1$ be given. At the fixed time step $n$, predictor-corrector steps for the adaptive procedure is performed using the following steps:
 	\begin{enumerate}
 		\item Compute the Global-Local solution from Algorithm \ref{alg_3} and set ${\BfrakP^0_{n}}=\texttt{GL}({\BfrakP_{n-1}})$,
 		\item Find $\bm x_L\in\Gamma_L:$ such that $d_L({\bm
 			x}_L)<\texttt{TOL}_{d}$ on $\Gamma_{L_2}$,
 		\item Find $E^{q}_{L,i}\in\Gamma_L$ such that ${\bm x_L}\in E^{q}_{L,i}$,
 		\item Find $E^{q}_{G,i}=\calP^{-1}E^{q}_{L,i}$ (corresponding edge {in} $\calB_G$),
 		\item Find $e_{G}\in\calB_G$ and $e_{G}\not\in\calB_F$ such that $E^{q}_{G,i}\in e_{G}$,
 		\item Set new local adaptive domain denoted by $\calB^{adapt}_{L_2}$ by adding refined global $e_G$ solution results from previous step. We denote interpolated solution in the $\calB^{adapt}_{L_2}$ with $\BfrakP^{adapt}_{n}$,
 		\item We define ${\BfrakP^*_{n}}:={\BfrakP^0_{n}}$ if $\Bx\in\calB^{old}_{L_2}$ otherwise  ${\BfrakP^*_{n}}:={\BfrakP^{adapt}_{n}}$ if $\Bx\in\calB^{adapt}_{L_2}$.
 	\end{enumerate}
  Output: $\calB^{new}_{L_2}=\calB^{\;old}_{L_2}\cup\calB^{\;adapt}_{L_2}$ and $\BfrakP^*_{n}$,
 \end{algorithm}

\sectpa[Section5]{Numerical Examples}

This section demonstrates the performance of the proposed adaptive multilevel Global-Local approach within phase-field ductile fracture. Four numerical examples are investigated. The material parameters listed in Table~\ref{material-parameters} are based on \cite{ambati2015phase,miehe2016phase}. In the numerical simulation all variables for both the global and local domains, are discretized by bilinear quadrilateral $Q_1$ finite elements. An alternate minimization scheme is used for solving the local BVP indicated in (L). Thus, we alternately {solve} for $d_L$ by fixing ${\bm u}_L$ and then solving for $({\bm u}_L,{\bm u}_{\Gamma},{\bm \lambda_L},\alpha_L)$ by fixing $d_L$ until convergence is reached.  
The proposed adaptive multilevel Global-Local approach for the ductile phase-field fracture is described based on set of given threshold values. These values are explained in Table  \ref{model-parameters}. 

\setcounter{table}{0}
\begin{table}[!ht]
	\caption{Material parameters employed in the numerical experiments according to \cite{ambati2015phase,miehe2016phase}.}	\vspace{1mm}
	\centering
	\begin{tabular}{ccllll}
		No.  &Parameter & Name                   & Example 1,3,4   & Example 2  & Unit            \\[2mm]\hline 
		1.   &$\mu$        & shear modulus       & $27,280$     & $70,300$  & $\mathrm{MPa}$ \\[2mm]
		2.   &$K$      & bulk modulus       & $71,660$   & $1,36,500$  & $\mathrm{MPa}$                 \\[2mm]
		3.   &$H$        & Hardening modulus        & $250$  & $300$ & $\mathrm{MPa}$ \\[2mm]
		4.   &$Y_0$        & Yield stress    & $345$ & $443$   & $\mathrm{MPa}$ \\[2mm]
		5.   &$Y_\infty$        &  Infinite yield stress    & $345$ & $443$   & $\mathrm{MPa}$ \\[2mm]
		6.   &$\psi_c$     & Specific fracture energy  & $25$ & $25$  & $\mathrm{MPa}$ \\[2mm]
		7.  &$\eta_f$    & Crack viscosity     & $10^{-14}$ & $10^{-14}$  & $\mathrm{N/m^{2}s}$\\[2mm] 
		8.  &$\kappa$    & Stabilization parameter     & $10^{-8}$ & $10^{-8}$ & --\\[2mm]
		\hline
		\label{material-parameters}
	\end{tabular}
\end{table}

\begin{table}[!ht]
	\caption{Numerical parameters employed in the following examples.}	\vspace{1mm}
	\centering
	\begin{tabular}{ccllll}
		No.  &Parameter & Tolerance for                   & Tolerance value              \\[2mm]\hline 
		1.   &$\texttt{TOL}_{N-R}\;$        & Newton-Raphson       & $10^{-8}$      \\[2mm]
		2.   &$\texttt{TOL}_{d}\quad\;\;$     & Adaptivity of phase-field       & $0.05$                    \\[2mm]
		3.   &$\texttt{TOL}_{\alpha}\quad\;$        & Adaptivity of plasticity       & $10^{-4}$   \\[2mm]
		4.   &$\texttt{TOL}_{stag}\;$       & Alternate minimzation    & $10^{-4}$   \\[2mm]
		5.   &$\texttt{TOL}_{Proj}$        &  Nearest point    & $10^{-6}$  \\[2mm]
		6.   &$\texttt{TOL}_{\text{GL}}\;$     & Global-Local coupling  & $ 10^{-3}$  \\[2mm]
		\hline
		\label{model-parameters}
	\end{tabular}
\end{table}
The overall response of the Global-Local approach in terms of accuracy/robustness
and efficiency was verified using single-scale solutions. In detail, we investigate: 
\begin{itemize}
	\item Load-displacement curves to evaluate the up-scaling procedure (i.e., a transition of local non-linearity and imperfections to the global level);
	\item Evaluating global equivalent plastic strain $\alpha_G$ as an effective local hardening;
	\item Evolution of local phase-field patterns $d_L$ in order to evaluate the down-scaling procedure (i.e. transition of external loading from the global to the local level);
	\item Evaluating local equivalent plastic strain $\alpha_L$;
	\item The efficiency of Global-Local formulation through the total accumulated computational time.
\end{itemize}
%
\sectpb[exm1]{Example 1: Double-notched specimen under tensile loading}

\begin{figure}[!t]
	\centering
	{\includegraphics[clip,trim=1cm 2.8cm 0cm 3cm, width=11cm]{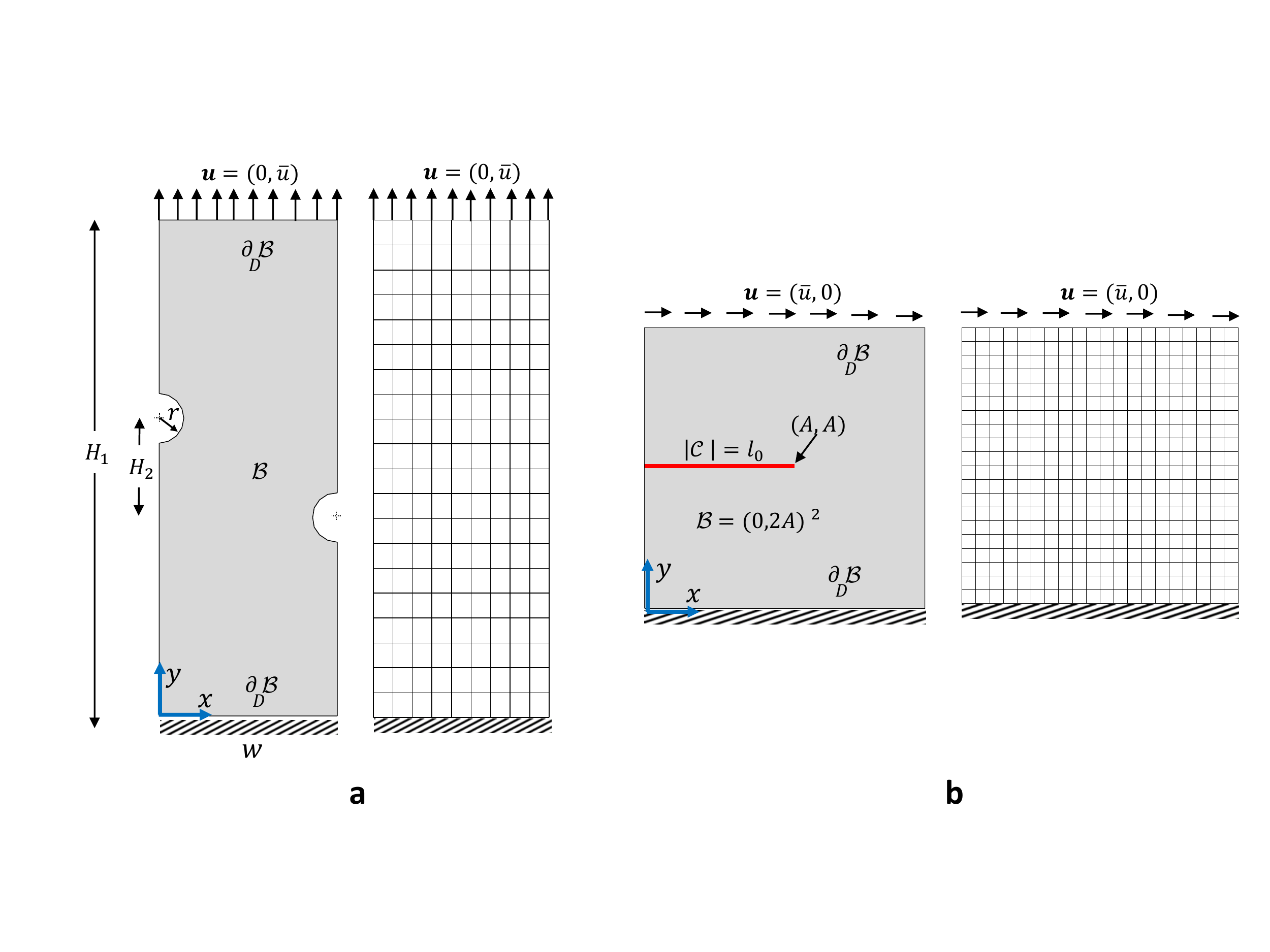}}  
	\caption{Geometry and loading setup for (a) Example 1. Double-notched specimen under tensile loading, and (b) Example 2. Single-edge-notched shear test with their global discretization.
	}
	\label{example1-2}
\end{figure}

To gain the first insight into the performance of the Global-Local approach, the following numerical example is concerned with the simulation of the double-notched specimen under tension.  The configuration is shown in Fig. \ref{example1-2}a. The top edge is constrained horizontally while the bottom edge is fixed for displacement in $x-y$ directions. The geometrical dimensions for Fig. \ref{example1-2}a are set as $H_1=50\;mm$, $H_2=5\;mm$, and $w=18\;mm$ with radius of two notches as  $r=2.5\;mm$.

\begin{figure}[!ht]
	\centering
	{\includegraphics[clip,trim=3cm 3cm 5cm 4cm, width=15cm]{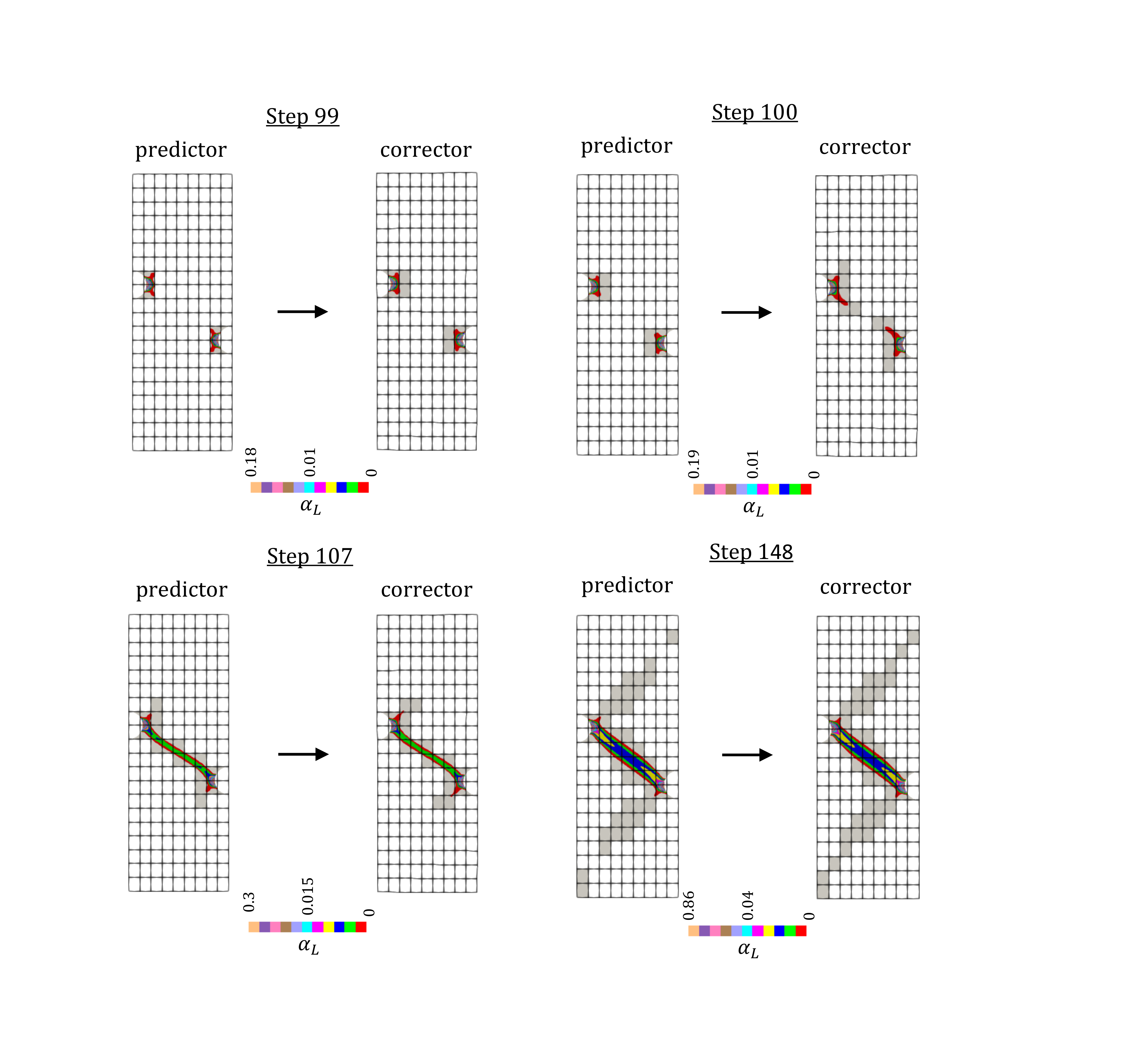}}  
	\caption{Example 1. Global-Local approach augmented with the predictor-corrector adaptive
		scheme; Evolution of the local hardening value of double-notched specimen for ductile fracture with different loading steps.}
	\label{exm1_profile_PC}
\end{figure}

A monotonic displacement increment  ${\Delta \bar{u}}_y=2\times10^{-3}\;mm$ is applied in a vertical direction at the top boundary of the specimen for 215 time steps. The minimum finite element size in the single-scale and local domains is $0.3\;mm$, which, in turn, the heuristic requirement $h<l/2$ inside the localization zone is fulfilled. The
single-scale domain partition contains 8715  elements while the global
domain contains 180 elements. Due to the adaptivity procedure, the number of elements at each time step may be vary in space discretization. The material and numerical parameters are those given in Table \ref{material-parameters} and Table \ref{model-parameters}, respectively. Accordingly, the degrees of freedom are shown in Table \ref{GL-SS-eG1}. 
In this example, a rather coarse mesh is used to highlight the effect of the predictor-corrector scheme on the adaptivity procedure; see Fig. \ref{exm1_profile_PC}.  By applying predictor-corrector steps, a better estimation for (i) the elastic-plastic response, and (ii) the fracture state before proceeding to the next time step are achieved  (since more elements are locally resolved; see Fig. \ref{exm1_profile_PC}). 
It was observed that the corrector scheme applied to the predictor step improved
the Global-Local results for the ductile fracture. 

\begin{figure}[!ht]
	\centering
	{\includegraphics[clip,trim=3.3cm 6cm 3.5cm 21cm, width=15cm]{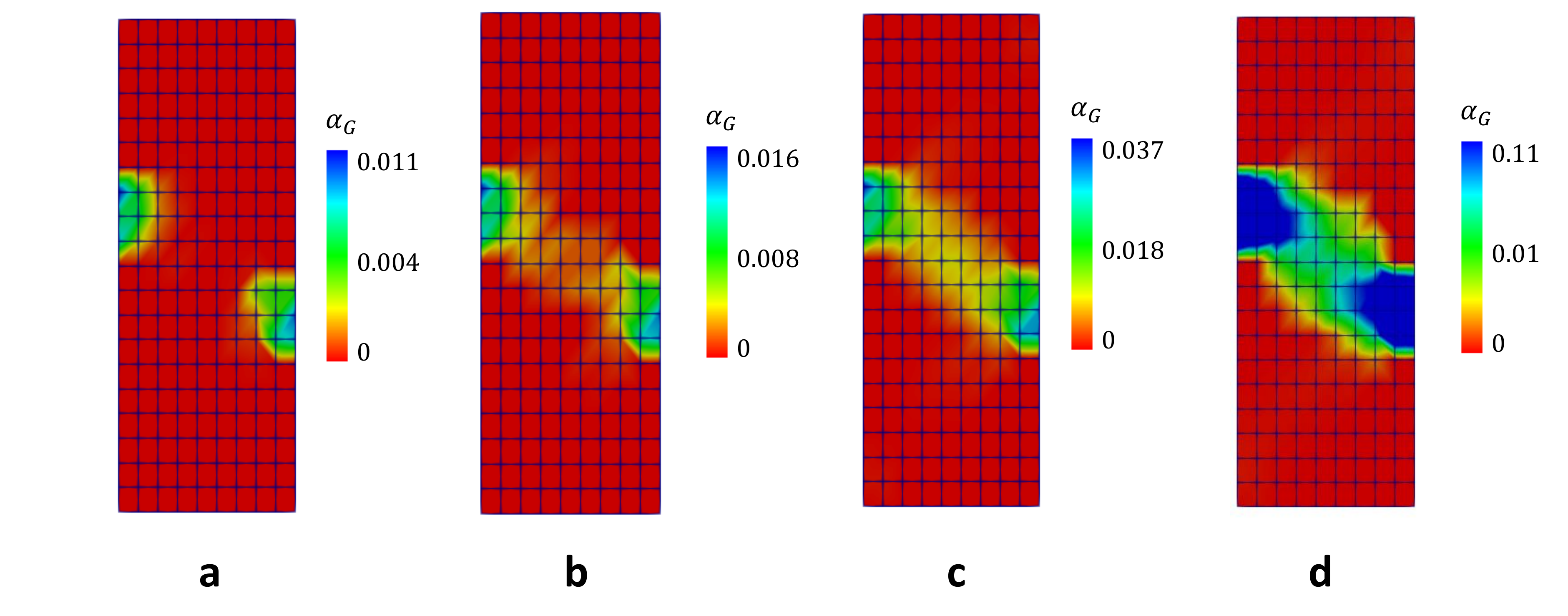}}  
	\caption{Example 1 ($g/l-1$). Evolution of the global hardening value $\alpha_G$ for different deformation stages up to
		complete failure at $\bar{u}_y=[0.192,0.212,0.278,0.42]~mm$.}
	\label{exm1_profile_gl1_hardningG}
\end{figure}

The computational analysis starts by illustrating the $g/l-1$ solutions for different deformations states up to complete failure.
The evolution of the global hardening $\alpha_G$ is demonstrated in Fig. \ref{exm1_profile_gl1_hardningG} for four-time steps at $\bar{u}_y=[0.192,0.212,0.278,0.42]~mm$. Even though there is no global imperfection, the influence of locally exists notched results in a consistent global plasticity flow with local hardening evolution $\alpha_L$.  That is mainly explained due to the consistency between the two scales. Hereby, the maximum global plasticity appears where the two local notches exist; see Fig. \ref{exm1_profile_gl1_hardningG}a. Thus, $\alpha_G$ can be interpreted as an effective hardening quantity which roots from its local source $\alpha_L$.

\begin{figure}[!ht]
	\centering
	{\includegraphics[clip,trim=5cm 6cm 4.1cm 6cm, width=14cm]{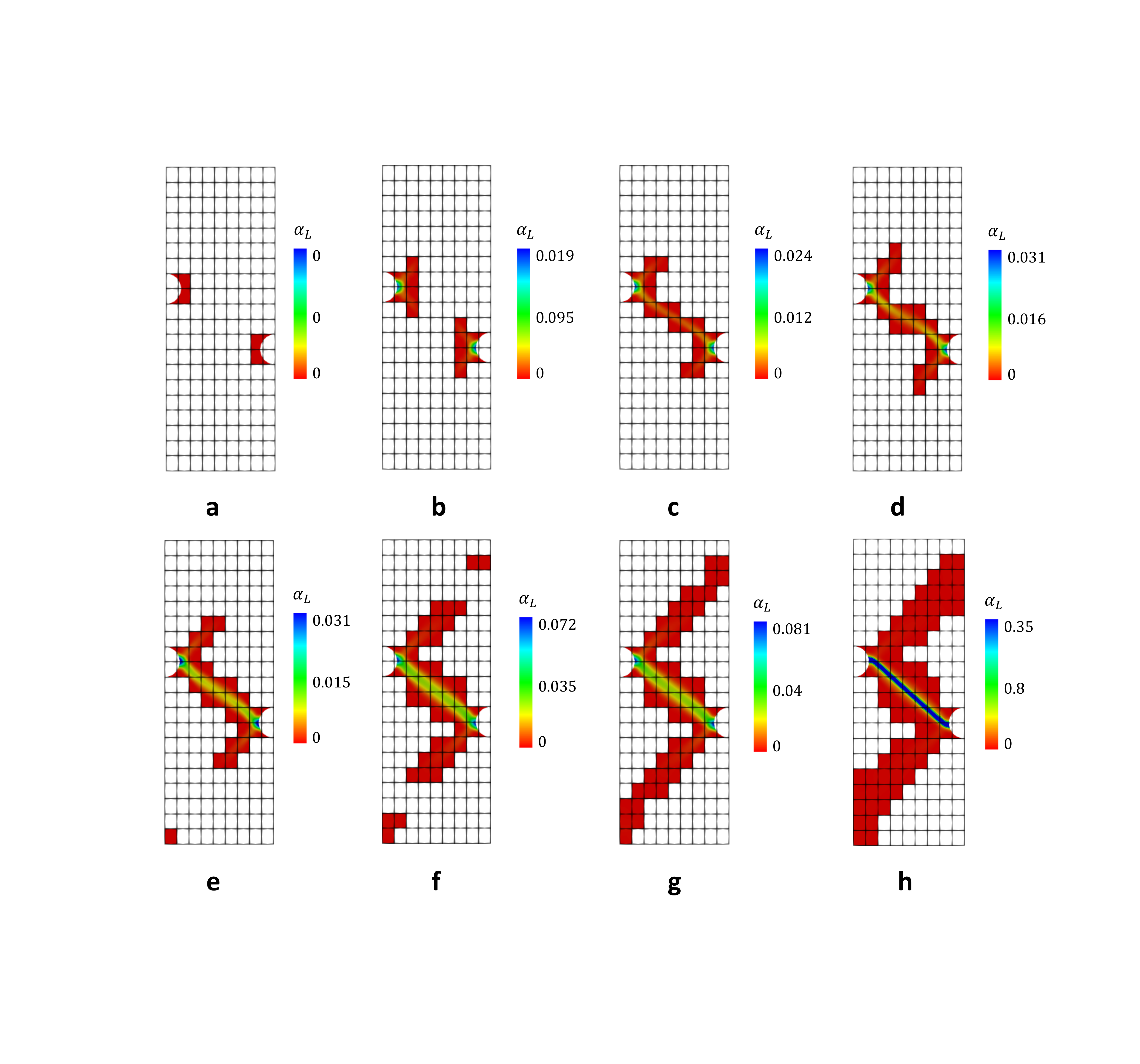}}  
	\caption{Example 1 ($g/l-1$). Evolution of the local hardening value $\alpha_L$ for different deformation stages up to
		complete failure at $\bar{u}_y=[0,0.192,0.2,0.212,0.254,0.278,0.29,0.42]~mm$.}
	\label{exm1_profile_gl1_hardningL}
\end{figure}

\begin{figure}[!ht]
	\centering
	{\includegraphics[clip,trim=4cm 21cm 5cm 5cm, width=15cm]{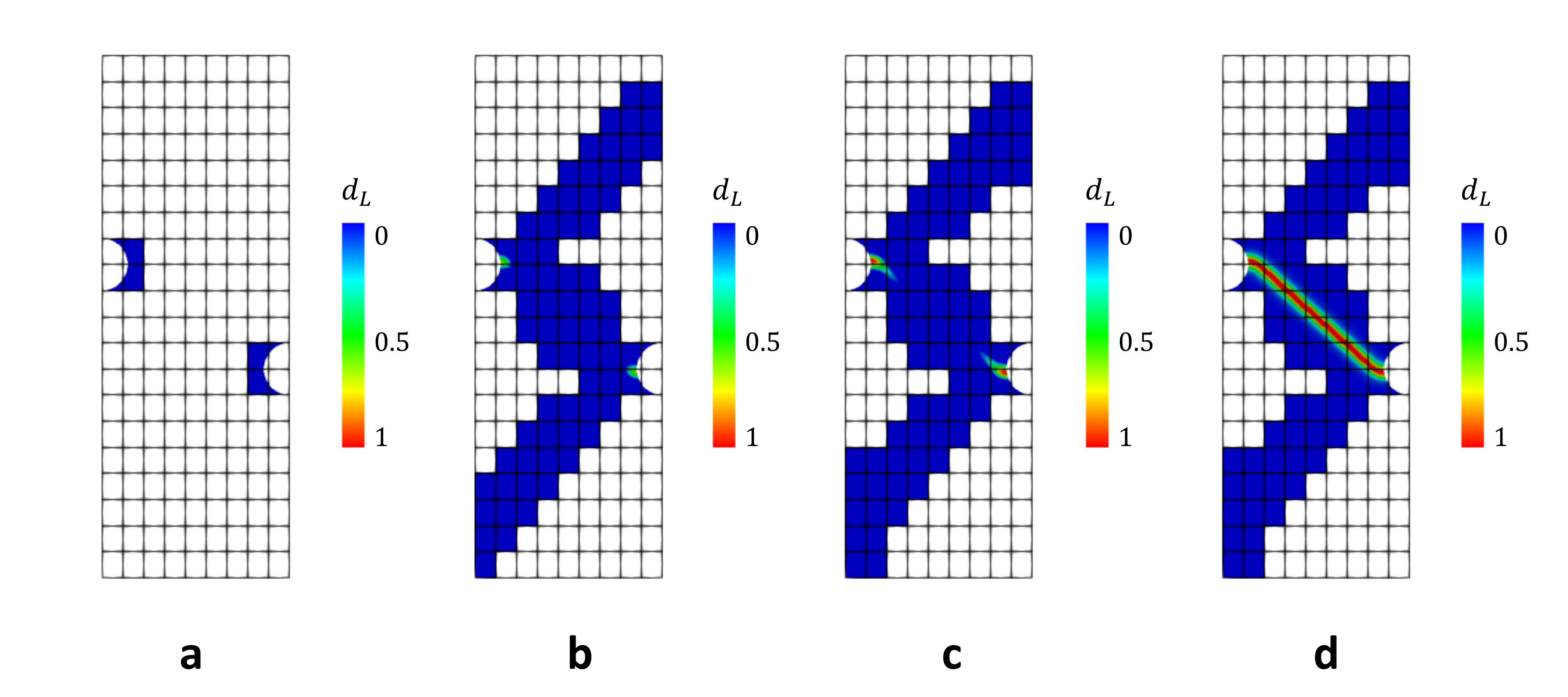}}  
	\caption{Example 1 ($g/l-1$). Evolution of the local phase-field $d_L$ for different deformation stages up to
		complete failure at $\bar{u}_y=[0,0.36,0.394,0.42]~mm$.}
	\label{exm1_profile_gl1_dL}
\end{figure}

Figure \ref{exm1_profile_gl1_hardningL} and \ref{exm1_profile_gl1_dL} illustrate the evolution response for the local solutions corresponds to $\alpha_L$ and $d_L$, respectively, for different deformation stages. As already mentioned, in $g/l-1$, the adaptivity criterium  is devised through $\eta_{e_G}$ which obtained from the global maximum equivalent plastic strain. In fact, the adaptive elements resolved here results from the global localization branches which form at an angle   about $45^\circ$. These localization bands corresponds to the shear band dictated from the given BVP. 

Accordingly, the fracture path initiates within the maximum equivalent plastic region, in which they appear near the notches. Next, the crack propagates in the plastic localization band, in which two cracks are merging at the specimen center; see Fig. \ref{exm1_profile_gl1_dL}. 

\begin{figure}[!ht]
	\centering
	{\includegraphics[clip,trim=4cm 21cm 4.5cm 5cm, width=15cm]{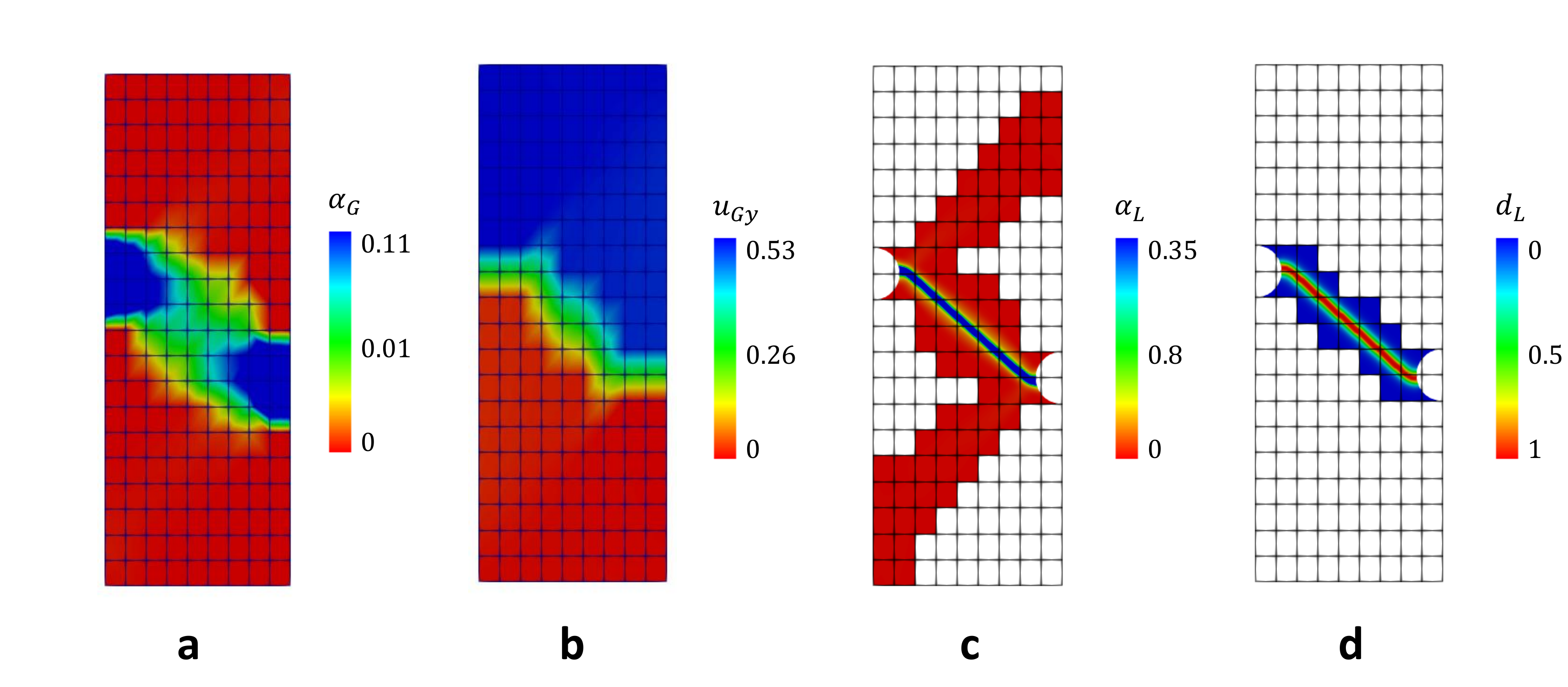}}  
	\caption{Example 1 ($g/l-2$). Approximated solution obtained through $g/l-2$ at the complete failure. (a) Global hardening value $\alpha_{G}$, (b) global vertical displacement $u_{Gy}$, (c) local hardening value $\alpha_{L}$ at $\calB_{L_1}$, and (d) local phase-field $d_{L}$ at $\calB_{L_2}$.}
	\label{exm1_profile_gl2}
\end{figure}

Now, we examine the proposed multilevel Global-Local solutions denoted as $g/l-2$. Figure \ref{exm1_profile_gl2} explains the  approximated solution for the global hardening value $\alpha_{G}$,  the global vertical displacement $u_{Gy}$,  the local hardening value $\alpha_{L}\in\calB_{L_1}$, and the local crack phase-field $d_{L}\in\calB_{L_2}$. The first important observation is that the results are in well-agreement with $g/l-1$, yet in $\calB_{L_2}$ much less elements are required to be locally refined. This leads to further reduction of computational time but preserving the numerical accuracy.

\begin{figure}[!ht]
	\centering
	{\includegraphics[clip,trim=0cm 5.5cm 0cm 1cm, width=16cm]{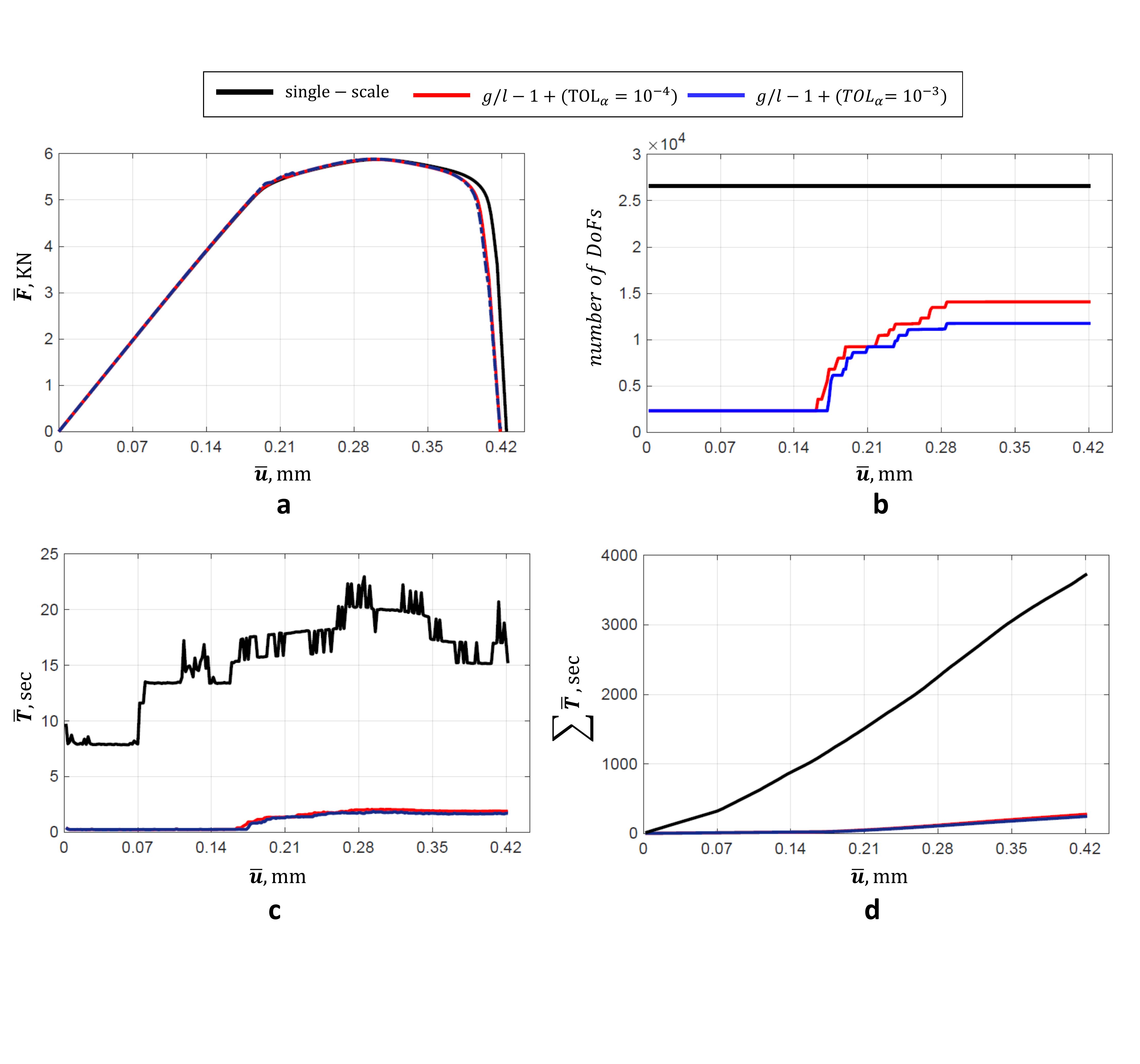}}  
	\caption{Example 1 ($g/l-1$). Computed response for double-notched specimen under tensile loading through single-scale  solution and $g/l-1$. (a) Comparison of the
		load-displacement curves, (b) number of degrees of freedom, (c) approximated time at fixed loading steps, and (d) accumulated time-displacement curves.}
	\label{exm1_curve_gl1}
\end{figure}

A comparison of the load-displacement curves for a single-scale response as well as different Global-Local schemes are shown in Fig. \ref{exm1_curve_gl1}a. Figure \ref{exm1_curve_gl1}b-d describe the efficiency of the proposed Global-Local approach. Here, the number of unknowns, approximated computational time per load steps, and total accumulated time-displacement curves are presented. It can be observed that the total accumulated time for the $g/l-1$ took 245 s whereas the single-scale simulation took 3732 s. Hence, Global-Local formulations perform 15.2 times faster; see Table \ref{GL-SS-eG1}. Accordingly, $g/l-2$ took 130 s thus it perfumes 47$\%$ faster than $g/l-1$ which turns to be 28.7 times faster than single-scale solutions; see Table \ref{GL-SS-eG1}-\ref{GL-SS-eG2}. 

\begin{figure}[!ht]
	\centering
	{\includegraphics[clip,trim=0cm 5.5cm 0cm 0cm, width=16cm]{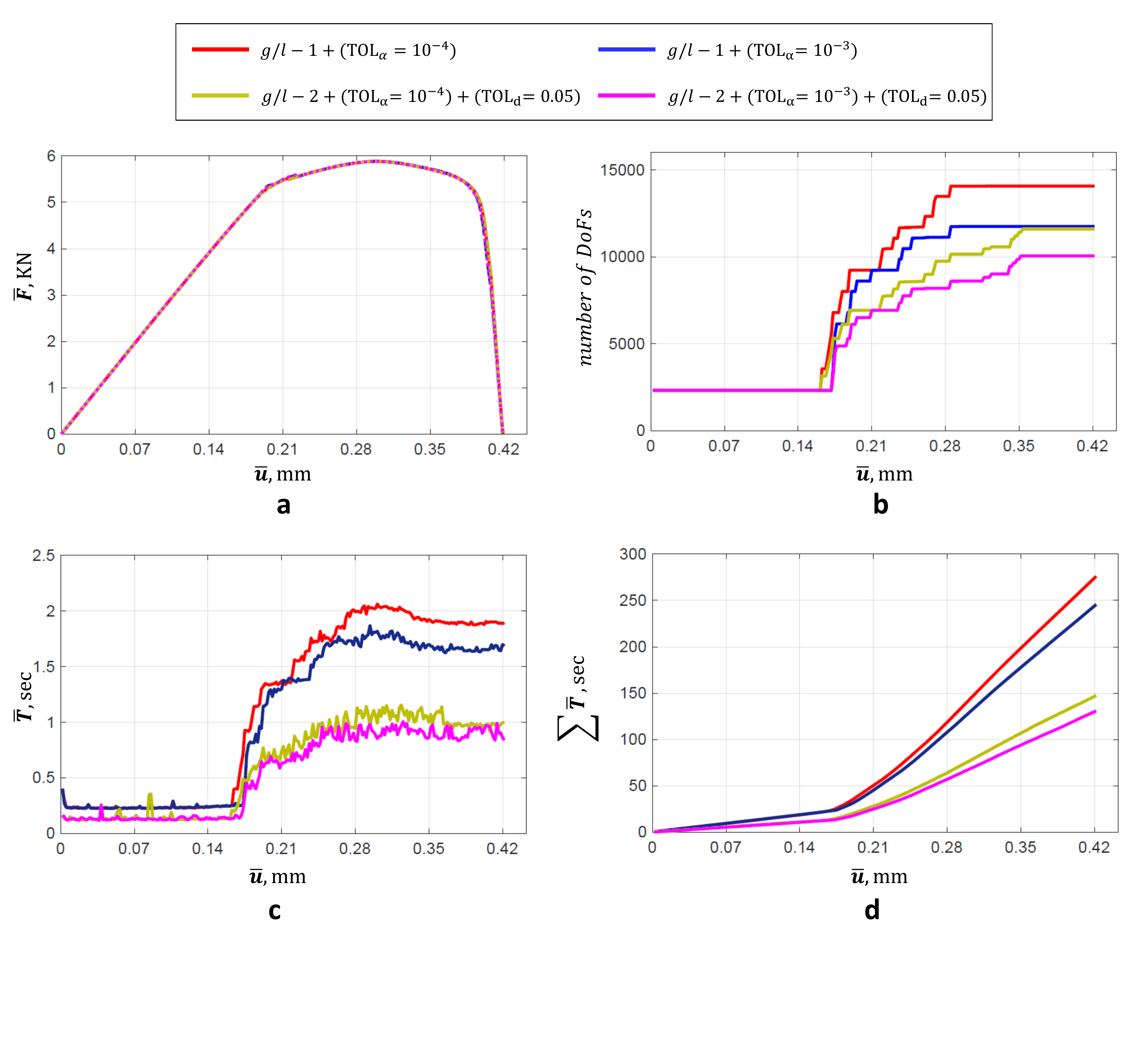}}  
	\caption{Example 1. Computed response for double-notched specimen through different global-local schemes. (a) Comparison of the
		load-displacement curves, (b) number of degrees of freedom, (c) approximated time at fixed loading steps, and (d) accumulated time-displacement curves.}
	\label{exm1_curve_gl2}
\end{figure}

Furthermore, the role played by the $\texttt{TOL}_\alpha$ which enters in $\eta_{e_G}$ in \req{eta_G} is investigated. It can be grasped that by choosing a larger value for the  $\texttt{TOL}_\alpha$, the computational time is drastically reduced; as demonstrated in Fig. \ref{exm1_curve_gl2}b-d. Yet good accuracy in comparison with a single-scale solution is preserved; see Fig. \ref{exm1_curve_gl2}a.

In conclusion, the proposed adaptivity procedure for both $g/l-1$ and $g/l-2$ while keeping the computational cost reasonably low yields an excellent agreement compared with the single-scale solution.

\sectpb[exm2]{Example 2: Isotropic single-edge-notched shear test}
The second example is concerned with a single-edge-notched shear test, abbreviated here as SENT. A boundary value problem applied to the square plate which is shown in Fig. \ref{example1-2}b. We set $A=0.5\;mm$ hence $\calB=(0,1)^2$ $mm^2$ that includes a predefined single notch from  the left edge to the body center, as depicted in Fig. \ref{example1-2}b. The predefined crack is in the $y=A$ plane and is restricted in $0\le|\mathcal{C}|\le l_0$ with $l_0=A=0.5$. 

The numerical example is performed by applying a monotonic displacement increment  ${\Delta \bar{u}}_x=2\times10^{-4}\;mm$ in a horizontal direction at the top boundary of the specimen for 400 time steps. The minimum finite element size in the single-scale and local domains is $0.05\;mm$. The
single-scale domain partition contains 30007  elements while the global domain contains 400 elements.

\begin{figure}[!ht]
	\centering
	{\includegraphics[clip,trim=1cm 17cm 3cm 20cm, width=15cm]{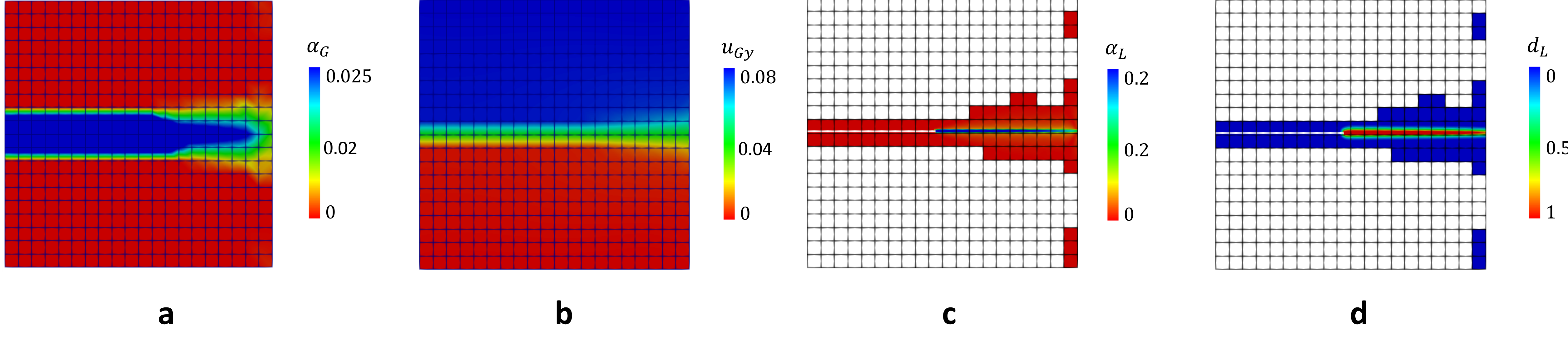}}  
	\caption{Example 2  ($g/l-1$). Approximated solution obtained through $g/l-1$ at the complete failure. (a) Global hardening value $\alpha_{G}$, (b) global vertical displacement $u_{Gy}$, (c) local hardening value $\alpha_{L}$ at $\calB_{L}$, and (d) local phase-field $d_{L}$ at $\calB_{L}$.}
	\label{exm2_profile_gl1}
\end{figure}

\begin{figure}[!ht]
	\centering
	{\includegraphics[clip,trim=0cm 1cm 1cm 12.1cm, width=15cm]{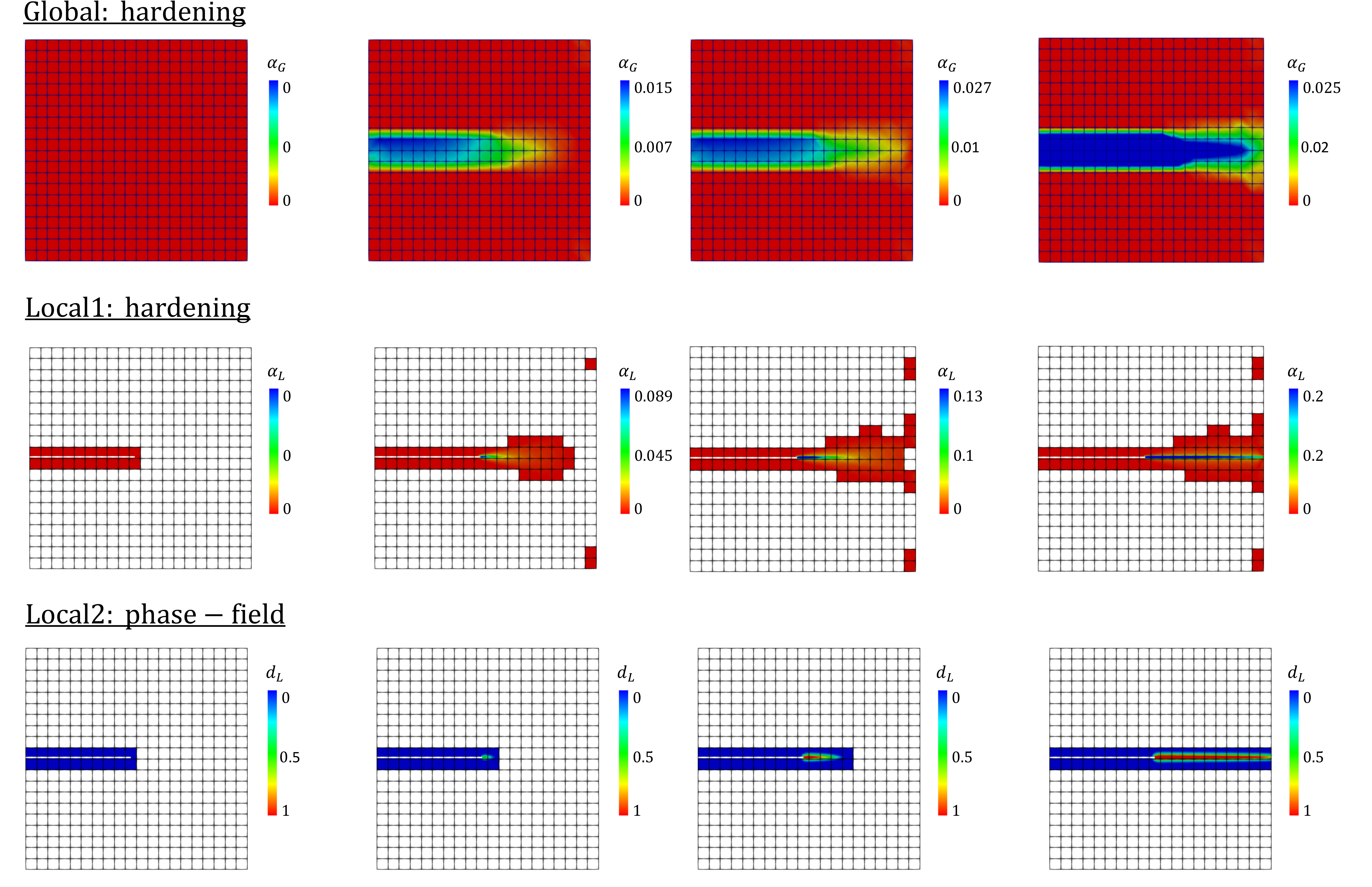}}  
	\caption{Example 2 ($g/l-2$). Approximated solution obtained through $g/l-2$ for the shear test on SENT. 
		Evolution of the (first row) global hardening value $\alpha_{G}$, (second row) local hardening value $\alpha_{L}$ at $\calB_{L_1}$, and (third row) local phase-field $d_{L}$ at $\calB_{L_2}$  for different deformation stages up to complete failure at $\bar{u}_y=[0,0.043,0.06,0.08]~mm$.}
	\label{exm2_profile_gl2}
\end{figure}

A qualitative representation of the $g/l-1$ at the complete failure is shown in Fig. \ref{exm2_profile_gl1}. Accordingly, Fig. \ref{exm2_profile_gl2} illustrates that the global equivalent plastic strain $\alpha_{G}$, local equivalent plastic strain $\alpha_{L}$, and local crack phase-field $d_{L}$  for different deformation stages up to complete failure at $\bar{u}_y=[0,0.043,0.06,0.08]~mm$. Note that, in $g/l-2$, a domain $\calB_{L_1}$ has different geometrical space compared to $\calB_{L_2}$, since it is required more elements to be refined for the elastic-plastic response while less elements are needed to capture localized crack phase-field; see Fig. \ref{exm2_profile_gl2}.
\begin{figure}[!b]
	\centering
	{\includegraphics[clip,trim=0cm 4cm 0cm 1cm, width=16cm]{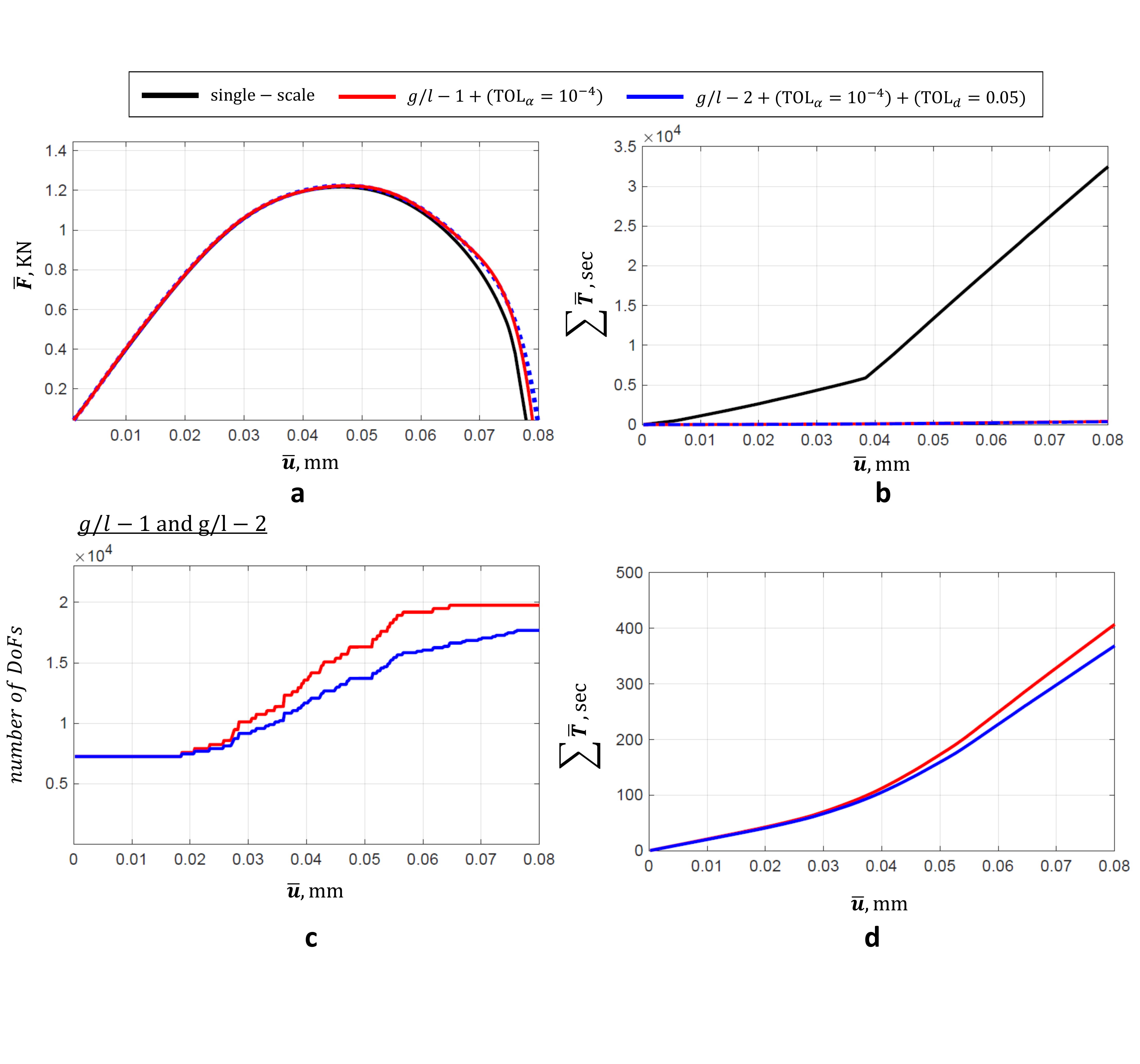}}  
	\caption{Example 2. Computed response for the SENT shear test through single-scale solution and different  Global-Local schemes. (a) Comparison of the
		load-displacement curves, (b) accumulated time-displacement curves, (c) number of degrees of freedom, and (d) accumulated time-displacement curves between $g/l-1$ and $g/l-2$.}
	\label{exm2_curve_gl1}
\end{figure}

A computed load-displacement curves shown in Fig. \ref{exm2_curve_gl1}a which demonstrates the Global-Local formulation (regardless of its type), results in an excellent agreement compared to the single-scale problem. Note that, at every jump which appears in Fig. \ref{exm2_curve_gl1}c, the predictor-corrector adaptive scheme is applied to the Global-Local scheme hence the number of degrees of
freedom is increased.

Resulting from the single-scale simulation indicates that corresponding accumulative computational time turns out to be high, whereas, the Global-Local formulation required much less  computational effort. More precisely, the total accumulated time for the $g/l-1$ took 245 s whereas the single-scale simulation took 3732 s. Hence, Global-Local formulations perform 15.2 times faster; see Table \ref{GL-SS-eG1}. Accordingly, $g/l-2$ took 130 s thus it perfumes 47$\%$ faster than $g/l-1$; see  Table \ref{GL-SS-eG2}. 

We should note  that in our setting the crack phase-field equation is linear (despite of the elastic-plastic equation). Thus, if the crack phase-field equation behaves as a non-linear equation (e.g., models used in \cite{ambati+kruse+lorenzis16,hesch+weinberg14}), using $g/l-2$ versus $g/l-1$, reduces the computational cost drastically, since for solving phase-field equation, an iterative Newton-Raphson method is also required.
\begin{figure}[!b]
	\centering
	{\includegraphics[clip,trim=0.5cm 3cm 1cm 0cm, width=17.8cm]{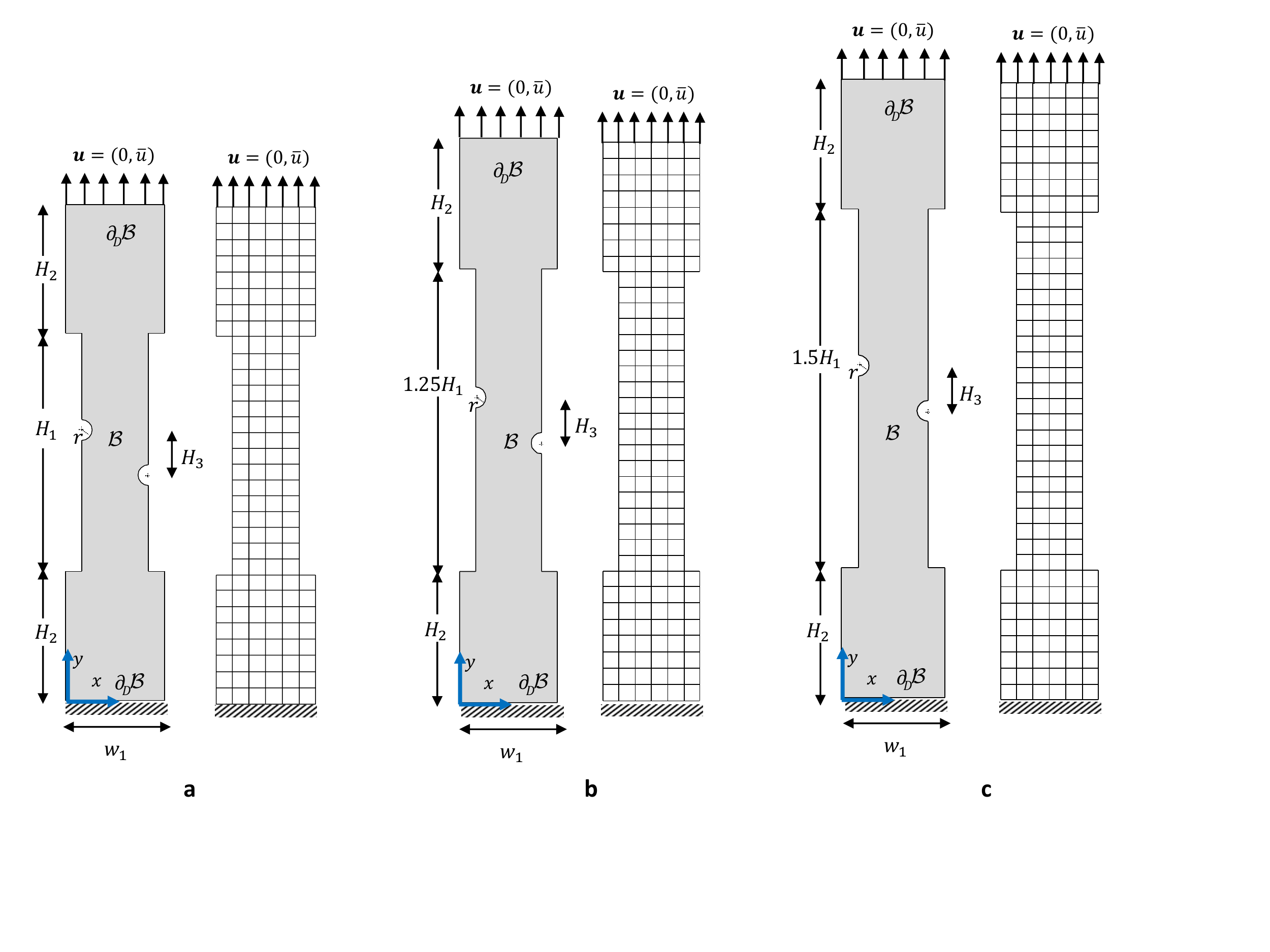}}  
	\caption{Example 3. I-shaped tensile specimen of ductile fracture for three different geometries with their boundary conditions. (a) Small , (b) medium, and (c) large size specimen with central height $H_1,\ 1.25H_1,\ 1.5H_1$, respectively, with their global discretization.
	}
	\label{example3}
\end{figure}
\sectpb[exm3]{Example 3:  I-shaped specimen under tensile loading}

The third example is aimed to illustrate the objectivity of the Global-Local formulation. 
By objectivity of the Global-Local formulation we mean, if the inelastic response is bounded in a localized region, thereafter changing the specimen size does not (typically) change  the localization band, which, in turn, the simulation time in Global- Local formulation (approximately) remain same.

To this end, an I-shaped specimen with three different sizes under tensile loading are considered. These  geometries are denoted as  small, medium, and large specimens, as shown in Fig. \ref{example3}.  The geometrical dimensions are set as $H_1=52.8\;mm$, $H_2=28.6\;mm$, $H_3=10\;mm$, and $w_1=22\;mm$ with radius of two notches as $r=2.5\;mm$.

A monotonic displacement increment  ${\Delta \bar{u}}_y=4\times10^{-4}\;mm$ is applied in a vertical direction at the top boundary of the specimens. Accordingly, we set 1200, 1350, and 1450 time steps for small, medium, and large sizes, respectively. To remove the rigid body motion, the bottom edge is fixed in $x-y$ directions. The minimum finite element size in the single-scale and local domains is $0.3\;mm$. The single-scale domain partition contains 30651 [small], 33588 [medium], and 40603 [large] elements, while the global domain contains 156 [small], 172 [medium], and 188 [large] elements correspond to different geometry size shown in Fig. \ref{example3}. 

The specific geometry used here dictates that the global equivalent plastic strain is bounded in the localized region, which, in turn, changing the geometry size will not drastically change the plastic strain localization band. Thus, the Global-Local formulation results in different global domains while having approximately the same local domain (for the elastic-plastic/fracture response). Note that, the strain localization band will be obtained through the proposed adaptivity procedure, thus it is not known in priory. The results are compared with a single-scale solution to demonstrate the powerful performance of the proposed Global-Local method. 

\begin{figure}[!ht]
	\centering
	{\includegraphics[clip,trim=0cm 0cm 0cm 0cm, width=16cm]{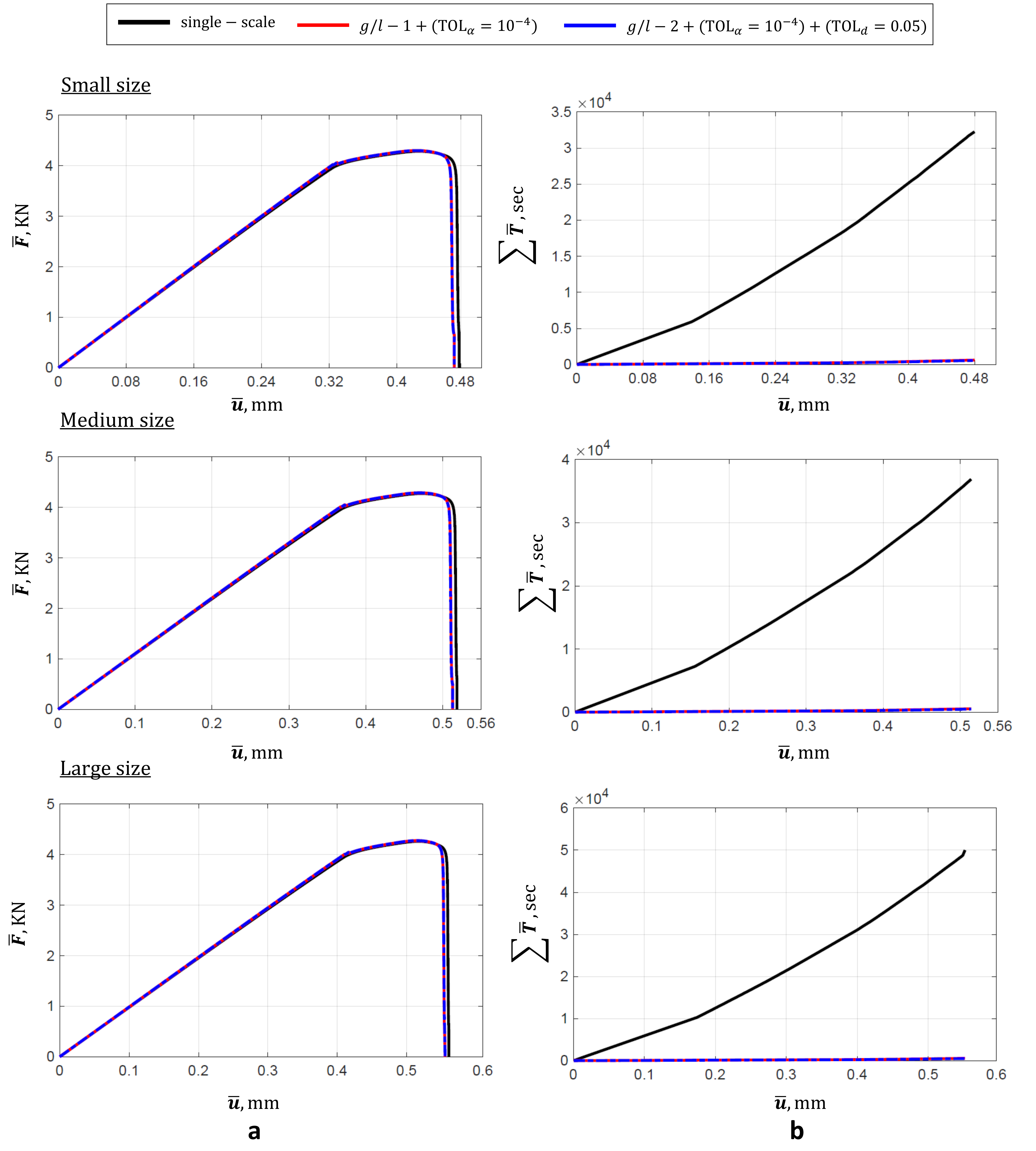}}  
	\caption{Example 3.  Computed response for different size I-shaped tensile specimen through single-scale  solution, $g/l-1$, and $g/l-2$. (a) Comparison of the
		load-displacement curves, and (b) accumulated time-displacement curves.}
	\label{exm3_curve_small_gl1}
\end{figure}

Figure \ref{exm3_curve_small_gl1} presents a comparison of the load-displacement curves for a signle-scale response as well as different Global-Local schemes. A very good agreement between the single-scale and the Global-Local solutions demonstrates the precise transition of the local non-linear constitutive mode as well as the local imperfections toward the global level; see Fig. \ref{exm3_curve_small_gl1}a. Additionally, following Fig. \ref{exm3_curve_small_gl1}b as well as Tables \ref{GL-SS-eG1}-\ref{GL-SS-eG2}, the computational time corresponds to the Global-Local formulation reduced drastically, which highlights the role of predictor-corrector adaptivity.

\begin{figure}[!t]
	\centering
	{\includegraphics[clip,trim=0cm 5.3cm 0cm 18cm, width=16cm]{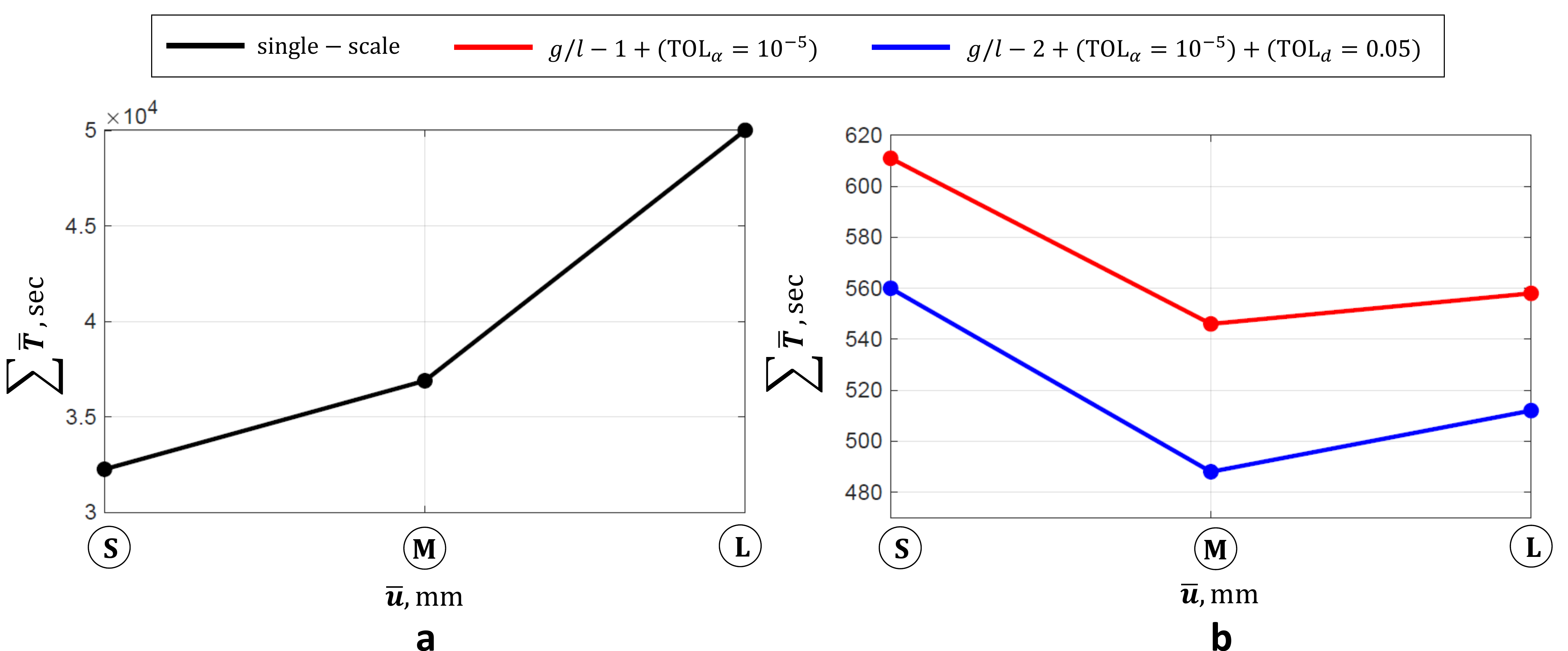}}  
	\caption{Example 3. Final simulation time at complete failure for small (S), medium (M), and large (L) I-shaped specimen. (a) Single-scale, and (b) different Global-Local schemes.}
	\label{comp_I_time}
\end{figure}

For a better insight into the computational cost and the objectivity of the Global-Local formulation, Fig. \ref{comp_I_time} is presented. The first important observation is that  extending the I-shaped domain will significantly increase the
computational cost for the single-scale problem (due to the increased elements number). Nevertheless, this does not
change the computational cost for the Global-Local formulation, thus applicable for real large structures; see Fig. \ref{comp_I_time}. In a summary, the results obtained from Global-Local formulation are 52.8, 67.6, and 89.6 times faster than single-scale simulations for small, medium, and large specimen, respectively; see Table \ref{GL-SS-eG1}-\ref{GL-SS-eG2}.

\begin{figure}[!b]
	\centering
	{\includegraphics[clip,trim=0cm 3cm 5cm 3cm, width=15cm]{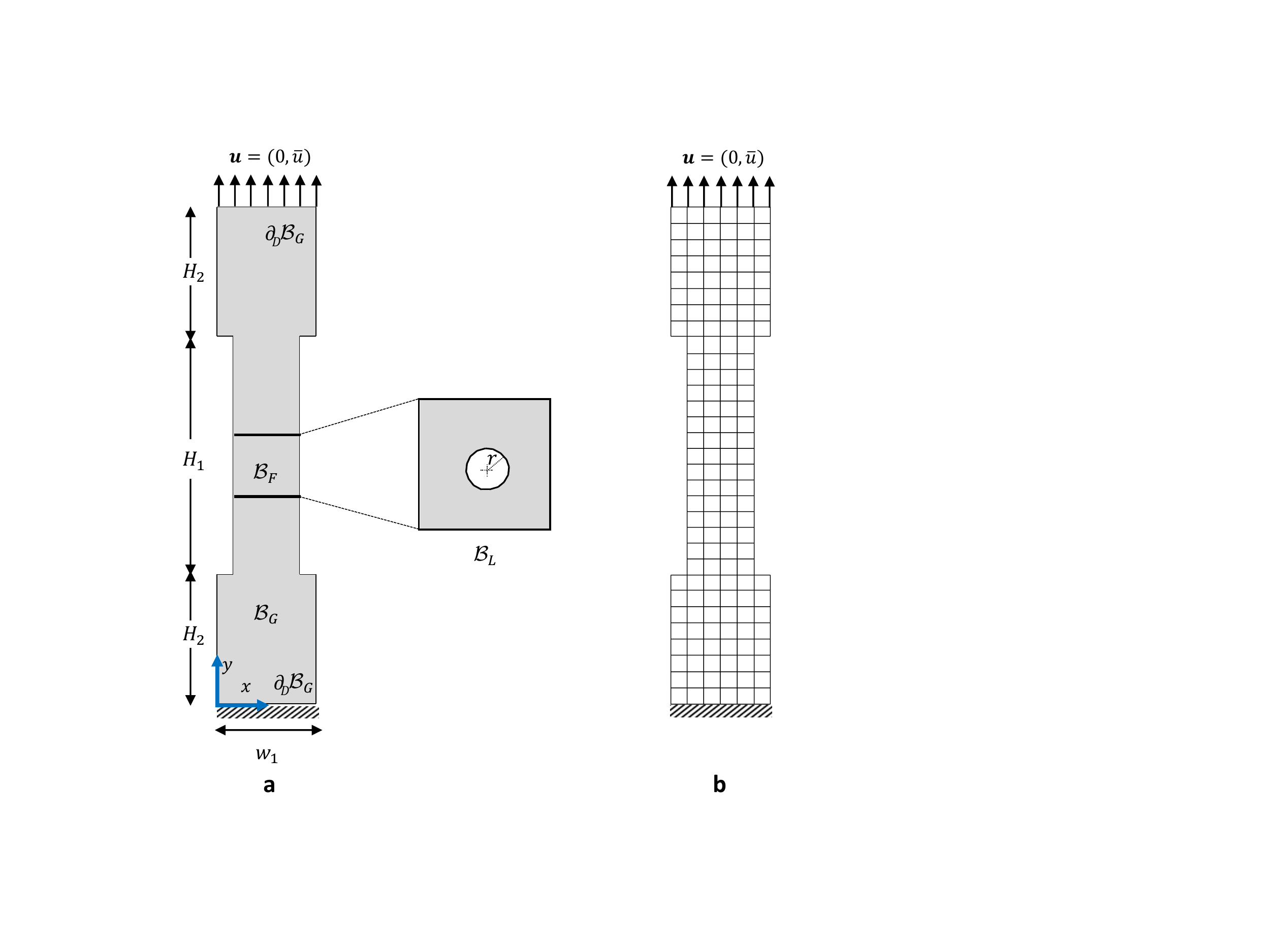}}  
	\caption{Example 4. Cyclic loading for an I-shaped tensile specimen of ductile fracture. (a) Geometry of specimen with their boundary condition and representative local imperfection, and (b) global discritization.
	}
	\label{example4}
\end{figure}
\sectpb[exm4]{Example 4:  I-shaped tensile specimen with cyclic loading}

The main objective of the final example is an adoption of the Global-Local formulation for cyclic loading applied to the ductile fracture. A boundary value problem is depicted in Fig. \ref{example4}, which is an I-shaped specimen with a circular void in the center of domain. The geometrical dimensions in Fig. \ref{example4}a are set as $H_1=52.8\;mm$, $H_2=28.6\;mm$, and $w_1=22\;mm$ with radius of void as $r=2.5\;mm$.

\begin{figure}[!t]
	\centering
	{\includegraphics[clip,trim=0cm 5.5cm 0cm 1cm, width=16cm]{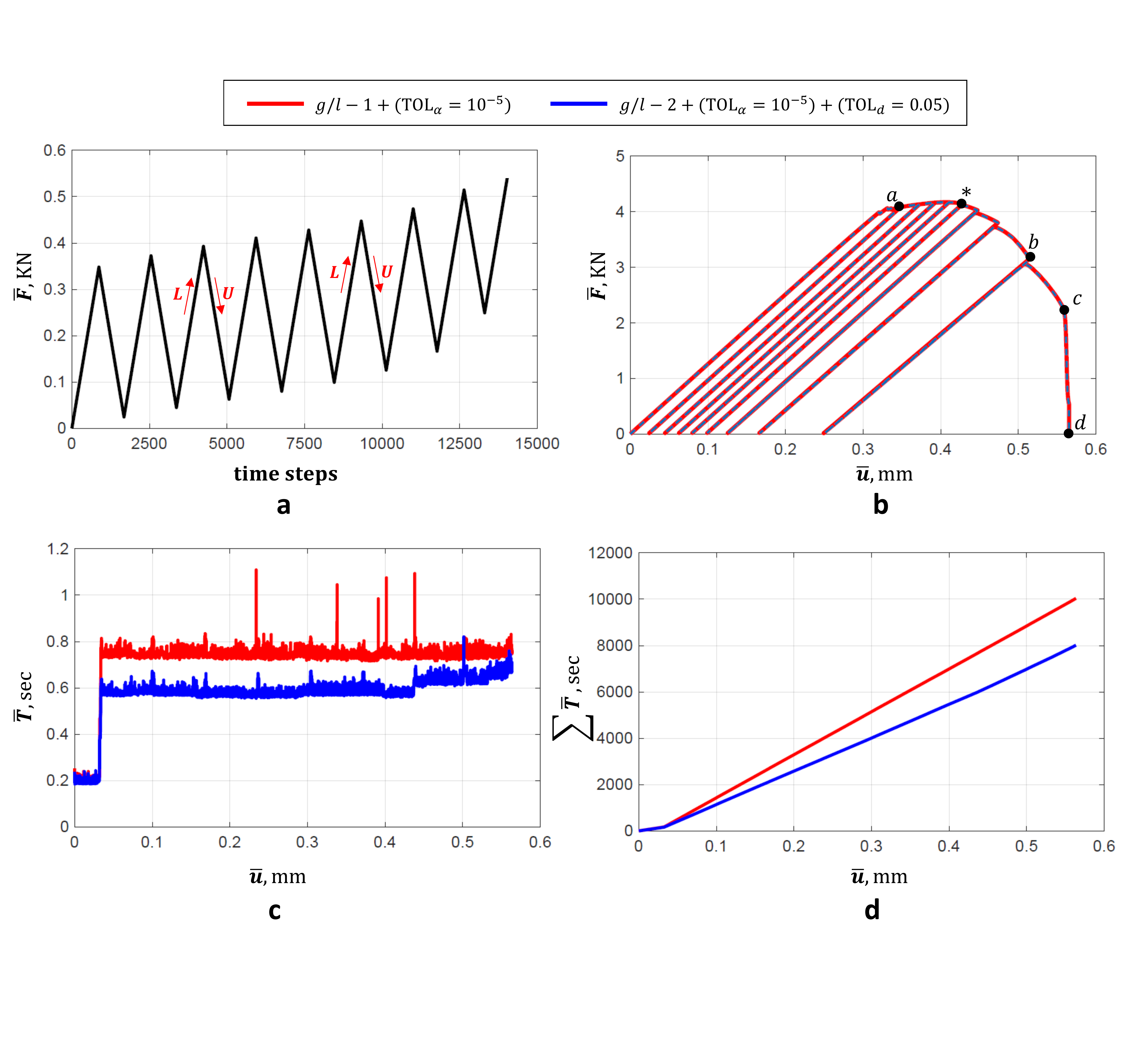}}  
	\caption{Example 4. Cyclic loading applied to the I-shaped tensile specimen with different Global-Local schemes. (a) History of load-unloading evolution applied to the specimen, (b) comparison of the
		load-displacement curves with the points of interest $(a-d)$  at $\bar{u}_y=[0.348,0.515,0.558,0.565]~mm$, (c) approximated time at fixed loading steps, and (d) accumulated time-displacement curves. Note that the point shown by ($\ast$) in  load-displacement curve stands for the onset of fracture.}
	\label{exm4_curve_gl2}
\end{figure}

The numerical example is performed by applying a loading-unloading displacement increment  $\pm{{\Delta \bar{u}}_y}=4\times10^{-4}\;mm$; see Fig. \ref{exm4_curve_gl2}a, in a vertical direction at the top boundary of the specimen for 14200 time steps. To remove the rigid body motion, the bottom edge is fixed in $x-y$ directions. The minimum finite element size in the single-scale and local domains is $0.3\;mm$. The single-scale domain partition contains 20296  elements, while the global domain contains 156 elements. 

Here, we examined the load-displacement curves,  computational time at fixed loading step, and the
total accumulated time for the $g/l-1$ and $g/l-2$. The results depicted in Fig. \ref{exm4_curve_gl2}, shows that $g/l-2$ requires 20$\%$ less computational time in comparison with $g/l-1$, which highlights the role of multilevel Global-Local formulation.  Thus, from  Fig. \ref{exm4_curve_gl2}c-d, it can be concluded that the desired improvement of
efficiency of the $g/l-1$ toward $g/l-2$ has indeed been achieved. Meanwhile, two different Global-Local formulations result in an identical representation of the load-displacement curve; see Fig. \ref{exm4_curve_gl2}b.

\begin{figure}[!b]
	\centering
	{\includegraphics[clip,trim=1cm 2cm 1cm 1cm, width=16cm]{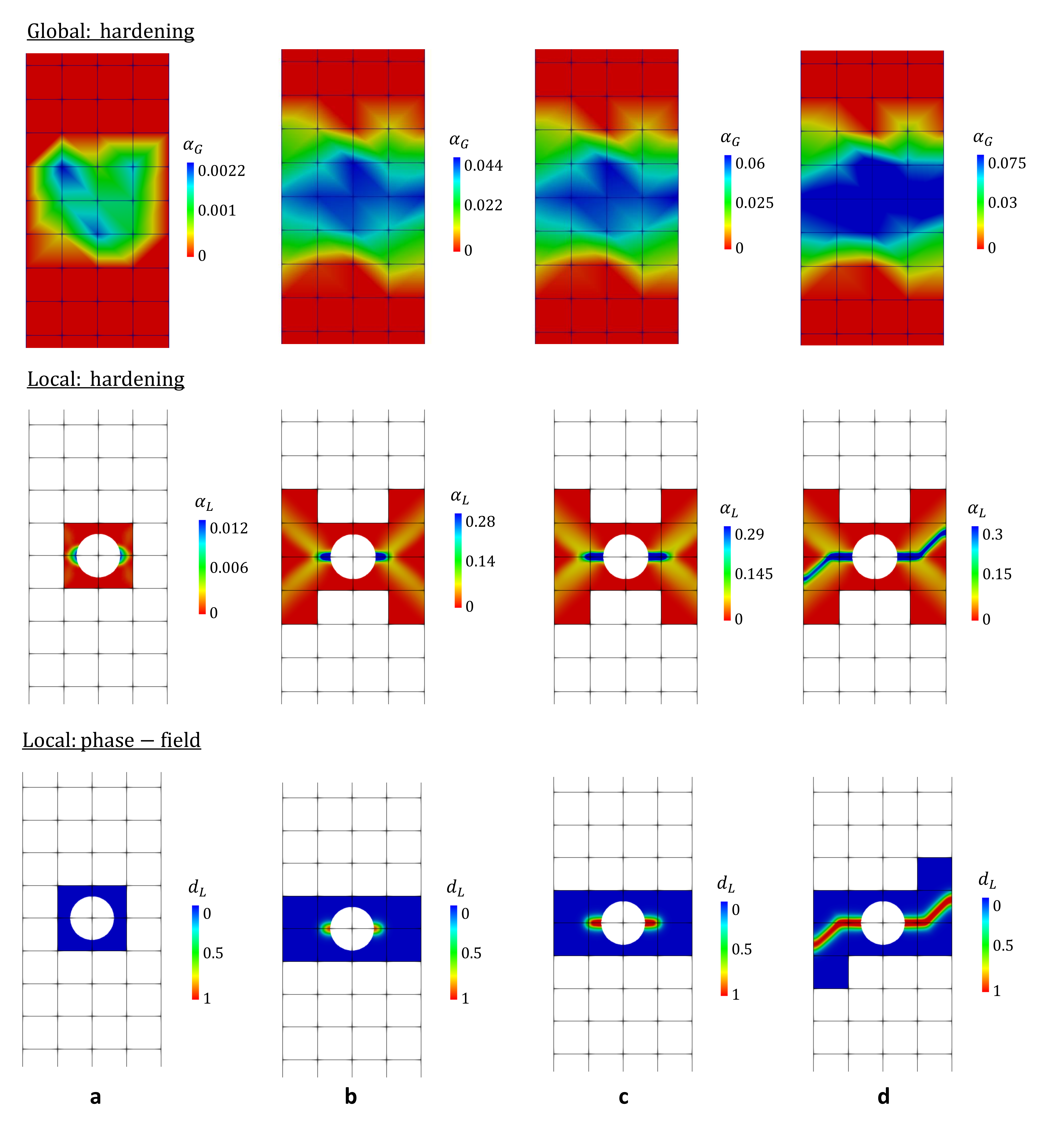}}  
	\caption{Example 4 ($g/l-2$). Evolution of the (first row) global hardening value $\alpha_G$, (second row) local hardening value $\alpha_L$, and (third row) local crack phase-field $d_L$  for different deformation stages up to complete failure (the pints $(a-d)$ from Fig. \ref{exm4_curve_gl2}b)}
	\label{exm4_profileL}
\end{figure}

The evolution of the equivalent plastic strain $\alpha_G$ and $\alpha_L$ as well as local crack phase-field $d_L$ for the $g/l-2$ are provided in Fig. \ref{exm4_profileL} at four deformation stages up to final failure. The crack initiates at the tip of the circular void at both sides and continues to propagate straight toward the edges of the specimen till the end of computation. This also holds for the local equivalent plastic strain $\alpha_L$. Another impacting factor that should be noted is that the global hardening value $\alpha_G$ as an effective quantity accurately captures the influence of local void at the upper level. Thus, the maximum global plastic flow is observed at the middle of the specimen; see Fig. \ref{exm4_profileL}(first row). Herein, the evolution of $\alpha_G$ is aligned with local hardening flow, thus underlines the precise coupling transition between two scales.

We note that, while in $g/l-1$ the local domain remains identical for both elastic-plastic behavior
and the crack phase-field, $g/l-2$, requires much less global elements to be refined for resolving the crack phase-field in $\calB_{L_2}$. More precisely, prior to point $b$ in Fig. \ref{exm4_curve_gl2}b, the local $\calB_{L_2}$ remains identical compares to the first time step in $g/l-2$, thus huge computational cost for solving the crack phase-field is avoided; see Fig. \ref{exm4_profileL}c-d.

\begin{table}[!b]\small
	\caption{Comparison accumulated time and degrees of freedom between the single-scale and $g/l-1$ formulation for different numerical examples.}
	\vspace{3mm}
	\label{GL-SS-eG1}
	\centering
	\begin{tabular}{|c|c|c|c|c|c|l}
		\hline
		Numerical examples  & \multicolumn{3}{c|}{Accumulated time, $sec.$} &  \multicolumn{2}{c|}{Total Degrees of freedom}  \\
		& single-scale & $g/l-1$  & ratio$^{\star}$ & single-scale & $g/l-1$ \\
		\hline
		Exm. 1   & 3732    & 245   & 15.2  &26577   & 11751 \\
		Exm. 2  & 32352   & 405   &79.9  &90834  &19761 \\
		Exm. 3: Small I-shaped \;  \; & 32261   & 611   & 52.8  &93126   &14705$^{\star\star}$ \\
		Exm. 3: Medium I-shaped  & 36896   & 546   & 67.6  &102087  &12804$^{\star\star}$ \\
		Exm. 3: Large I-shaped   \;  \; & 50018   &558  & 89.6 &123375  &12778$^{\star\star}$ \\
		\hline
	\end{tabular}\\[1.8mm]
	{\parbox{10in}{{${\star} ~ratio:=\frac{time_{single-scale}}{time_{g/l-1}}$~.}}}
	{\parbox{10in}{{${\star\star} ~$\textit{Note that, in medium and large specimens, an equivalent plastic strain is bounded in a\\ narrower  region compared to the small size, thus resulting in less computational cost}}}}
\end{table}
\begin{table}[!ht]\small
	\caption{Comparison accumulated time and degrees of freedom between $g/l-1$ and $g/l-2$ schemes for different numerical examples.}
	\vspace{3mm}
	\label{GL-SS-eG2}
	\centering
	\begin{tabular}{|c|c|c|c|c|c|l}
		\hline
		Numerical examples  & \multicolumn{3}{c|}{Accumulated time, $sec.$} &  \multicolumn{2}{c|}{Total Degrees of freedom}  \\
		& $g/l-1$ & $g/l-2$  & $\%$ratio$^{\star}$ & $g/l-1$ & $g/l-2$  \\
		\hline
		Exm. 1   & 245   & 130   & 47  &11751   &10058 \\
		Exm. 2  & 405   & 367   & 9.3 &19761  &17689 \\
		Exm. 3: Small I-shaped \;  \; & 611   & 560   & 8.3  &14705   &16223 \\
		Exm. 3: Medium I-shaped  & 546   & 488   & 10.6  &12804  &14124 \\
		Exm. 3: Large I-shaped   \;  \; & 558  & 512   & 8&12778  &18058 \\
		Exm. 4  & 10037  & 8013   & 20 &12086  &11384 \\
		\hline
	\end{tabular}\\[1.8mm]
	{\parbox{10in}{{${\star} ~\%ratio:=\frac{time_{g/l-1}-time_{g/l-2}}{time_{g/l-1}}\times 100$~.}}}
\end{table}
Additionally, the results obtained from Global-Local formulation demonstrate that the loading-unloading conditions (i.e, \req{kktP}, \req{kktD}, and \req{kktPG}) which are imposed to the Global-Local formulation is precisely fulfilled. Specifically, that means during unloading state, we do not have a crack phase-field evolution; i.e. prior to the point $b$ in Fig. \ref{exm4_curve_gl2}b. Meanwhile, the slope of the unloading stage is aligned with an elastic stage (if we are in the plastic phase and prior to the onset of fracture). In turn, an equivalent plastic strain remains
constant during unloading process. Thus, the KKT conditions for both plasticity and fracture response in the global and local levels are accurately performed. 

In conclusion, the developed models showed its proficiency for cyclic loading while  keeping the computational cost reasonably low.  Hence, a further application of the proposed Global-Local framework can be applied toward a fatigue failure problem, which is the subject of our ongoing research work. 

\sectpa[Section7]{Conclusion}

In this work, we outlined a robust and efficient Global-Local approach for phase-field ductile fracture problems. Hereby, a fine mesh is required to approximate the sharp crack topology  resulting in a $\textit{huge}$ computational cost. To overcome this difficulty, two different Global-Local formulations are proposed. In the first model ($g/l-1$), a global constitutive model behaves as an elastic-plastic response, while it is enhanced with a single local domain. Thereafter, we
developed the key goal of this contribution, by describing the second Global-Local formulation ($g/l-2$). The main objective of this extension was to introduce an adoption of the Global-Local approach toward the \textit{multilevel} local setting.  Because, the strain localization band is not known in a priory, an adaptivity procedure is proposed. A predictor-corrector adaptivity scheme was devised through the evolution of the
effective global plastic flow (for the first Global-Local approach), and through the evolution of both effective global  equivalent plastic strain and the local crack phase-field state (for the second Global-Local approach). 

In our numerical simulations, we have shown that the Global-Local approach has the potential to tackle practical field problems in which large structures might be considered. Additionally, it requires significantly less degrees of freedom than the single-scale formulation, leading to a remarkable reduction of the computational time. It is observed that an average accumulated simulation time for the $g/l-1$ approach was up to 60 times faster than the standard phase-field formulation (single-scale solution). Meanwhile, the $g/l-2$ approach  was up to  15$\%$ faster than the $g/l-1$ formulation, yet, an excellent performance for both Global-Local approaches of the proposed framework are observed. We also examined Kuhn-Tucker conditions (i.e., KKT for global plasticity and local fracture and plasticity) in the proposed Global-Local formulation by imposing loading-unloading conditions. Our numerical results demonstrate, the KKT conditions for both plasticity and fracture states at the global and local levels are accurately performed. In conclusion, the proposed adaptivity procedure for both Global-Local formulations while keeping the computational cost reasonably low yields an excellent agreement compared with single-scale solutions. In a future study, a further application of the proposed Global-Local framework to tackle a fatigue ductile failure problem will be considered.

\subsection*{Acknowledgment}
F. Aldakheel and N. Noii were founded by the Priority Program \texttt{DFG-SPP 2020} within its second funding phase. T. Wick and P. Wriggers were funded by the Deutsche Forschungsgemeinschaft (DFG, German Research Foundation) under Germany's Excellence Strategy within the Cluster of Excellence PhoenixD, \texttt{EXC 2122} (project number: 390833453). O. Allix  would like to  thank the Alexander Foundation for its support through the the Gay-Lussac-Humboldt prize which made it possible to closely interact with the colleagues from the Institute of Continuum Mechanics at Leibniz Universit\"at Hannover.
\cleardoublepage
\begin{Appendix}
	\setcounter{equation}{0}
	\renewcommand{\theequation}{A.\arabic{equation}}
\subsection*{Appendix A. Compact \textsc{MATLAB} Open-Source Code for Global-Local Approach: One-Dimensional Elasticity.}

This Appendix provided the compact \textsc{MATLAB} open-source  code which can be used for one-dimensional elasticity analysis through the Global-Local formulation, given in Section \ref{Section2}.
\begin{lstlisting}
%=============================
% Initialization
%=============================
L=1;A=1;          % Length//the cross-sectional area
E1=10;l1=L;       % local Young's modulus  for 0<x<L
E2=2*E1;l2=L;     % local Young's modulus  for L<x<2L 
E3=3*E1;l3=8*L;   % global Young's modulus        
tau=0.1;          % applied traction force at the global level  
TOL=1E-12;        % tolerance of iterative GL         
%-----------------------------
uR=[];uG=[];uL=[];
resi=[];    
residual=1;
i=0;
k1=(E1*A)/L;         
k2=(E2*A)/L;         
k3=(E3*A)/6*L; 
%=============================
% Reference
%=============================
KR=[k1 -k1 0 0;...
-k1 k1+k2 -k2 0;...
0 -k2 k2+k3 -k3;...
0  0 -k3 k3];                    % Jacobian
uR0=KR(2:end,2:end)\[0 0 tau]';  % reference solution displacement
uR=[0;uR0];
%=============================
% Assembeled Jacobian
%=============================
KG=[3*k3 -3*k3 0;...
-3*k3 4*k3 -k3;...
0 -k3 k3];                    % global Jacobian
KL=[k1 -k1 0;
-k1 k1+k2 -k2;
0 -k2 k2];                    % local Jacobian
KF=[3*k3 -3*k3;-3*k3 3*k3];   % fictitious Jacobian
uG0=KG(2:end,2:end)\[0;tau];  % initial global solution
uGamma=uG0(1);                % initial uGamma 
%=============================
% To perform iterative GL formulation
%=============================
while residual>TOL
i=i+1;    
%-----------------------------
% Local with Dirichlet B.C.
%-----------------------------
uL2=uGamma*(k2/(k1+k2));   % updated local solution at local node 2
uL=[0,uL2,uGamma]';        % updated local solution
rL=KL*uL;                  % updated local reaction force
lambdaL=rL(3);
%-----------------------------
% Fictitious with Dirichlet B.C.
%-----------------------------
uF=[0,uGamma]';            % updated fictitious solution
rF=KF*uF;                  % updated fictitious reaction force
lambdaF=rF(2);
rT=lambdaF-lambdaL;        % updated interface residual
%-----------------------------
% Global with Neumann B.C.
%-----------------------------
uG0=KG(2:end,2:end)\[rT;tau];
uG=[0;uG0];                % updated global solution
uGamma=uG0(1);             % updated interface solution
%-----------------------------
% GL residual indicator
%-----------------------------
residual=norm(uGamma-uR(3));
resi(i)=residual;
%-----------------------------
% Saving the results
%-----------------------------
mat_uL(:,i)=uL;
mat_uG(:,i)=uG;
end
%=============================
% Plot residual behaviour//reference, global and local displacements
%=============================
figure;plot(log10(resi),'LineWidth',2);title('Convergence behavior of the GL')
%-------------------
figure;plot([0,1,2,8],uR,'DisplayName','Reference u','LineWidth',2);hold on;
plot([0,1,2],mat_uL(:,1),'DisplayName','Local u','LineWidth',2);hold on;
plot([0,2,8],mat_uG(:,1),'DisplayName','Global u','LineWidth',2);title('initial displacement solutions')
%-------------------
figure;plot([0,1,2,8],uR,'DisplayName','Reference u','LineWidth',2);hold on;
plot([0,1,2],mat_uL(:,i),'DisplayName','Local u','LineWidth',2);hold on;
plot([0,2,8],mat_uG(:,i),'DisplayName','Global u','LineWidth',2);title('Converged displacement solutions')
%=============================
\end{lstlisting}
\end{Appendix}
\bibliographystyle{bibls1}
\bibliography{./lit.bib}

\end{document}